\pdfoutput=1
\RequirePackage{fix-cm}
\documentclass[12pt,reqno]{amsart}
\usepackage{amsmath,amssymb,amsfonts,latexsym,amsthm,enumerate,mathabx}

\usepackage[T1]{fontenc}
\usepackage{lmodern}

\usepackage{xfrac}
\usepackage{mathrsfs}

\DeclareFontFamily{U}{rsfs}{\skewchar\font127 }
\DeclareFontShape{U}{rsfs}{m}{n}{%
   <-6> rsfs5
   <6-8> rsfs7
   <8-> rsfs10
}{}

\DeclareFontFamily{U}{matha}{\skewchar\font127 }
\DeclareFontShape{U}{matha}{m}{n}{%
   <-6> matha5
   <6-8> matha7
   <8-> matha10
}{}

\DeclareFontFamily{U}{mathb}{\skewchar\font127 }
\DeclareFontShape{U}{mathb}{m}{n}{%
   <-6> mathb5
   <6-8> mathb7
   <8-> mathb10
}{}

\DeclareFontFamily{U}{mathx}{\skewchar\font127 }
\DeclareFontShape{U}{mathx}{m}{n}{%
   <-6> mathx5
   <6-8> mathx7
   <8-> mathx10
}{}

\usepackage[margin=1in]{geometry}
\usepackage{stackrel}

\usepackage{pkg3}
\usepackage{kbordermatrix}

\newcommand{\mis}{\bm{A}} 
\newcommand{\misn}{{\bar{\bm{A}}}_n} 
\newcommand{\aone}{\smash{\al_\one}} 
\newcommand{\refl}{\text{\textsf{R}}}
\newcommand{\tpm}{\text{\sc tpm}} 
\newcommand{\nodd}{\mathcal{O}} 
\newcommand{\cc}{c_0} 
\newcommand{\aubd}{\al_{\mathrm{ubd}}} 
\newcommand{\albd}{\al_{\mathrm{lbd}}} 

\newcommand{\FC}{\mathfrak{F}}
\newcommand{\fc}{\mathfrak{f}}

\newcommand{\ii}{\mathbf{i}} 

\newcommand{\jfree}{\mathbf{j}_\free}
\newcommand{\kfree}{\mathbf{k}_\free}
\newcommand{\gone}{\mathbf{g}_\one}
\newcommand{\gfree}{\mathbf{g}_\free}
\newcommand{\gfreebar}{\mathbf{\overline{g}}_\free}
\newcommand{\gzro}{\mathbf{g}_\zro}
\newcommand{\gzroum}{\smash{\mathbf{g}^\bullet}_{\hspace{-5pt}\zro}}
\newcommand{\gzz}{\mathbf{g}_\zz}

\newcommand{\gzrominval}{\smash{\mathbf{g}_\zro}^{\hspace{-4pt}\minval}}
\newcommand{\gzzminval}{\smash{\mathbf{g}_\zz}^{\hspace{-8pt}\minval\hspace{4pt}}}

\newcommand{\erho}{\ep_0}

\newcommand{\hamming}{\bm{d}}

\newcommand{\qo}{q_\one}
\newcommand{\qz}{q_\zro}
\newcommand{\qf}{q_\free}

\newcommand{\forz}{\ttt{Z}}
\newcommand{\og}{{\one\forz}}
\newcommand{\fg}{{\free\forz}}
\newcommand{\gf}{{\forz\free}}
\newcommand{\go}{{\forz\one}}
\newcommand{\fzfz}{{\forz\forz}}
\newcommand{\gz}{{\forz\zro}}
\newcommand{\zg}{{\zro\forz}}
\newcommand{\fzdot}{{\forz\cdot}}

\newcommand{\zros}{\zro^\ttt{s}}
\newcommand{\zror}{\zro^\ttt{r}}

\def\longrightharpoonup{\relbar\joinrel\rightharpoonup}
\def\longleftharpoondown{\leftharpoondown\joinrel\relbar}

\def\longrightleftharpoons{
  \mathop{
    \vcenter{
      \hbox{
	\ooalign{
	  \raise1pt\hbox{$\longrightharpoonup\joinrel$}\crcr
	  \lower1pt\hbox{$\longleftharpoondown\joinrel$}
	}
      }
    }
  }
}

\newcommand{\ffne}{{\free!\free}}

\newcommand{\bdV}{V^\pd}
\newcommand{\bdF}{F^\pd}
\newcommand{\edges}{F}

\newcommand{\bE}{\bm{E}}
\newcommand{\DD}{\mathscr{D}}

\newcommand{\extY}{\mathcal{Y}^\star}
\newcommand{\intY}{\mathcal{Y}}
\newcommand{\intR}{\mathcal{R}}
\newcommand{\intNE}{\bm{\rho}}
\newcommand{\extNE}{\bm{\ep}}
\newcommand{\intN}{\bm{r}}
\newcommand{\extN}{\bm{y}}
\newcommand{\by}{\bm{y}}

\newcommand{\numzf}{\smash{n^\bullet}_{\hspace{-5pt}\zf}}
\newcommand{\numfz}{\smash{n^\bullet}_{\hspace{-5pt}\fz}}
\newcommand{\Vumzf}{\smash{V^\bullet}_{\hspace{-5pt}\zf}}
\newcommand{\Vumfz}{\smash{V^\bullet}_{\hspace{-5pt}\fz}}
\newcommand{\pumzf}{\smash{\pi^\bullet}_{\hspace{-5pt}\zf}}
\newcommand{\pumfz}{\smash{\pi^\bullet}_{\hspace{-5pt}\fz}}

\newcommand{\numzfint}{\smash{n^\circ}_{\hspace{-5pt}\zf}}
\newcommand{\numfzint}{\smash{n^\circ}_{\hspace{-5pt}\fz}}

\newcommand{\Vucff}{\widetilde V_\ff}
\newcommand{\Vucfd}{\widetilde V_\fd}
\newcommand{\nucff}{\widetilde n_\ff}
\newcommand{\nucfd}{\widetilde n_\fd}
\newcommand{\extNucff}{\widetilde\extN_\ff}
\newcommand{\extNucfd}{\widetilde\extN_\fd}
\newcommand{\intYucff}{\widetilde\intY_\ff}
\newcommand{\intYucfd}{\widetilde\intY_\fd}

\newcommand{\Deucfd}{\widetilde\De_\fd}
\newcommand{\extNEucfd}{\widetilde\extNE_\fd}

\newcommand{\tree}{\smash{T_d}}
\newcommand{\treebip}{\smash{T_{d,2}}}
\newcommand{\etreebip}{\smash{\acute{T}_{d,2}}}

\newcommand{\nuaux}{\bm{\nu}}
\newcommand{\bemax}{\be_{\max}}

\newcommand{\ue}[1][]{\bm{E}_{#1}}

\newcommand{\alfm}{\al_\square} 
\newcommand{\alfmt}{\wt{\al}_\square} 

\newcommand{\EFC}{\bm{e}}
\newcommand{\QFC}{\mathfrak{q}}

\title{Maximum independent sets on random regular graphs}

\author[J.~Ding]{${}^{*}$Jian Ding}
\author[A.~Sly]{${}^{\dagger}$Allan Sly}
\author[N.~Sun]{${}^{\ddagger}$Nike Sun}
\thanks{Research supported by 
${}^{*}$NSF grant DMS-1313596;
${}^{\dagger}$Sloan Research Fellowship; ${}^{\ddagger}$NDSEG and NSF GRF}

\date{\today}

\begin{document}

\begin{abstract}
We determine the asymptotics of the independence number of the random $d$-regular graph for all $d \ge d_0$. It is highly concentrated, with constant-order fluctuations around $n\alpha_\star-c_\star\log n$ for explicit constants $\alpha_\star(d)$ and $c_\star(d)$. Our proof rigorously confirms the one-step replica symmetry breaking heuristics for this problem, and we believe the techniques will be more broadly applicable to the study of other combinatorial properties of random graphs.
\end{abstract}

\maketitle

\section{Introduction}

An \emph{independent set} in a graph is a subset of the vertices of which no two are neighbors. Establishing asymptotics of the maximum size of an independent set (the \emph{independence number}) on random graphs is a classical problem in probabilistic combinatorics. On the random $d$-regular graph $\Gnd$, the independence number grows linearly in the number $n$ of vertices. Upper bounds were established by Bollob\'as~\cite{MR624948} and McKay~\cite{MR890138}, and lower bounds by Frieze--Suen~\cite{MR1379227}, Frieze--{\L}uczak~\cite{MR1142268} and Wormald~\cite{MR1384372}, using a combination of techniques, including first and second moment bounds, differential equations, and switchings. The bounds are quite close, with the maximal density of occupied vertices (the \emph{independence ratio}) roughly asymptotic to $2(\log d)/d$ in the limit of large $d$ --- however, for every fixed $d$ there remains a constant-size gap in the bounds on the independence ratio. For a more complete history and discussion of many related topics see the survey of Wormald~\cite{MR1725006}.

In fact, a long-standing open problem (see \cite{AldousOpen, MR2023650}) was to determine if there even exists a limiting independence ratio. A standard martingale bound implies that the independence number has $O(n^{1/2})$ fluctuations about its mean, so an equivalent question was to prove convergence of the expected independence ratio. This conjecture was recently resolved by Bayati--Gamarnik--Tetali~\cite{MR2743259} using interpolation methods from statistical physics.  Their method is based on a sub-additivity argument which does not yield information on the limiting independence ratio or the order of fluctuations.

In this paper, we establish for all sufficiently large $d$ the asymptotic independence ratio $\al_\star = \lim_n\mis_n/n$, and determine also the lower-order logarithmic correction. Further, we prove tightness of the \emph{non-normalized} independence number, proving that the random variable is much more strongly concentrated than suggested by the classical $O(n^{1/2})$ bound:

\bThm\label{t:main}
The maximum size $\mis_n$ of an independent set in the random $d$-regular graph $\Gnd$ has constant fluctuations about
\beq\label{e:max.is.n}
\misn \equiv n\al_\star-c_\star\log n
\eeq
for $\al_\star$ and $c_\star$ explicit functions of $d$, provided $d$ exceeds an absolute constant $d_0$.
\eThm

The values $\al_\star$ and $c_\star$ are given in terms of a function $\bPhistar\equiv\bPhistar_d$ defined on the interval $(5/3)(\log d)/d\le\al\le2(\log d)/d$: the function is smoothly decreasing with a unique zero $\al_\star$, and we set $c_\star\equiv -[2\cdot\bPhistar'(\al_\star)]^{-1}$ (positive since the function is decreasing). Explicitly,
\beq\label{e:explicit.intro}
\bPhistar(\al)
=-\log[1-q(1-1/\lm)]
	-(d/2-1)\log[1-q^2(1-1/\lm)]
	-\al\log\lm
\eeq
where $q$ is determined from $\al$ by solving the equation
\[\al= q \f{1+(d/2-1)  (1-q)^{d-1}}{ 1+q-(1-q)^{d-1} }
\text{ on the interval }1.6(\log d)/d\le q\le3(\log d)/d,\]
and $\lm = q[1-(1-q)^{d-1}]/(1-q)^d$. For $\lm$ determined from $\al$ in this manner we will see that $\bPhistar'(\al)=-\log\lm$, therefore $c_\star=[2\log\lm_\star]^{-1}$ for $\lm_\star$ corresponding to $\al_\star$.

A natural question is whether the same behavior holds for regular graphs of low degree. Though it is certainly possible to determine an explicit $d_0$ from our proof, we have not done so because the calculations in the paper are already daunting, and have not been carried out with a view towards optimizing $d_0$. More importantly, our result is in line with the \emph{one-step replica symmetry breaking} (\textsc{1rsb}) prediction, which is believed to fail on low-degree graphs where physicists expect \emph{full replica symmetry breaking} \cite{BKZZ}. In the latter regime no formula is predicted even at a heuristic level.

\subsection{Replica symmetry breaking}

Ideas from statistical physics have greatly advanced our understanding of random constraint satisfaction and combinatorial optimization problems~\cite{PNAS19062007, MR1026102, MR2518205}. This deep, but for the most part non-rigorous, theory has led to a detailed picture of the geometry of the space of solutions for a broad class of such problems, including exact predictions for their satisfiability threshold. While some aspects of this rich picture have been established, including celebrated results such as Aldous's solution to the random assignment problem~\cite{MR1839499} and Talagrand's proof of Parisi's formula for the Sherrington--Kirkpatrick spin-glass model~\cite{MR2195134}, many of the most important ideas remain at the level of conjecture.  We believe that recent developments, including our own previous work~\cite{DMS,DMSS}, make it possible to establish thresholds predicted by this theory for many such models.

The natural approach to studying the independence ratio is the (first and second) moment method applied to the number of $Z_{n\al}$ of independent sets of fixed density $\al$. Indeed, an analogous approach correctly determines the asymptotic independence number for the \emph{dense} Erd\doubleacute{o}s-R\'{e}nyi random graph~\cite{MR0369129}. On \emph{sparse} random graph ensembles, however, the second moment approach fails to locate the sharp transition. Due to the sparsity of the graph, almost every independent set can be locally perturbed in a linear number of places: thus the existence of a single independent set implies the existence of a \emph{cluster} of exponentially many independent sets, all related by sequences of local perturbations. Moreover the expected cluster size remains exponentially large even beyond the first moment threshold --- thus there is a regime below the first moment threshold where it overcomes the first moment, causing the second moment to be exponentially large compared with the first moment squared.

From statistical physics, the (mostly heuristic) understanding of this phenomenon is that as $\al$ exceeds roughly $(\log d)/d$, the solution space of independent sets becomes \emph{shattered} into exponentially many well-separated clusters~\cite{Coja-Oghlan:2011:ISR:2133036.2133048}. This geometry persists up to a further (conjectured) structural transition where the solution space \emph{condensates} onto the largest clusters. In the non-trivial regime between the condensation and satisfiability transitions, most independent sets are concentrated within a bounded number of clusters according to the theory from statistical physics. This within-cluster correlation then dominates the moment calculation, causing the failure of the second moment method.

In this paper, we determine the exact threshold by a novel approach which rigorizes the \textsc{1rsb} heuristic from statistical physics, which suggests that we count clusters of independent sets rather than the sets themselves. Our proof has several new ideas which we now describe.

The first is that most clusters of maximal (or locally maximal) independent sets can be given the following simple combinatorial description. A cluster is encoded by a configuration of what we call the ``frozen model'', which takes values in $\set{\zro,\one,\free}^V$ where the $\free$ (``free'') spins indicate the local perturbations within the cluster, and come in matched pairs corresponding to swapping occupied/unoccupied states across an edge. Each $\one$ indicates an occupied vertex, and can only neighbor $\zro$ (unoccupied) vertices; and each $\zro$ must neighbor at least two $\one$'s. We shall prove that this effectively encodes the clusters in this model. The expression $\bPhistar$ in equation~\eqref{e:explicit.intro} is exactly the free energy of this model, and the asymptotic density $\al_\star$ corresponds to the largest density for which it is positive.

Secondly, we note that this new model is itself a Gibbs measure on a random hypergraph, and its properties of local rigidity hint that applying the second moment in this model \emph{does} locate the exact threshold. However, the actual moment calculation appears at first intractable, involving maximizations over high-dimensional simplices. By a certain ``Bethe variational principle''~\cite{DMS,DMSS} we are able to characterize local maximizers via fixed points of certain tree recursions, reducing the optimization to (in the second moment) $81$ real variables. With delicate {\it a priori} estimates we are able to establish symmetry relations among these variables which drastically reduce the dimensionality and allows us finally to pinpoint the global maximizers.

The second moment method itself only establishes the existence of clusters with asymptotically positive probability. Our final innovation is a method to improve positive probability bounds to high probability, which in this model yields the constant fluctuations. The approach is based on controlling the incremental fluctuations of the Doob martingale of a certain log-transform of the partition function.

As an illustration of the robustness of these methods, in a companion paper~\cite{dss-naesat} we apply the same techniques to establish the exact satisfiability threshold for the random regular not-all-equal-\textsc{sat} problem. This gives the first threshold for a sparse constraint satisfaction model with replica symmetry breaking. We expect ultimately that these methods may be extended to other combinatorial properties such as the chromatic number or maximum cut, and to the sparse Erd\doubleacute{o}s-R\'{e}nyi random graphs.

\subsection{Notation}

Throughout this paper, a \emph{graph} $G=(V,M,H)$ consists of a vertex set $V$, a collection $H$ of labelled half-edges (with each half-edge incident to a specified vertex), and a perfect matching $M$ on $H$. Omitting the labelling of half-edges gives the undirected graph $G=(V,E)$ where $(x,y)\in E$ if and only if $(e_x,e_y)\in M$ for half-edges $e_x$ and $e_y$ incident to $x$ and $y$ respectively. A subset $S\subseteq V$ is an \emph{independent set} of $G=(V,E)$ if $(x,y)\notin E$ for all $x,y\in S$.\footnote{If $E$ contains a self-loop $(x,x)$ then $x$ cannot belong to an independent set.} We equivalently view an independent set as a configuration in $\{\zro,\one\}^V$ where $\one$ indicates a vertex occupied by the set.

We work with the $d$-regular configuration model: a \emph{$d$-regular graph} with vertex set $[n]\equiv\set{1,\ldots,n}$ is a perfect matching on the set $[nd]$ of labelled half-edges, where half-edge $i$ is incident to vertex $\ceil{i/d}$. Assume $nd$ is even; the number of $d$-regular graphs on $[n]$ is the double factorial
$$(nd-1)!!=\smp{j=0}{nd/2} (nd-2j+1)
=\smf{(nd)!}{(nd/2)! 2^{nd/2}}
=e^{O[1/(nd)]}\sqrt{2}(nd/e)^{nd/2}.$$
Under the configuration model, the \emph{random $d$-regular graph} $\Gnd$ corresponds to the uniformly random perfect matching on $[nd]$.

A graph is \emph{simple} if it has no self-loops  or multi-edges. For each simple $d$-regular graph with unlabelled half-edges, there are $(d!)^n$ $d$-regular graphs with labelled half-edges, so the graph $\Gnd$ conditioned to be simple is the \emph{(uniformly) random simple $d$-regular graph} $\Gndsimp$. It is well known that the probability that $\Gnd$ is simple is bounded below as $n\to\infty$ by a constant (depending on $d$), so Thm.~\ref{t:main} also holds with $\Gndsimp$ in place of $\Gnd$.

\subsection*{Acknowledgements} We are grateful to Sourav Chatterjee, Amir Dembo, Persi Diaconis, Elchanan Mossel, and Andrea Montanari for helpful conversations.

\section{Independent sets and coarsening algorithm}\label{s:coarsen}

For non-negative functions $f(d,n)$ and $g(d,n)$ we use any of the equivalent notations $f=O_d(g)$, $g=\Om_d(f)$, $f\lesssim_d g$, $g \gtrsim_d f$ to indicate $f\le C(d)\, g$ for a finite constant $C(d)$ depending on $d$ but not on $n$. We drop the subscript $d$ to indicate that we can take the same constant $C(d)\equiv C$ for all $d\ge d_0$.

\subsection{First moment of independent sets}

Let $Z_{n\al}$ count the number of independent sets of cardinality $n\al$ on graph $G$.

\blem\label{l:is.first.moment}
$\E Z_{n\al}\asymp_d n^{-1/2}\exp\{n\Phi_d(\al)\}$ for $\Phi_d$ an explicit function not depending on $n$. The first moment threshold $\alfm\equiv\alfm(d)\equiv\inf\set{\al>0:\Phi_d(\al)<0}$ satisfies
\beq\label{e:alpha0}
\alfm
=(2/d) \cdot[ \log d-\log\log d+\log(e/2)+O(\tf{\log\log d}{\log d}) ].
\eeq
Further
\beq\label{e:first.moment.2}
\f{\E Z_{n(\alfm-\de)}}
{ \exp\{ n\de[\log d-\log\log d] \} }
\begin{cases}
\le e^{O(n\de)}
	& \text{for all }0\le\de\le\alfm;\\
\ge e^{ O(n\de) }
	& \text{for }0\le \de\lesssim d^{-1}.
\end{cases}
\eeq

\bpf
The expected number of independent sets of size $n\al$ on $\Gnd$ is
\[
\E Z_{n\al}
=\smb{n}{n\al}\smp{i=0}{nd\al-1}\smf{nd(1-\al)-i}{nd-1-2i}
=\smb{n}{n\al}
	\smf{ ( nd(1-\al) )_{nd\al} }
	{ \fdf{nd}{nd\al} 
	},
\]
where $(A)_b$ and $\fdf{A}{b}$ denote respectively the falling factorial
and falling double factorial:
\beq\label{e:falling.factorial}
(A)_b\equiv
	\smp{i=0}{b-1}(A-i),\quad
\fdf{A}{b}\equiv
	\smp{i=0}{b-1}(A-1-2i).
\eeq
We assume $\al\le\slf12$ since otherwise $\E Z_{n\al}=0$. For $\al$ constant in $n$, Stirling's formula gives $\E Z_{n\al}\asymp_d n^{-1/2}\exp\{n\Phi(\al)\}$ where
\beq\label{e:first.moment}
\begin{array}{rl}
\Phi(\al)\equiv\Phi_d(\al)
\hspace{-6pt}&\equiv
	H(\al) - d[\tf12(1-2\al)\log(1-2\al) - (1-\al)\log(1-\al) ]\\
&=\al[F(\al) + O(d\al^2)]
\quad\text{with }F(\al)\equiv\log(e/\al) - (d+1) (\al/2).
\end{array}
\eeq
It is straightforward to check that $\Phi'(\al)>0$ for small $\al$ and $\Phi''(\al)<0$ for all $\al<\tf12$, so $\Phi$ has a unique zero-crossing
$\alfm$.
The function $F$ is decreasing in $\al$ with unique zero
$$\alfmt = \smf{2}{(d+1)}  W\Big( \smf{(d+1) e}{2} \Big)
=[2+o_d(1)]\,\smf{\log d}{d}$$
where $W$ denotes the principal branch of the Lambert $W$ function
defined by $z=W(z)e^{W(z)}$
(see \cite{MR1414285,MR1809988} and references therein).
Near $z=\infty$, the function $W$ has absolutely convergent series expansion
$$
W(z)
= \log z-\log\log z
+\overbrace{
	\sum_{k\ge0} 
	\smf{1}{(-\log z)^k}
	\sum_{m\ge1} \smf{1}{m!}
	\Big[\text{\footnotesize$\DS\begin{matrix}
		k+m\\k+1\end{matrix}$}\Big]
	\Big(\smf{\log\log z}{\log z}\Big)^m
	}^{ O( (\log\log z)/(\log z) ) },$$
where $[\spint{n}{m}]$ are the Stirling cycle numbers (or unsigned Stirling numbers of the first kind), generated by $\log(1+z)^m = m!\sum_{n\ge0} (-1)^{n+m} [\spint{n}{m}] z^n/n!$. By estimating $F'$ and $F''$ near $\alfmt$ we see that $\alfm=\alfmt+O(d^{-2}\log d)$; 
\eqref{e:alpha0} then follows from the above estimate on $W(z)$. The bounds \eqref{e:first.moment.2} are easily obtained by estimating $\Phi'$ near $\alfm$ and recalling that $\Phi''<0$.
\epf
\elem

\subsection{Coarsening algorithm and frozen model}

We hereafter encode an independent set as a configuration $\ux\in\set{\zro,\one}^V$ with $\one$ indicating a vertex which is occupied by the independent set. We define the following algorithm to map an independent set $\ux\in\{\zro,\one\}^V$ to a \emph{coarsened} configuration $\ueta\equiv\ueta(\ux)\in\set{\zro,\one,\free}^V$. In the coarsened model the spin $\free$ indicates vertices which are ``free,'' as follows:

\medskip\noindent\textbf{Coarsening algorithm.} Set $\ueta^0\equiv\ux$.
\bnm[1.]
\item \emph{Step 1 (iterate for $0\le s<t$): formation of free pairs.}\\
If there exists $v\in V$ such that $\eta^s_v=\zro$ and $v$ has a unique neighbor $u$ with $\ueta^s_u=\one$, then take the first\footnote{First with respect to the ordering on $V=[n]$.} such $v$, set $\eta^{s+1}_v=\eta^{s+1}_u=\free$, and declare $(uv)$ to be a matched edge. Set $\eta^{s+1}_w=\eta^s_w$ for all $w\ne u,v$.\\
Iterate until the first time $t$ that no such vertex $v$ remains.
\item \emph{Step 2: formation of single frees.}\\
Set $\eta^{t+1}_v=\free$ whenever $\eta^t_v=\zro$ and $\set{u\in\pd v:\eta^t_u=\one}=\emptyset$, otherwise set $\eta^{t+1}_v=\eta^t_v$.
\enm
Denote the terminal configuration $\ueta\equiv\ueta(\ux)\equiv\ueta^{t+1}$. Write $\match\subseteq E$ for the set of all matched free pairs formed during Step 1 of the coarsening process.

The idea is that the pre-image of any $\ueta$ under the coarsening algorithm constitutes a \emph{cluster} of independent set configurations --- a set of configurations connected by (sequences~of) local changes, in our setting by making neighboring $\zro$/$\one$ swaps. An important property of a coarsened configuration is that every $\zro$-vertex has at least two $\one$-neighbors: as this is a rigid local configuration (a $\zro$/$\one$ swap across an edge cannot be made without violating the hard-core constraint), we have some indication (non-rigorously) that different clusters will be well separated in some sense. Note that Step 2 of the coarsening algorithm is needed to ensure this property even when the initial configuration is a maximum independent set: consider for example $\zro$ --- $\one$ --- $\zro$ --- $\one$ --- $\zro$ arranged in a $5$-cycle such that every neighbor not on the cycle is a $\zro$-vertex with many $\one$-neighbors. This can be part of a maximal configuration, but Step 1 results with $\free$ $\stackrel{m}{\text{---}}$ $\free$ --- $\free$ $\stackrel{m}{\text{---}}$ $\free$ --- $\zro$ on the cycle where $m$ indicates a matched pair and the last $\zro$ has no $\one$-neighbors. The purpose of this section is show that we may discard these odd-cycle scenarios and still recover the sharp asymptotics for $\mis_n$.

Let $Z_{n\al,\ge n\be}$ denote the contribution to $Z_{n\al}$ from independent set configurations $\ux\in\{\zro,\one\}^V$ such that the coarsened configuration $\ueta(\ux)\in\{\zro,\one,\free\}^V$ has more than $n\be$ free variables.

\bppn\label{p:coarsening}
There exists an absolute constant $C$ (not depending on $d$) such that if $\al=y\logdbyd$ for $1\le y\le 2$, then $\E Z_{n\al,\ge n\be} \le C^{-n\be/2} \E Z_{n\al}$ for $\be=d^{-y} (C \log d)$. In particular if $\al=\alfm-O(d^{-3/2})$ then $\E Z_{n\al,\ge n\be}\le e^{-nc_d}$ for $\be=C(\log d)^3/d^2$.

\bpf
Suppose $\ux\in\set{\zro,\one}^V$ such that $V_\one\equiv\set{v:x_v=\one}$ has size $n\al$, and write $V_\zro\equiv V\setminus V_\one$. By definition of the coarsening algorithm, any subset $S$ of the $\free$-vertices must satisfy $|\pd S\cap V_\one|\le |S|$. Conditioned on $\ux$ being an independent set of cardinality $n\al$, the random graph $G$ can be obtained by first matching all half-edges originating from $V_\one$ to a subset of half-edges originating from $V_\zro$ uniformly at random, then placing a uniformly random matching on the remaining half-edges of $V_\zro$. For $v\in V_\zro$ let $D_v$ denote the number of neighbors of $v$ in $V_\one$; then $(D_v)_{v\in V_\zro}$ is distributed as a vector of i.i.d.\ $\Bin(d,\slf{\al}{(1-\al)} )$ random variables conditioned on $\sum_{v\in V_\zro}D_v=nd\al$. Fixing $S$ any subset of $V_\zro$ of size $n\be$, and recalling that $\al=y\logdbyd$ with $1\le y\le 2$, we have
\[\begin{array}{rl}
(\E Z_{n\al})^{-1}\E Z_{n\al,\ge n\be}
\hspace{-6pt}&\le \tbinom{n(1-\al)}{n\be}
	\P(
	\sum_{v\in S} D_v\le n\be
	\giv
	\sum_{v\in V_\zro} D_v=nd\al
	)
= n^{O(1)}\exp\{ nf_d(\al,\be)\},\vspace{2pt}\\
f_d(\al,\be)
\hspace{-6pt}&\equiv
	(1-\al) H(\tf{\be}{1-\al})
	-d\be \,
	H(\tf{1}{d} \,|\,\tf{\al}{1-\al})
	=\be [ \log\tf{y\log d}{\be d^y} + O(1)]
\end{array}\]
with the $O(1)$ uniform in $d$. Taking $\be= d^{-y} (C \log d)$ for $C$ a large absolute constant then proves the first claim, $\E Z_{n\al,\ge n\be} \le C^{-n\be/2} \E Z_{n\al}$. The second claim follows simply by noting that if $\al=\alfm-O(d^{-3/2})$ then by Lem.~\ref{l:is.first.moment} we have $\al=y\logdbyd$ with $y\log d=2\log d-2\log\log d+O(1)$ with the $O(1)$ uniform in $d$.
\epf
\eppn

\bdfn\label{d:frozen}
We say $\ueta \in \{\zro,\one,\free\}^V$ is a valid \emph{unweighted frozen model} configuration on $G$ if the following hold:
\bnm[(a)]
\item \label{d:frozen.a}
Every $v\in V$ satisfies
\bnm[(i)]
\item If $\eta_v=\one$ then $v$ neighbors only $\zro$-vertices;
\item If $\eta_v=\zro$ then at least two edges join $v$ to $\one$-vertices;
\item If $\eta_v = \free$ then $v$
	neighbors only $\zro$- or $\free$-vertices, with at least one
	neighboring $\free$;
\enm
\item In the subgraph $\FC(\ueta)\subseteq G$ induced by the $\free$-vertices, every connected component either is not a tree, or is a tree with a (necessarily unique) perfect matching.\footnote{To be precise, since we regard $G$ as a matching on the set $[nd]$ of labelled half-edges incident to vertices, $\FC$ shall be regarded as a matching on a subset of $[nd]$.}
\enm
\edfn

\bdfn\label{d:frozen.matched}
A \emph{(weighted) frozen model configuration} $\uz\equiv(\ueta,\match)$ on $G$ is a unweighted frozen model configuration $\ueta$ together with a perfect matching $\match$ on $\FC(\ueta)$: equivalently, every $v\in V$ satisfies
\bnm[(i)]
\item \label{d:frozen.matched.i}
If $\eta_v=\one$ then $v$ neighbors only $\zro$-vertices;
\item \label{d:frozen.matched.ii}
	If $\eta_v=\zro$ then at least two edges join $v$ to $\one$-vertices;
\item \label{d:frozen.matched.iii}
	If $\eta_v = \free$ then $v$
	neighbors only $\zro$- or $\free$-vertices,
	and is matched to an $\free$-neighbor by $\match$.
\enm
The \emph{intensity} of $\uz$ is the number of $\one$-vertices plus the number of matched $\free$-pairs:
\beq\label{e:intens}
\ii(\uz)\equiv|\set{v:\eta_v=\one}|
	+\tf12|\set{v:\eta_v=\free}|.
\eeq
\edfn

Let $\tZ_{n\aone,n\be}$ denote the \emph{unweighted} frozen model partition function on $G$ restricted to configurations with exactly $n\aone$ $\one$-vertices and $n\be$ $\free$-vertices; define similarly $\ZZ_{n\aone,n\be}$ with respect to the \emph{weighted} frozen model. In view of Lem.~\ref{l:is.first.moment} and Propn.~\ref{p:coarsening}, we restrict all consideration to the \emph{truncated} frozen model partition function
\beq\label{e:truncated.frozen.model}
\ZZ_{n\al}
\equiv
\sum_{\substack{
	\aone,\be\,:\, \be\le \bemax\\
	2\aone+\be=2\al
	}} \ZZ_{n\aone,n\be},\quad
	\bemax\equiv d^{-3/2}.
\eeq
We shall always assume that the normalized intensity $\al$ lies in a restricted regime:
\[\albd\equiv \smf{(5/3)\log d}{d}\le\al\le \smf{2\log d}{d}
	\equiv \aubd,\quad
\ZZ_{\ge n\al}\equiv\sum_{n\al\le n\al'\le n\aubd}\ZZ_{n\al'}.\]
We now describe our reduction from the independent set model to the frozen model.

\bthm\label{t:frozen.first.moment}
For $d\ge d_0$, it holds uniformly over $\albd\le\al\le\aubd$ that
\[\lim_{n\to\infty} n^{1/2}\,
	\exp\{ -n\,\bPhistar(\al) \}\,\E\ZZ_{n\al}= \mathscr{C}(\al)\]
for $\bPhistar\equiv\bPhistar_d$, $\mathscr{C}\equiv\mathscr{C}_d$ smooth functions of $\al$. The explicit form of $\bPhistar$ is given in \S\ref{ss:explicit} below; in particular, on the interval $\albd\le\al\le\aubd$ it is strictly decreasing with a unique zero $\al_\star$ in the interval's interior.
\ethm

The proof is given in \S\ref{s:nd} (making use of \S\ref{s:coarsen}~and~\S\ref{s:first}).

\bcor\label{c:first.moment.ubd}
With $\al_\star$ the unique zero of $\bPhistar$ on the interval $\albd\le\al\le\aubd$ (Thm.~\ref{t:frozen.first.moment}), let
$c_\star \equiv -[ 2 \cdot \bPhistar'(\al_\star) ]^{-1}>0$. Then, for $\misn\equiv n\al_\star - c_\star\log n$ as in~\eqref{e:max.is.n}, $\E\ZZ_{\ge\misn+C}$ tends to zero as $C$ tends to infinity, uniformly over $n \ge n_0(d)$ .

\bpf
Taylor expansion gives $\E\ZZ_{\misn}\asymp_d1$ while $\E\ZZ_{\misn+5c_\star\log n}\le n^{-2}$. Recall from Thm.~\ref{t:frozen.first.moment} that $\bPhistar$ is decreasing on the interval $\albd\le\al\le\aubd$: since $\E\ZZ_{\ge n\al}$ is trivially a sum over at most $n$ terms $\E\ZZ_{n\al'}$, this immediately implies $\smash{\E\ZZ_{\ge\misn+5c_\star\log n} \le n^{-1}}$. For $0\le\de\le 5c_\star\log n$ we have $\smash{\E\ZZ_{\misn+\de}
	\asymp_d
	(\lm_\star)^{-\de}
	\exp\{ O_d( n^{-1}(\log n)^2 ) \}}$
where $\smash{\lm_\star\equiv e^{1/(2c_\star)}}$. Combining these estimates gives
$$\TS\E\ZZ_{\ge\misn+C}
\asymp_d
	(\lm_\star)^{-C}
	\exp\{ O_d( n^{-1}(\log n)^2 ) \}
	\sum_{\de\ge0} (\lm_\star)^{-\de}
 + n^{-1}
\lesssim (\lm_\star)^{-C},$$
implying the result.
\epf
\ecor

\bThm\label{t:lbd}
For $d\ge d_0$,
it holds uniformly over $n\ge n_0(d)$ that
$\P(\ZZ_{\misn-C}>0)$ tends to one as $C$ tends to infinity.
\eThm

\bpf[Proof of Thm.~\ref{t:main}]
We prove the theorem relying on Propns.~\ref{p:big}~and~\ref{p:tree} which will be proved in the remainder of this section.

\medskip\noindent\emph{Upper bound.}
Let $\Crs_{\ge n\al}$ denote the number of configurations $\ueta$ obtained by coarsening a \emph{maximum}\footnote{It would suffice for the independent set to be ``locally maximal'' in the sense that its coarsening is a configuration of the unweighted frozen model.} independent set of size $\ge n\al$. Write $\nodd(\FC)$ for the number of odd-size components in $\FC$; recall Defn.~\ref{d:frozen} that each odd component must contain at least one cycle. Define $\tZ_{n\aone,n\be}(O)$ to be the contribution to $\tZ_{n\aone,n\be}$ from configurations $\ueta$ with $\nodd[\FC(\ueta)]=O$, and write $\tZeven_{n\aone,n\be}\equiv \tZ_{n\aone,n\be}(O=0)$. In any coarsening $\ueta$ of a \emph{maximal} independent set $\ux$, there can be no isolated $\free$-vertices, and all tree $\free$-components must have a perfect matching. Thus, recalling Propn.~\ref{p:coarsening}, we have that with high probability $\Crs_{\ge n\al}$ is bounded above by
\beq\label{e:unweighted.frozen.ge}
\tZ_{\ge n\al}
\equiv
	\bigg[\sum_{\substack{\aone,\be\,:\,\be\le\bemax,\\
	2n\al\le 2n\aone+n\be}}
	\tZeven_{n\aone,n\be}\bigg]
+
	\bigg[\sum_{O>0} \sum_{
	\substack{
	\aone,\be\,:\,
	\be\le\bemax,\\
	2n\al+O
	\le 2n\aone+n\be
	}} \tZ_{n\aone,n\be}(O)\bigg]
\equiv \tZeven_{\ge n\al} + \tZodd_{\ge n\al}
\eeq
For $k\ge4$ even we shall set a threshold $\ell_k\equiv\ell_k(n\be)\ge1$ (to be determined), and further decompose
\beq\label{e:unweighted.frozen.even.decomp}
\tZeven_{\ge n\al}
\equiv
\tZtree_{\ge n\al}
+\tZbig_{\ge n\al}
+\overbrace{
	\sum_{A\ge1}\tZcyc_{\ge n\al}(A)
	}^{\equiv\,\tZcyc_{\ge n\al}}
\quad\text{where}
\eeq
\bnm[1.]
\item $\tZbig_{\ge n\al}$ is the contribution from configurations $\ueta$ 
such that the number of $\free$-components of size $k$ is
$\ge\ell_k(|\FC(\ueta)|)$ for some $k\ge4$ even;

\item $\tZcyc_{\ge n\al}(A)$ is the contribution from configurations $\ueta$ such that for all $k\ge4$ even there are fewer than $\ell_k$ $\free$-components of size $k$, where a minmum of $A$ edges must be removed from $\FC$ to make it a forest; and

\item $\tZtree_{\ge n\al}$ is the remainder, consisting of configurations $\ueta$ where every $\free$-component is a tree with (unique) perfect matching, and with $<\ell_k$ size-$k$ $\free$-components for all $k\ge4$ even.
\enm
All statements that follow are made 
uniformly over $\albd\le\al\le\aubd$
for $d\ge d_0$, $n\ge n_0(d)$.
We shall prove (Propns.~\ref{p:big}~and~\ref{p:tree}) that
$\E[\tZodd_{\ge n\al}
+\tZbig_{\ge n\al}
+\tZcyc_{\ge n\al}]
\le \tf12\E\tZ_{\ge n\al}$, consequently
\beq
\label{e:chain.exp}
\tf12\E\tZ_{\ge n\al}
\le
\E[\tZtree_{\ge n\al}]
=\E\ZZtree_{\ge n\al}
\le \E\ZZ_{\ge n\al}.
\eeq
Recall Cor.~\ref{c:first.moment.ubd} that $\E\ZZ_{\ge \misn+C}\to0$ uniformly in $n$ as $C\to\infty$, and therefore
$$
\P(\mis_n\ge \misn+C)
\le \P(\Crs_{\ge \misn+C}>0)
\le o_n(1)+\P(\tZ_{\ge \misn+C}>0)
$$
converges to zero as $C\to\infty$, uniformly in $n$. (In the above, the first inequality is by the coarsening algorithm, the second inequality is by Propn.~\ref{p:coarsening} with $o_n(1)$ indicating an error which tends to zero as $n\to\infty$, and the rightmost expression tends to zero by Markov's inequality and the chain of inequalities \eqref{e:chain.exp}.)

\medskip
\noindent\emph{Lower bound.}
In the weighted frozen model, by definition there are no odd-size components, so $\ZZ_{\ge n\al}$ has a similar decomposition as that of $\tZeven$ in \eqref{e:unweighted.frozen.even.decomp},
$$\TS
\ZZ_{\ge n\al}
\equiv
\ZZtree_{\ge n\al}+\ZZbig{\ge n\al}+\ZZcyc_{\ge n\al},\quad
\ZZcyc_{\ge n\al}\equiv
\sum_{A\ge1}\ZZcyc_{\ge n\al}(A)$$
We shall show (Propn.~\ref{p:big}(b) and Propn.~\ref{p:tree}(b)) that $\E\ZZbig{\ge n\al}\le n^{-1/2}\E\ZZ_{\ge n\al}$ and $\E\ZZcyc_{\ge n\al}(A)\le [d(\log n)/n]^A \E\ZZ_{\ge n\al}$. If $\uz=(\ueta,\match)$ contributes to $\ZZcyc_{\ge n\al}(A)$ then one can set $\eta_u=\eta_v=\zro$ on $A$ edges $(u,v)\in\match$ (removing those edges from the matching) to obtain a configuration $\uz'$ contributing to $\ZZtree_{\ge n\al-A}$, which immediately gives rise to independent sets of size $n\al-A$. Thus, with $\lm_\star\equiv \exp\{ -\bPhistar'(\al)\}$ as in the proof of Cor.~\ref{c:first.moment.ubd}, we have
\[\begin{array}{rl}
\P(\mis_n \ge \misn-2C)
\hspace{-6pt}&\ge
\P(\ZZ_{\ge \misn-C}>0)
-\E\ZZbig{\ge \misn-C}
-\sum_{A>C}\E\ZZcyc_{\misn-C}(A)\\
&\ge
\P(\ZZ_{\ge \misn-C}>0)
- e^{O_d(1)} (\lm_\star)^C
	(n^{-1/2} + [d(\log n)/n]^C ),
\end{array}\]
and by Thm.~\ref{t:lbd} this converges to one as $C\to\infty$, uniformly in $n$. Combining with the upper bound proves that $\mis_n-\misn$ is a tight random variable as claimed.
\epf

The above equivalence between the maximum independent set size and the threshold of the frozen model is based on the following

\bppn\label{p:big}
The following hold uniformly over 
$\albd\le\al\le\aubd$ for $d\ge d_0$, $n\ge n_0(d)$, with $\cc$ an absolute constant.
\bnm[(a)]
\item \label{p:big.a}
$\E[\tZodd_{\ge n\al}] \le d^{-1/2}\E\tZ_{\ge n\al}$.
\item \label{p:big.b} If for $k\ge4$ even we set $\ell_k\equiv[(\cc d\be)^{k/8}n/(dk)]\vee1$, then
$$\E[\tZbig_{\ge n\al}]\le n^{-1/2}\,\E\tZ_{\ge n\al}
\quad\text{and}\quad
\E\ZZbig{\ge n\al}\le n^{-1/2}\,\E\ZZ_{\ge n\al}.$$
\enm
\eppn

\bppn\label{p:tree}
The following hold uniformly over $\albd\le\al\le\aubd$ for $d\ge d_0$, $n\ge n_0(d)$, with $\ell_k$ as in the statement of Propn.~\ref{p:big}.
\bnm[(a)]
\item \label{p:tree.a}
	$\E[\tZcyc_{\ge n\al}]\le d^{-1/5}\E\tZ_{\ge n\al}$.
\item \label{p:tree.b}
	$\E\ZZcyc_{\ge n\al}(A) \le
	[ dn^{-1}\log n ]^A\E\ZZ_{\ge n\al}$.
\enm
\eppn

\subsection{Large components and trees with matchings}

We now prove Propns.~\ref{p:big}~and~\ref{p:tree}. Write $n_\one\equiv n\aone$, $n_\free\equiv n\be$, and $n_\zro\equiv n-n_\one-n_\free$. We decompose
\beq\label{e:Z.decomp}
\TS
\E\tZ_{n_\one,n_\free}
=\sum_{|\FC|=n_\free} \E[ \tZ_{n_\one}(\FC) ],\quad
\E\ZZ_{n_\one,n_\free}
=\sum_{|\FC|=n_\free} \mathfrak{m}(\FC)\,\E[ \tZ_{n_\one}(\FC) ]
\eeq
where $\tZ_{n_\one}(\FC)$ denotes the partition function of unweighted frozen model configurations $\ueta$ with $|\set{v:\eta_v=\one}|=n_\one$ and with $\free$-subgraph $\FC(\ueta)=\FC$, and the sum is taken over all possible $\free$-subgraphs $\FC$ containing $|\FC|=n_\free$ vertices. In keeping with Defn.~\ref{d:frozen}, we regard each $\FC$ as a matching on a subset of $[nd]$: thus the quantity $\tZ_{n_\one}(\FC)$ cannot be positive unless $\FC$ is present as an induced subgraph in the random graph $\Gnd$.

Observe that $\E[ \tZ_{n_\one}(\FC) ]$ depends on $\FC$ only through $n_\free$ and the number $\EFC\equiv\EFC(\FC)$ of internal edges in $\FC$:
writing $E_\eta\equiv n_\eta d$ and $\bE\equiv nd$, we calculate
\beq\label{e:unweighted.zof}
\E[ \tZ_{n_\one}(\FC) ]
=\mathbf{c}(n_\one,n_\free)
\,\gone(E,E_\one,E_\free)
\,\gzro( n_\zro,E_\one )
\,\gfreebar(E,E_\one,E_\free,\EFC)
\equiv \mathbf{g}(n_\one,n_\free,\EFC)
\eeq
where $\mathbf{c}(n_\one,n_\free)\equiv\tbinom{n-n_\free}{n_\one}$ enumerates the number of partitions of $V\setminus\FC$ into subsets $V_\zro,V_\one$ of sizes $n_\zro,n_\one$, and the remaining factors calculate probabilities with respect to a \emph{fixed} partition
$V_\zro,V_\one,V_\free$ of $V$
(with $V_\free$ the vertices involved in $\FC$), as follows:
\bnm[1.]
\item $\gone(E,E_\one,E_\free)$ is the probability all  half-edges from $V_\one$ are matched to half-edges from $V_\zro$;
\item $\gzro(n_\zro,E_\one)$ is the probability
each vertex in $V_\zro$ receives at least two half-edges from $V_\one$, with respect to a uniformly random assignment of the $E_\one$ incoming half-edges to the $n_\zro d$ available from $V_\zro$;
\item $\gfreebar(E,E_\one,E_\free,\EFC)$ is the probability, conditioned on the valid assignment of $E_\one$ edges between $V_\one$ and $V_\zro$, that the $\EFC$ internal edges of $\FC$ are present in the graph, while the remaining $E_\free-2\EFC$ half-edges from $V_\free$ are matched to half-edges from $V_\zro$.
\enm
After the assignment of $E_\one$ edges between $V_\one$ and $V_\zro$, there are $\bE-2E_\one=E_\zro-E_\one+E_\free$ half-edges remaining to be assigned. Recalling the notation \eqref{e:falling.factorial},
we have explicitly
\beq\label{e:unweighted.one}
\begin{array}{rl}
\gone(n_\zro,E_\one)
\hspace{-6pt}&= \smf{ (\bE-E_\one-E_\free)_{E_\one} }
	{ \fdf{\bE}{E_\one} }
	=\smf{(E_\zro)_{E_\one}}{\fdf{\bE}{E_\one}};\vspace{2pt}\\
\gfreebar(\bE,E_\one,E_\free,\EFC)
\hspace{-6pt}&=\smf{1}{\fdf{ \bE-2E_\one }{\EFC} }\times
\smf{ (
	E_\zro-E_\one
	)_{ E_\free-2\EFC } }
{\fdf{E_\zro-E_\one + E_\free-E_\one}{ E_\free-2\EFC } }
= \smf{ ( \bE-2E_\one-E_\free )_{E_\free-2\EFC} }
	{ \fdf{\bE-2E_\one}{E_\free-\EFC} }.
\end{array}
\eeq
Since $\mathbf{c},\gone,\gzro$ do not depend on $\EFC$ we find crudely that the function $\mathbf{g}$ of \eqref{e:unweighted.zof} satisfies
\beq\label{e:internal.edges}
\mathbf{g}(n_\one,n_\free,\EFC+\De)
\le [4/(nd)]^\De
	\, \mathbf{g}(n_\one,n_\free,\EFC)
\eeq
throughout the regime $\albd\le\al\le\aubd$, indicating that excess internal edges in $\FC$ are costly. On the other hand let us note that $\mathbf{g}$ is much less sensitive to small shifts in mass from $V_\zro$ to $V_\free$ or vice versa: for $\de \le n/d$ we estimate
\beq\label{e:set.free.to.zero}
\left.\begin{array}{rrl}
(\circ)
	&\mathbf{c}(n_\one,n_\free+\de)
	\hspace{-6pt}&= 
	\mathbf{c}(n_\one,n_\free)\,e^{O(\de)},\\
(\circ)
	&\gone(E,E_\one,E_\free+d\de)
	\hspace{-6pt}&=
	\gone(E,E_\one,E_\free) \, d^{O(\de)},\\
(\circ)
	&\gfreebar(E,E_\one,E_\free+d\de,\EFC)
	\hspace{-6pt}&= \gfreebar(E,E_\one,E_\free,\EFC)\,e^{O(\de)},\\
(\star) &
\gzro(n_\zro-\de,E_\one)
	\hspace{-6pt}&= 
	\gzro(n_\zro,E_\one)\,e^{O(\de)}
\end{array}\right\},
\text{ thus }
\f{\mathbf{g}(n_\one,n_\free+\de,\EFC)}
	{\mathbf{g}(n_\one,n_\free,\EFC)}
= d^{O(\de)}.
\eeq
The three estimates marked ($\circ$) follow straightforwardly from the explicit expressions given above; the last estimate ($\star$) is deferred to a later section (Propn.~\ref{p:forcing}).

\bpf[Proof of Propn.~\ref{p:big}]
Let $k\ge3$ be fixed, and let $\fc\equiv\fc(\FC)$ denote the
components of $\FC$ whose size is neither $2$ nor $k$. Write $T\equiv T(\FC)\equiv |\FC\setminus\fc|$, $t\equiv t(\FC)$ for the number of vertices in $\FC\setminus\fc$ in components of size $2$, and $\ell\equiv\ell(\FC)$ for the number of components of size $k$, so that $T=t+\ell k$. Let $\QFC\equiv\QFC(\FC)$ denote the number of internal edges among the components of size $k$. Any odd-size component must contain at least one cycle, so we define the number of ``excess'' edges among these components to be
\beq\label{e:excess.A}
A\equiv A(\FC)\equiv \QFC-\ell(k-\Ind{k\text{ even}})\ge0.\eeq
Write $\Om_{\ell,A}(n_\free;\fc)$ for the collection of all $\FC$ with $n_\free$ vertices and with the given $\fc,\ell,A$.

\medskip\noindent\emph{Bounds for $T$ even.}
Let $\Om(n_\free;\fc)\equiv\Om_{0,0}(n_\free;\fc)$
denote the collection of all $\FC$ with $n_\free$ vertices and the given $\fc$, such that $\FC\setminus\fc$ consists of $T$ vertices matched up into isolated edges, where we assume $T$ is even. Then all
$\FC\in\Om(n_\free;\fc)$ have the same number of internal edges $\EFC$ while all $\FC\in\Om_{\ell,A}(n_\free;\fc)$ have the same number of internal edges
\[\EFC'=\EFC + \QFC-\ell k/2
=\EFC + A + \ell( k/2-\Ind{k\text{ even}} ).\]
We calculate $|\Om(n_\free;\fc)| = \binom{n-|\fc|}{T} d^T (T-1)!!$, while
$$
|\Om_{\ell,A}(n_\free;\fc)|
\le 
\smb{n-|\fc|}{T}
	\smb{T}{t}
	d^t (t-1)!!
	\smf{(\ell k)!}{\ell!(k!)^\ell}
	d^{2\QFC} \smb{\ell k(k-1)/2}{\QFC}.$$
Combining with \eqref{e:internal.edges} and simplifying gives
$$
\mathbf{\bar{R}}^{n_\one,n_\free}_{\ell,A}(\fc)
\equiv\f{|\Om_{\ell,A}(n_\free;\fc)|}{|\Om(n_\free;\fc)|}
	\f{\mathbf{g}(n_\one,n_\free,\EFC')}{\mathbf{g}(n_\one,n_\free,\EFC)}
\le
\text{\footnotesize$\DS\f{(T/2)! 2^{\ell k/2}}{(t/2)!}$}
	\cdot
	\smf{1}{\ell! (k!)^\ell}
	\cdot
	\smb{\ell k(k-1)/2}{\QFC}
	\cdot
	\Big( \smf{4d}{n}\Big)^{\QFC-\ell k/2}.
$$
The first factor is $\le n_\free^{\ell k/2}$ since $t+\ell k= T\le n_\free$; the second factor is $\le e^{O(k\ell)}/[\ell k^k]^\ell$; and
the third factor is $\le e^{O(\QFC)} k^\QFC$ using $\QFC\ge\ell(k-1)$. Combining with \eqref{e:excess.A} gives
\beq\label{e:big.even.ubd}
\mathbf{\bar{R}}^{n_\one,n_\free}_{\ell,A}(\fc)
\le \Big(
	\smf{(\cc d\be)^{k/2}}{\ell}
	\Big( \smf{n}{dk} \Big)^{\Ind{k\text{ even}}}
	\Big)^\ell
	\Big( \smf{\cc dk}{n} \Big)^A
\quad\text{($T$ even, $\mathfrak{m}(\fc)>0$)}.
\eeq
for some constant $\cc$ uniform in $d$.
We consider also the effect of reweighting by the number of matchings on $\FC$, which is non-zero only if $\mathfrak{m}(\fc)>0$ and $k$ is even. If $\FC\in\Om(n_\free;\fc)$ then clearly $\mathfrak{m}(\FC)=\mathfrak{m}(\fc)$ since the collection $\FC\setminus\fc$ of isolated edges has only one matching. If $\FC\in\Om_{\ell,A}(n_\free;\fc)$ then the trivial bound $\mathfrak{m}(\FC)\le \mathfrak{m}(\fc) \, 2^\QFC$ implies
$$\mathbf{R}^{n_\one,n_\free}_{\ell,A}(\fc)
\equiv
\mathbf{\bar{R}}^{n_\one,n_\free}_{\ell,A}(\fc)
	\f{\max\set{\mathfrak{m}(\FC) : \FC\in\Om_{\ell,A}(n_\free;\fc)}}
	{\min\set{\mathfrak{m}(\FC) : \FC\in\Om(n_\free;\fc)}}
\le 
\Big( \smf{(\cc d\be)^{k/2}n}{dkl}\Big)^\ell
	\Big(\smf{\cc dk}{n}\Big)^A.$$

\medskip\noindent\emph{Bounds for $T$ odd.}
For $T$ odd we instead compare $\FC'\in\Om_{\ell,A}(n_\free;\fc)$
with $\FC\in\Om(n_\free-1;\fc)|$, which has cardinality
$\smash{\tbinom{n-|\fc|}{T-1} d^{T-1}(T-2)!!}$.
The difference in the number of internal edges between $\FC'$ and $\FC$
is $\EFC'-\EFC=\tf12(\ell k+1)+A$, so applying \eqref{e:internal.edges}
and \eqref{e:set.free.to.zero} and calculating as above gives
\beq\label{e:big.odd.ubd}
\mathbf{\bar{R}}^{n_\one,n_\free}_{\ell,A}(\fc)
\equiv
\f{|\Om_{\ell,A}(n_\free;\fc)|}{|\Om(n_\free-1;\fc)||}
\f{\mathbf{g}(n_\one,n_\free,\EFC')}{\mathbf{g}(n_\one,n_\free-1,\EFC)}
\le
d^{O(1)}
\Big( \smf{(\cc d\be )^{k/2}}{\ell}\Big)^\ell
\Big(\smf{\cc dk}{n}\Big)^A
\quad\text{($T$ odd)}.
\eeq

\medskip\noindent\emph{Odd-size components.}
Recall \eqref{e:unweighted.frozen.ge}
that $\tZ_{n_\one,n_\free}(O)$
denotes the contribution to $\tZ_{n_\one,n_\free}$
from frozen model configurations $\ueta$
with $\nodd(\FC)=O$ odd-size components in $\FC\equiv\FC(\ueta)$.
For $k\ge3$ odd and $1\le\ell\le O$
let $\tZ_{n_\one,n_\free}(O;L_k=\ell)$
denote the contribution to $\tZ_{n_\one,n_\free}(O)$ from configurations with $\ell$ $\free$-components of size $k$. Recalling \eqref{e:Z.decomp} and \eqref{e:unweighted.zof},
$$\E[\tZ_{n_\one,n_\free}(O;L_k=\ell)]
= \sum_{\substack{
	\fc: \nodd(\fc)+\ell=O, \\
	A\ge0 }
	}
\sum_{\FC\in \Om_{\ell,A}(n_\free;\fc)}
	\E[ \tZ_{n_\one}(\FC)]
=\sum_{\substack{
	\fc: \nodd(\fc)+\ell=O, \\
	A\ge0 }
	}
	|\Om_{\ell,A}(n_\free;\fc)|\,
	\mathbf{g}(n_\one,n_\free,\EFC')$$
Applying \eqref{e:big.even.ubd} ($\ell$ even) and \eqref{e:big.odd.ubd} ($\ell$ odd) gives, with $\overline{n}_\free\equiv n_\free-\Ind{\ell\text{ odd}}$,
\beq\label{e:decrease.odd.cycles}
\E[\tZ_{n_\one,n_\free}(O;L_k=\ell)]
\le 
\E[\tZ_{n_\one,\overline{n}_\free}(L-\ell;L_k=0)]
\,d^{O(1)}[(\cc d\be)^{k/2}\ell^{-1}]^\ell
	\quad(\ell\ge1).\eeq
For fixed $k\ge3$ odd and $\ell\ge1$,
summing both sides of \eqref{e:decrease.odd.cycles} over
$n_\one,n_\free,O$ with $O\ge\ell$ and
$2n\al + O \le 2n_\one+n_\free$ gives 
$\E[\tZ_{\ge n\al};L_k=\ell]
\le d^{-k\ell/5}\E\tZ_{\ge n\al}$
(cf.~\eqref{e:unweighted.frozen.ge}).
Summing over $k\ge3$ odd, $\ell\ge1$ then gives
$\E[\tZodd_{\ge n\al}] \le d^{-1/2}\E\tZ_{\ge n\al}$, concluding the proof of \eqref{p:big.a}.

\medskip\noindent\emph{Even-size components.}
Recall $\tZeven_{n_\one,n_\free}$ denotes the contribution to $\tZ_{n_\one,n_\free}$ from frozen model configurations $\ueta$ such that $\FC(\ueta)$ has \emph{no} odd-size components.
If $L_k$ denotes the number of size-$k$ components in $\FC$,
then
\beq\label{e:decrease.even.components}
\E[\tZeven_{n_\one,n_\free};L_k=\ell]
=\sum_{\fc:\nodd(\fc)=0}\sum_{A\ge0}|\Om_{\ell,A}(n_\free;\fc)|\,\mathbf{g}(n_\one,n_\free,\EFC').\eeq
Let $\ell_k\equiv[(\cc d\be)^{k/8}n/(dk)]\vee1$ as in the statement of the lemma: then for all $\ell\ge\ell_k$ we have from \eqref{e:big.even.ubd} that
$\mathbf{\bar{R}}^{n_\one,n_\free}_{\ell,A}(\fc)
\le (\cc d\be)^{(3/8)k\ell}  (\cc dk/n )^A$.
Further
$$(\cc d\be)^{(3/8)k\ell} \le \Big\{\hspace{-3pt}
	\begin{array}{lll}
	n^{-(3/16)\log d + O(1)} & \text{if }k\ell_k>\log n, \\
	{} [d (\log n)/n ]^{3\ell} & \text{if }k\ell_k\le\log n;
	\end{array}
$$ so in any case
we have crudely that
$(\cc d\be)^{(3/8)k\ell} < n^{-5/2}$.\footnote{As always, the statement should be understood to apply for $d\ge d_0$, $n\ge n_0(d)$.}
Thus from \eqref{e:decrease.even.components} and the definition of $\mathbf{\bar{R}}$ we have
$$\E[\tZeven_{n_\one,n_\free};L_k=\ell]
\le n^{-5/2} \E[\tZeven_{n_\one,n_\free};L_k=0]
\quad(k\ge4\text{ even}).$$
Summing over $\ell\ge1$, $k\ge4$,
and $n_\one,n_\free$ with $2n_\one+n_\free\ge2n\al$ proves
$\E[\tZbig_{\ge n\al}]\le n^{-1/2}\E\tZ_{\ge n\al}$. Essentially the same  argument (with $k$ even, $\mathbf{R}$ in place of $\mathbf{\bar{R}}$) gives $\E\ZZbig{\ge n\al}\le n^{-1/2} \E\ZZ_{\ge n\al}$, concluding the proof of \eqref{p:big.b}.
\epf

\bpf[Proof of Propn.~\ref{p:tree}]
Fix $k\ge4$ even and let $\ell<\ell_k$ with $\ell_k$ be as in Propn.~\ref{p:big}; it suffices to consider $k\le 20\log_d n$, otherwise $\ell_k=1$ for all $\be\le \bemax$.

Let $\fc^*\equiv\fc^*(\FC)$ denote the $\free$-components in $\FC$ which are not of size $k$. Let $\Omcyc_{\ell,r,A'}(n_\free;\fc^*)$ denote the collection of $\FC$ with $n_\free$ vertices and the given $\fc^*$, such that $\FC\setminus\fc^*$ consists of $\ell\equiv\tf1k(n_\free-|\fc^*|)$ components with $k$ vertices each, of which $\ell-r$ components are trees with (unique) perfect matching, and the remaining $r$ components are not trees and have a total of $\QFC'\equiv rk+A'$ internal edges with $A'\ge0$. Let $\Omtree(n_\free;\fc^*)\equiv\Omcyc_{0,0,0}(n_\free;\fc^*)$ denote the collection of $\FC$ with $n_\free$ vertices and the given $\fc^*$, such that each of the $\ell$ components of $\FC\setminus\fc^*$ is a tree with (unique) perfect matching on $k$ vertices. Any $\FC\in\Omtree(n_\free;\fc^*)$ has $\EFC$ internal edges any $\FC\in\Omcyc_{\ell,r,A'}(n_\free;\fc^*)$ has $\EFC'$ internal edges with $\EFC'-\EFC=\QFC'-r(k-1)=r+A'$.

Regarding $[dk]$ as a set of half-edges with the $i$-th half-edge incident to vertex $\ceil{i/k}$, let $\tpm_{d,k}$ denote the number of ways that $2(k-1)$ elements of $[dk]$ can be used to form a graph on $[k]$ which is a spanning tree with perfect matching (meaning a perfect matching \emph{in~the~tree}, not to be confused with a perfect matching of half-edges). Then, using \eqref{e:internal.edges},
$$\mathbf{\bar{S}}^{n_\one,n_\free}_{\ell,r,A'}(\fc^*)
\equiv
\f{|\Omcyc_{\ell,r,A'}(n_\free;\fc^*)|}
	{|\Omtree(n_\free;\fc^*)|}
\f{\mathbf{g}(n_\one,n_\free,\EFC')}
	{\mathbf{g}(n_\one,n_\free,\EFC)}
\le
\f{\tbinom{\ell}{r}
	(\tpm_{d,k})^s
	d^{2\QFC'}
	\tbinom{r k(k-1)/2}{\QFC'}
	}{(\tpm_{d,k})^\ell}
\Big(\smf{4}{nd}\Big)^{r+A'}.
$$
We crudely lower bound the number of spanning trees by the number of line graphs (which clearly have a perfect matching as $k$ is even), $\tpm_{d,k} \ge \tf12  k! \tf{[d(d-1)]^k}{(d-1)^2} \ge e^{O(k)} k^k d^{2(k-1)}$ with the $O(k)$ uniform in $d$. Consequently there exists an absolute constant $C$ not depending on $d$ such that
$\mathbf{\bar{S}}^{n_\one,n_\free}_{\ell,r,A'}(\fc^*)
\le  (C^k \ell d/n)^r (Cdk/n)^{A'}$.
By the trivial inequality $\mathfrak{m}(\FC)\le \mathfrak{m}(\fc^*) \, 2^{\QFC'}$, the same bound holds (adjusting $C$ as needed) after reweighting by the number of matchings:
$$\mathbf{S}^{n_\one,n_\free}_{\ell,r,A'}(\fc^*)
\equiv
\mathbf{\bar{S}}^{n_\one,n_\free}_{\ell,r,A'}(\fc^*)
\text{\small$\DS\f
	{\max\set{\mathfrak{m}(\FC) : \FC\in\Omcyc_{\ell,r,A'}(n_\free;\fc^*)}}
	{\min\set{\mathfrak{m}(\FC) : \FC\in\Omtree(n_\free;\fc^*)}}$}
	\le \Big( \smf{C^k \ell d}{n} \Big)^r
		\Big( \smf{Cdk}{n}\Big)^{A'}.$$
Recalling 
$k\le20\log_d n$ and
$\ell<\ell_k\equiv[(\cc d\be)^{k/8}n/(dk)]\vee1$ gives
$$\mathbf{\bar{S}}^{n_\one,n_\free}_{\ell,r,A'}(\fc^*),
\mathbf{S}^{n_\one,n_\free}_{\ell,r,A'}(\fc^*)
\le (cd\be)^{kr/9}[d(\log n)/n]^{A'}.$$
Recalling \eqref{e:Z.decomp} and arguing similarly as in the proof of Propn.~\ref{p:big} we have
\begin{align*}
&\E[\tZcyc_{n_\one,n_\free}]
\le
\sum_{\substack{
	k\ge4\text{ even},\\
	1\le r<\ell_k
	}}
\sum_{\substack{ 
	\fc^*:\nodd(\fc^*)=0,\\
	A'\ge0
	}}
	|\Omcyc_{\ell,r,A'}(n_\free;\fc^*)|\,
	\mathbf{g}(n_\one,n_\free,\EFC')
	\\
&\lesssim
\sum_{\substack{
	k\ge4\text{ even},\\
	1\le r<\ell_k
	}} (cd\be)^{kr/9}
\sum_{\fc^*:\nodd(\fc^*)=0}
	|\Omtree(n_\free;\fc^*)|\,\mathbf{g}(n_\one,n_\free,\EFC)
\le d^{-1/5}\E[\tZeven_{n_\one,n_\free}],
\end{align*}
and summing over $\aone,\be$ proves $\E[\tZcyc_{\ge n\al}]\le d^{-1/5}\E[\tZeven_{\ge n\al}]$ which implies \eqref{p:tree.a}. Similarly
$$
\E\ZZcyc_{\ge n\al}(A)
\le
\sum_{\substack{
	k\ge4\text{ even},\\
	1\le r<\ell_k
	}}
\sum_{\substack{\fc:\nodd(\fc^*)=0,\\ A'\ge A}}
	\mathfrak{m}(\FC)\,
	|\Omcyc_{\ell,r,A'}(n_\free;\fc^*)|\,
	\mathbf{g}(n_\one,n_\free,\EFC')
\le [d(\log n)/n]^A
	\E\ZZ_{\ge n\al},$$
proving \eqref{p:tree.b}.
\epf

\subsection{Estimates on forcing constraints}\label{ss:zero.forcing}

In this subsection we estimate the probability cost of the constraint that each $\zro$-vertex is forced by at least two neighboring $\one$-vertices.

Slightly more generally, for $\minval$ a fixed positive integer,
write $\gzrominval(n,E_\one)$ for the probability,
with respect to a uniformly random assignment of $E_\one$ half-edges to $n$ vertices of degree $d$, that each vertex receives at least $\minval$ incoming half-edges (thus the function $\gzro$ of \eqref{e:unweighted.zof} is simply $\gzrominval$ with $\minval=2$). Let $X,X^1,\ldots,X^n$ be i.i.d.\ $\mathrm{Bin}(d,\vth)$ random variables, with joint law $\P_\vth$: then
$$\TS
\gzrominval(n,E_\one)
=\P_\te( X^i\ge\minval\,\forall 1\le i\le n \giv
	 \sum_{i=1}^n X^i = E_\one).
$$
where $0<\vth<1$ may be arbitrarily chosen.

We also give estimates on a two-dimensional analogue $\gzzminval$ which will be used in our second moment analysis \S\ref{s:second}:
let $\spins\equiv\set{\zro,\one}^2\setminus\set{\zz}$. If $\te\equiv (\te_\om)_{\om\in\spins}$ is a positive vector with sum less than $1$, then it defines a positive probability measure on $\set{\zro,\one}^2$; and we take $\vpi$ to be a $\set{\zro,\one}^2$-valued random variable with law given by this measure. Take $(\vpi_j)_{j\ge1}$ independent random variables identically distributed as $\vpi$, and define the multinomial random variable
\beq\label{e:multinomial}
\TS X\equiv(X_\om)_{\om\in\spins}
	\equiv (\sum_{j=1}^d \Ind{\vpi_j=\vpi})_{\om\in\spins}.
\eeq
Though we always use $\te$ and $X$ to denote three-dimensional vectors, we shall also write $\smash{\te_\zz\equiv 1-\sum_{\om\in\spins}\te_\om}$ and $\smash{X_\zz\equiv d-\sum_{\om\in\spins}X_\om}$. Let $X^i$ ($1\le i\le n$) be independent random variables identically distributed as $X$, and write $\P_\te$ for their (joint) law. Then, with $\badset\equiv\set{ x : x_\od\wedge x_\doto<\minval }$, we consider
$$\TS \gzzminval(n,\vec E) \equiv \P_\te( X^i \notin\badset \,\forall
	1\le i\le n \giv \sum_{i=1}^n X^i = \vec E ).$$

\bppn\label{p:forcing}
Let $\erho$ be a small constant uniform in $d$, and suppose
$$
\erho \logdbyd \le \set{\target[\od],\target[\doto]} \le 
	10\logdbyd.$$
\bnm[(a)]
\item \label{p:forcing.a}
Let $a,y$ be defined by
$\target[\od] \equiv a\logdbyd$
and $\target[\od]\wedge\target[\doto]\equiv y\logdbyd)$.
Then
$$
\begin{array}{rl}
\gzrominval(n,nd\target[\od])
	\hspace{-6pt}&=\exp\{ O( n d^{-a}(\log d)^{\minval} ) \},\\
\gzzminval(n,nd\target)
	\hspace{-6pt}&=\exp\{ O( nd^{-y} (\log d)^{\minval} ) \}.
\end{array}
$$
\item \label{p:forcing.b}
If $\targetp$ is another vector with
$|\targetp_\om/\ze_\om-1|\le e^{-1}$
for all $\om\in\spins$, then
\[\begin{array}{rl}
\gzrominval(n,nd\targetp_\od)
	\hspace{-6pt}&=
	\gzrominval(n,nd\target[\od])\,
	\exp\{ O( nd|\targetp[\od]-\target[\od]| ) \},\\
\gzzminval(n,nd\targetp)
	\hspace{-6pt}&=
	\gzzminval(n,nd\target)\,
	\exp\{ O( nd\|\targetp-\target\|_1 ) \}.
\end{array}\]
\enm
\eppn

\blem\label{l:lagrange.multipliers.bound}
For the multinomial random variable $X$ defined
by \eqref{e:multinomial}, 
consider the cumulant generating function
$\smash{\Lm_\origin(\gm)
	\equiv\log\E_\origin[ e^{ \ip{\gm}{X} } \giv X\notin\badset ]}$,
defined for $\gm\in\R^3$. 
For positive vectors $\origin,\target$ in the regime
$$\erho\logdbyd \le 
\set{\origin_\od,\origin_\doto,\target[\od] ,\target[\doto]}
\le 10\logdbyd,$$
there exists a unique $\gm\in\R^3$
with $\Lm_\origin'(\gm)=d\target$;
and it satisfies
$$\Big|
\gm_\om - \log\f{\origin_\zz \target[\om]}{\target[\zz]\origin_\om}
\Big|\lesssim \f{(\log d)^{\minval-1}}{d^y}\quad
\text{for all }\om\in\spins.$$

\bpf
If $\gm$ exists, then it is clearly unique by the strict convexity of
the cumulant generating function on $\R^3$. To see existence, first note that
$$
\Lm_\origin'(\gm) = \E_\reweigh[ X \giv X\notin\badset ]
\quad\text{for }
	\reweigh_\om
	= \f{\origin_\om e^{\gm_\om}}{ z_{\origin,\gm}}
	\equiv 
	\f{\origin_\om e^{\gm_\om}}
		{ \origin_\zz + \sum_{\pi\in\spins} \origin_\pi e^{\gm_\pi} }
	\text{ ($\om\in\spins$)}.
$$
In particular, each entry of $\Lm_\origin'(\gm)$ lies between $0$ and $d$. Consider the $\ep$-perturbed ($\ep>0$ small) cumulant generating function $\Lm_\origin(\gm) + \tf12 \ep \|\gm\|^2$, corresponding to the random variable $X+\ep^{1/2} Y$ for $Y$ a standard Gaussian in $\R^3$. For any fixed $\ep>0$ this is a smooth convex function on $\R^3$, with gradient $\Lm_\origin'(\gm) + \ep\gm$ tending in norm to $\infty$ as $\|\gm\|\to\infty$. It follows from Rockafellar's theorem (see e.g.\ \cite[Lem.~2.3.12]{MR2571413}) that there exists a unique $\gm_\ep\equiv(\gm_{\ep,\om})_{\om\in\spins}$ such that $\Lm_\origin'(\gm_\ep) +\ep\gm_\ep = d\target$. We shall show by some rough estimates that $\gm_\ep$ must remain within a compact region as $\ep$ tends to zero. We claim first that $\origin_\zz \ge \origin_\om e^{\gm_{\ep,\om}}$ for all $\om\in\spins$: if not, then for some $\om\in\spins$ we must have $z_{\origin,\gm} \le 4\origin_\om e^{\gm_{\ep,\om}}$, implying (in the stated regime of $\origin,\target$) that
$$\gm_{\ep,\om}
= \ep^{-1}[d\target[\om] - \pd_\om\Lm_\origin(\gm_\ep)]
\le \ep^{-1} [d\target[\om] - d/4] \ll0,$$
contradicting the hypothesis that
$\origin_\zz < \origin_\om e^{\gm_{\ep,\om}}$.
Thus $\limsup_{\ep\downarrow0}\gm_{\ep,\om}$
must be finite for each $\om\in\spins$.
In the other direction,
the trivial bound $z_{\origin,\gm}\ge\origin_\zz$ gives
$$
\gm_{\ep,\om}
=\ep^{-1}[ d\target[\om] - \pd_\om \Lm_\origin(\gm_\ep)]
\ge \ep^{-1} d\target[\om](1 - e^{\gm_{\ep,\om}}/\origin_\zz),
$$
so clearly $\liminf_{\ep\downarrow0}\gm_{\ep,\om}$ must also be finite.
It follows by an easy compactness argument that
$\gm_\ep$ converges in the limit $\ep\downarrow0$
to the required solution $\gm$ of $\Lm_\origin'(\gm)=d\target$.

To control the norm of the solution $\gm$ of $\Lm_\origin'(\gm)=d\target$, we shall argue that for $\target$ in the stated regime, $\reweigh_\om$ is close to $\target[\om]$ for each $\om$. The bound for $\om=\oo$ is easiest: consider $X$
as the $d$-th step of the random walk
$$\TS X_t\equiv
	(X_{t,\om})_{\om\in\spins}
	\equiv
	(\sum_{j=1}^t \Ind{\vpi_j=\om})_{\om\in\spins}.$$
Define the stopping time
$\tau\equiv\inf\set{t\ge0: X_t\notin\badset}$,
so $\set{X\notin\badset} = \set{\tau\le d}$.
Since $X_{\tau,\oo} \le\minval$,
applying the Markov property gives
\beq\label{e:ubd.one.one}
\begin{array}{rl}
d\target[\oo]
\hspace{-6pt}&=\E_\reweigh[X_\oo \giv X\notin\badset]
\le \minval + \E_\reweigh[ X_{d,\oo}-X_{\tau,\oo} \giv\tau\le d ]\\
&\le \minval + \E_\reweigh[X_{d-\tau,\oo}]
\le \minval + d\reweigh_\oo.
\end{array}
\eeq
Next observe that for the multinomial random variable $X$ defined by \eqref{e:multinomial},
for any $\om\ne\pi$ 
the conditional expectation
$\E_\te[X_\om\giv X_\pi=k]
=(d-k)\te_\om/(1-\te_\pi)$ is decreasing with $k$, so
\beq\label{e:multinomial.ubd}
\E_\te[X_\om\giv X_\pi\ge\ell]
= \sum_{k\ge\ell}
\P( X_\pi=k \giv X_\pi\ge\ell )
\E_\te[X_\om\giv X_\pi=k]
\le \E_\te X_\om = d\te_\om.
\eeq
Define the stopping times
$\tau_\od\equiv\inf\set{t\ge0 : X_\od(t) \ge\minval}$
and symmetrically $\tau_\doto$;
then $\tau = \tau_\od\vee\tau_\doto$.
Since $X_{\tau_\od,\oz}\le\minval$,
\begin{align*}
d\target[\oz]
&=\E_\reweigh[X_\oz\giv X\notin\badset]
= \f{\E_\reweigh[ X_\oz \Ind{\tau\le d} ]}{\P_\reweigh(\tau\le d)}
\le \minval
	+ \f{\E_\reweigh[ (X_\oz - X_{\tau_\od,\oz}) \Ind{\tau\le d} ]}{\P_\reweigh(\tau\le d)}\\
&=\minval
+ \sum_{k<d}\sum_x
\f{\P_\reweigh(\tau_\od=k, X_{\tau_\od}=x,
	\wt X_{d-k,\doto} \ge\minval-x_\doto)}
	 {\P_\reweigh(\tau\le d) }
	 \E_\reweigh[ \wt X_{d-k,\oz} \giv \wt X_{d-k,\doto} \ge \minval - x_\doto ],
\end{align*}
where $(\wt X_t)_{t\ge0}$ is an independent realization of the random walk $X$. Maximizing over all possible $k,x$ and applying \ref{e:multinomial.ubd} gives
$$\begin{array}{l}
d\target[\oz] \le
	\minval + \max_{k<d, \ell\le\minval}
	\E_\reweigh[ \wt X_{d-k,\oz} \giv \wt X_{d-k,\doto} \ge \ell ]
	\le
	\minval + d\reweigh_\oz,\\
\text{and symmetrically }
	d\target[\zo]\le\minval + d\reweigh_\zo.
\end{array}$$
Thus for $\target$ in the stated regime we conclude that
$d\reweigh_\om \ge d\target[\om]-O(1)$ for all $\om\in\spins$,
therefore $d\reweigh_\zz\le d\target[\zz] +O(1)$.

From these bounds we see clearly that,
with $y$ defined by 
$\target[\od]\wedge\target[\doto]\equiv y\logdbyd$,
$$\P_\reweigh(X\in\badset)
\le
\P_\reweigh(X_\od<\minval)+\P_\reweigh(X_\doto<\minval)
\lesssim d^{-y} (\log d)^{\minval-1}.
$$
Write $p_\reweigh$ for the law of $X\equiv X_d$
and $q_\reweigh$ for the law of $X_{d-1}$,
and observe that $x_\om p_\reweigh(x) = d\reweigh_\om q_\reweigh(x-\I_\om)$
(with both sides zero for $x_\om=0$). We then calculate
\beq\label{e:theta.mu.reln}
d\target[\om]
=\f{ \E_\reweigh X_\om - d\reweigh_\om \sum_{x\in\Om} q_\reweigh(x-\I_\om) }
	{1-\P_\reweigh(X\in\Om)}
= d\reweigh_\om\,[ 1 + O( d^{-y}(\log d)^{\minval-1} ) ]
\eeq
for all $\om\in\set{\zro,\one}^2$.
Summing over $\om\in\spins$ gives
\beq\label{e:reweight.normalizing}
z_{\origin,\gm}
=(\origin_\zz/\target[\zz])(\target[\zz]/\reweigh_\zz)
	= (\origin_\zz/\target[\zz])[1 + O( d^{-1-y}(\log d)^{\minval} )],
\eeq
implying the stated bound
on $\gm_\om=\log (z_{\origin,\gm} \reweigh_\om /\origin_\om)$,
$\om\in\spins$.
\epf
\elem

It follows (see e.g.\ \cite[Lem.~2.3.9]{MR2571413})
that with $\origin,\target$ in the stated regime,
 the Fenchel-Legendre transform
$\Lm_\origin^*(d\target)
\equiv\sup_\gm[ \ip{\gm}{d\target}-\Lm_\origin(\gm) ]$
of the cumulant generating function
is given by
\beq\label{e:fenchel.legendre}
\Lm_\origin^*(d\target)
= \ip{\gm}{d\target} - \Lm(\gm),\quad
\gm\text{ the solution of $\Lm_\origin'(\gm)=d\target$.}
\eeq
Since $\Lm_\origin$ is strictly convex,
we find by implicit differentiation 
that $\gm$ is differentiable with respect to $\target$
(in the stated regime).
We then see from \eqref{e:fenchel.legendre}
that $\Lm_\origin^*$ is differentiable
with respect to $\target$, with gradient
$(\Lm_\origin^*)'(d\target)=\gm$.

\bpf[Proof of Propn.~\ref{p:forcing}]
Let $\nu$ denote the normalized empirical measure of $X^1,\ldots,X^n$.
The $\P_\origin$-probability of realizing $\nu$ is $\smash{\tbinom{n}{n\nu}\prod_x p_\origin(x)^{n\nu(x)}}$
where $p_\origin$ denotes the law of $X^1$.
Denote the empirical mean
$\ol{\nu}\equiv\sum_x x\nu(x)$;
then $\gzzminval(n,nd\target)$
can be expressed as a constrained optimization over
 empirical measures $\nu$ with $\nu(\badset)=0$
 and  $\ol{\nu}=d\target$:
\begin{align*}
&\f{\gzzminval(n,nd\target)}{[1-\P_\origin(\badset)]^n}
= 
\f{\P_\origin( \sum_{i=1}^n X^i=nd\target
	, X^i\notin\badset\,\forall i)}
	{\P_\origin(X^i\notin\badset\,\forall i)\,
	\P_\origin( \sum_{i=1}^n X^i=nd\target)}\\
&=\f{ e^{ -n\ip{\gm}{d\target} }
	\sum_{\nu: \nu(\badset)=0,\ol{\nu}=d\target}
	\tbinom{n}{n\nu}
	\prod_x
	[p_\origin(x) e^{ \ip{\gm}{x} } ]^{n\nu(x)}
	}
	{ [1-\P_\origin(\badset)]^n
	\P_\origin( \sum_{i=1}^n X^i=nd\target) },
\end{align*}
where we introduced a Lagrangian term which clearly has no effect on the constrained space.
In the denominator,
Stirling's approximation gives
$$\TS
\P_\origin( \sum_{i=1}^n X^i=nd\target)
\asymp_d n^{-3/2}\exp\{ -nd\relent{\target}{\origin} \}$$
where 
$\relent{\target}{\origin}
\equiv\sum_{\om\in\set{\zro,\one}^2}
	\target[\om]\log(\target[\om]/\origin_\om)$
denotes the relative entropy between the probability measures induced by $\target$ and $\origin$ on $\set{\zro,\one}^2$.
In the numerator,
if we take $\gm$ to be the unique solution $\gm$
given by Lem.~\ref{l:lagrange.multipliers.bound} of
$\Lm_\origin'(\gm)=d\target$,
then we see that the sum over
$\set{\nu(\badset)=0, \ol{\nu}=d\target }$
equals (up to constants depending on $d$)
the sum over $\set{\nu(\badset)=0}$
divided by $n^{3/2}$:
thus, recalling \eqref{e:fenchel.legendre},
$$\gzzminval(n,nd\target)
\asymp_d
[1-\P_\origin(\badset)]^n\,
\exp\{ -n \Lm_\origin^*(d\target)
	+ nd\relent{\target}{\origin} \}.
$$
We can estimate this easily by taking $\origin=\target$.
Since $e^{\ip{\gm}{x}} p_\origin(x) = (z_{\origin,\gm})^d p_\reweigh(x)$,
it follows from \eqref{e:theta.mu.reln}~and~\eqref{e:reweight.normalizing} that
$$\Lm_{\target}(\gm)
= d\log z_{\target,\gm}
	+ \log \f{\P_\reweigh(X\notin\badset)}{\P_{\target}(X\notin\badset)}
= O(d^{-y} (\log d)^{\minval})
\quad\text{for }\Lm_{\target}'(\gm)=d\target.
$$
Combining with \eqref{e:fenchel.legendre}
and the bound of Lem.~\ref{l:lagrange.multipliers.bound}
gives
$\Lm_{\target}^*(d\target) = O(d^{-y} (\log d)^{\minval})$.
To compare $\gzzminval(n,nd\target)$ with $\gzzminval(n,nd\targetp)$
where $\targetp_\om = \target[\om][1+O(\de/n)]$,
we fix $\origin=\target$.
Recalling
$(\Lm_\origin^*)'(d\target)=\gm$,
the bound of Lem.~\ref{l:lagrange.multipliers.bound}
gives
$$
|n[\Lm_{\target}^*(d\targetp)-\Lm_{\target}^*(d\target)]|
\lesssim  n d^{-y}(\log d)^{\minval-1}
\cdot \|d\targetp-d\target\|_1. 
$$
For
$|\targetp_\om/\ze_\om-1|\le e^{-1}$
it is straightforward to estimate
$nd\relent{\targetp}{\target}
\lesssim nd\|\targetp-\target\|_1$,
and combining these estimates concludes the proof of the proposition.
\epf

\section{First moment of frozen model}\label{s:first}

In this section we identify the leading exponential order $\bPhistar(\al) = \lim_{n\to\infty} n^{-1}\log\E\ZZ_{n\al}$ of the first moment of the (truncated) frozen model partition function \eqref{e:truncated.frozen.model}. The random $d$-regular graph $\Gnd$ converges locally weakly (in the sense of \cite{MR1873300,MR2354165}) to the infinite $d$-regular tree $\tree$. Our calculation is based on a variational principle which relates the exponent $\bPhistar(\al)$ to a certain class of Gibbs measures for the frozen model on $\tree$ which are characterized by fixed-point recursions. In fact the recursions can have multiple solutions, and much of the work goes into identifying (via \emph{a~priori} estimates) the unique fixed point which gives rise to $\bPhistar(\al)$. We begin by introducing the Gibbs measures which will be relevant for the variational principle.

\subsection{Frozen model tree recursions}
\label{ss:frozen.tree}

We shall specify a Gibbs measure $\nu$ on $\tree$ by defining a consistent family of finite-dimensional distributions $\nu_t$ on the depth-$t$ subtrees $\tree(t)$. A typical manner of specifying $\nu_t$ is to specify a law on some boundary conditions at depth $t$, and then to define $\nu_t$ as the law of the configuration on $\tree(t)$ given the (random) boundary conditions.

In our setting some difficulty is imposed by the fact that the frozen model is not a factor model (or Markov random field) in the conventional sense that $\ueta|_A$ and $\ueta|_B$ are conditionally independent given the configuration $\ueta|_C$ on any subset $C$ separating $A$ from $B$ --- in particular, given the spins at level $t$ of $\tree$, whether a vertex at level $t-1$ is required to take spin $\one$ depends on whether its neighboring $\zro$'s in level $t$ are forced or not by $\one$'s in level $t+1$. Also, the frozen model spins $\ueta$ do not encode the matching $\match$ on the $\free$-vertices.

We shall instead specify Gibbs measures for the frozen model via a message-passing system, as follows. Each vertex $v$ will send a message $\si_{v\to w}$ to each neighboring vertex $w\in\pd v$ which represents the ``state of $v$ ignoring $w$''; this will be a function $\dmp_{d-1}$ of the other $d-1$ incoming messages $(\si_{u\to v})_{u\in\pd v\setminus w}$. The actual state $\eta_v$ of $v$ is then a function $\dmp_d$ of all its incoming messages $\dsi_{\pd v\to v}\equiv (\si_{u\to v})_{u\in\pd v}$. If on the boundary of $\tree(t)$ we are given a vector $\ueta^\uparrow\equiv(\si_{v\to w})_{v,w}$ (for $v$ at level $t$, $w$ the parent of $v$), then there is a \emph{unique} completion of $\ueta^\uparrow$ to a (bi-directional) message configuration on $\tree(t)$: iterating $\dmp$ gives all the messages upwards in the direction of the root, and once those are known we can recurse back down to determine the messages in the opposite direction. The measure $\nu_t$ can then be specified by giving the law of the boundary messages $\ueta^\uparrow$ (possibly with an additional reweighting): our choice will be to take $\ueta^\uparrow$ i.i.d.\ according to a law $q$; consistency of the family $(\nu_t)_t$ will then amount to fixed-point relations on $q$.

\begin{figure}[ht]
\centering
\begin{subfigure}[h]{0.35\textwidth}
	\centering
	\includegraphics[height=1.5in,trim=.6in .5in .8in .6in,clip]%
		{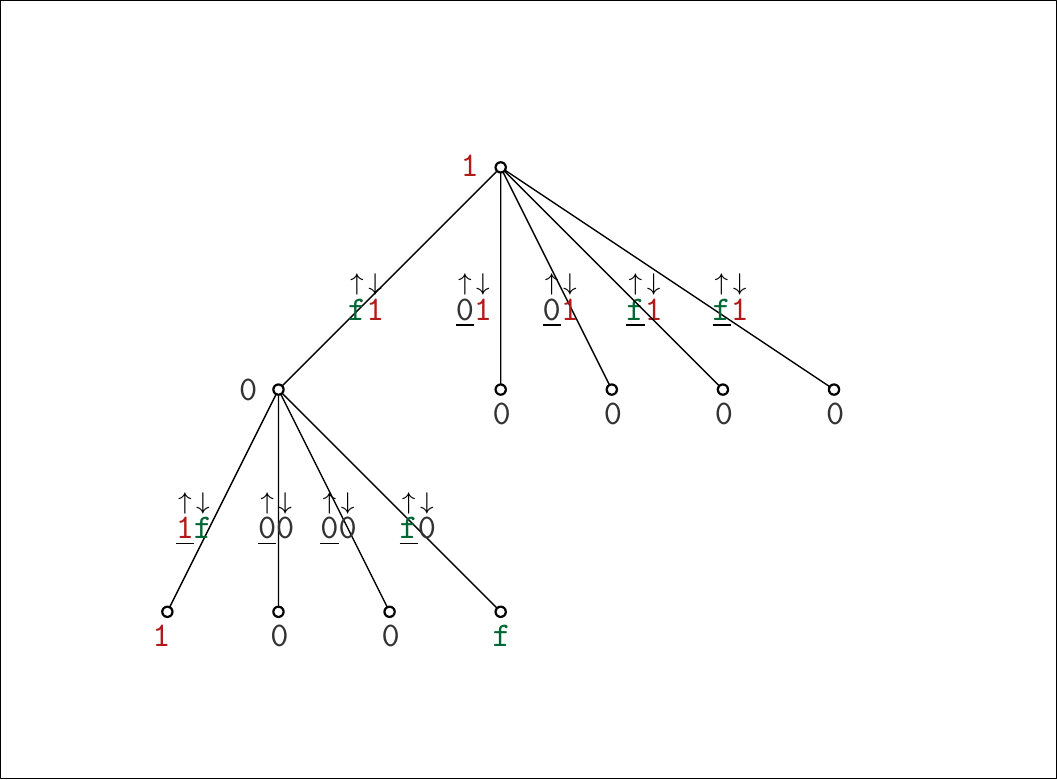}
	\caption{Root spin $\eta_\rt=\one$}
\end{subfigure}
\quad
\begin{subfigure}[h]{0.35\textwidth}
	\centering
	\includegraphics[height=1.5in,trim=.6in .5in .8in .6in,clip]%
		{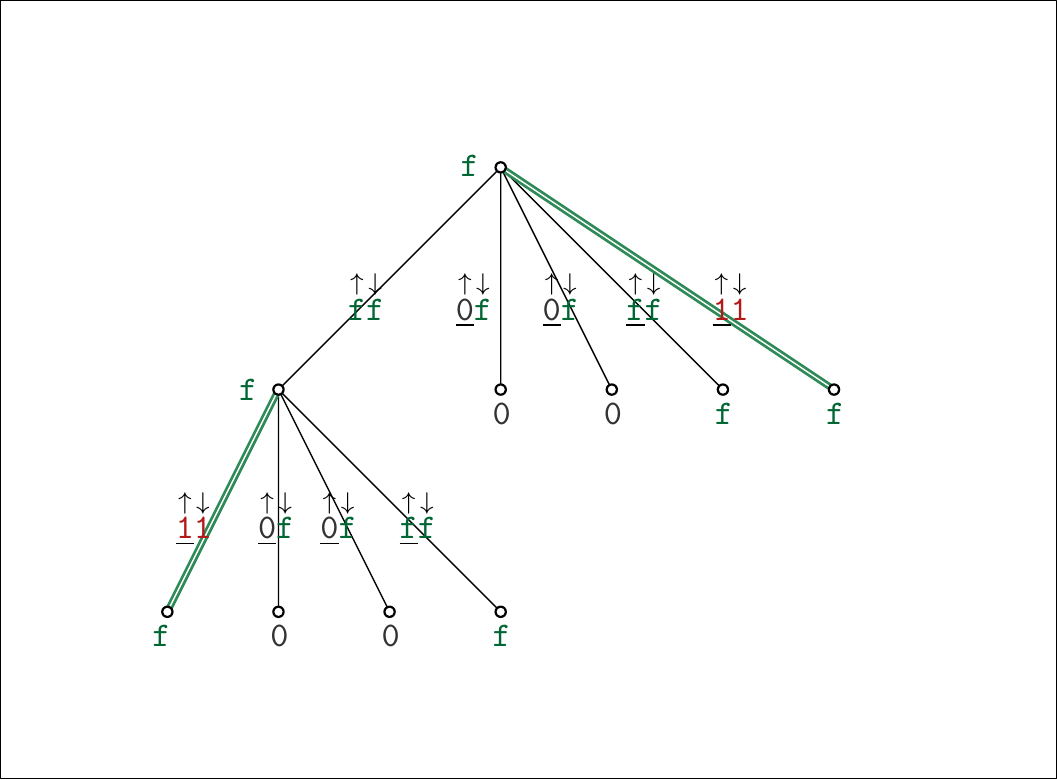}
	\caption{Root spin $\eta_\rt=\free$}
\end{subfigure}\\
\centering
\begin{subfigure}[h]{0.35\textwidth}
	\centering
	\includegraphics[height=1.5in,trim=.6in .5in .8in .6in,clip]%
		{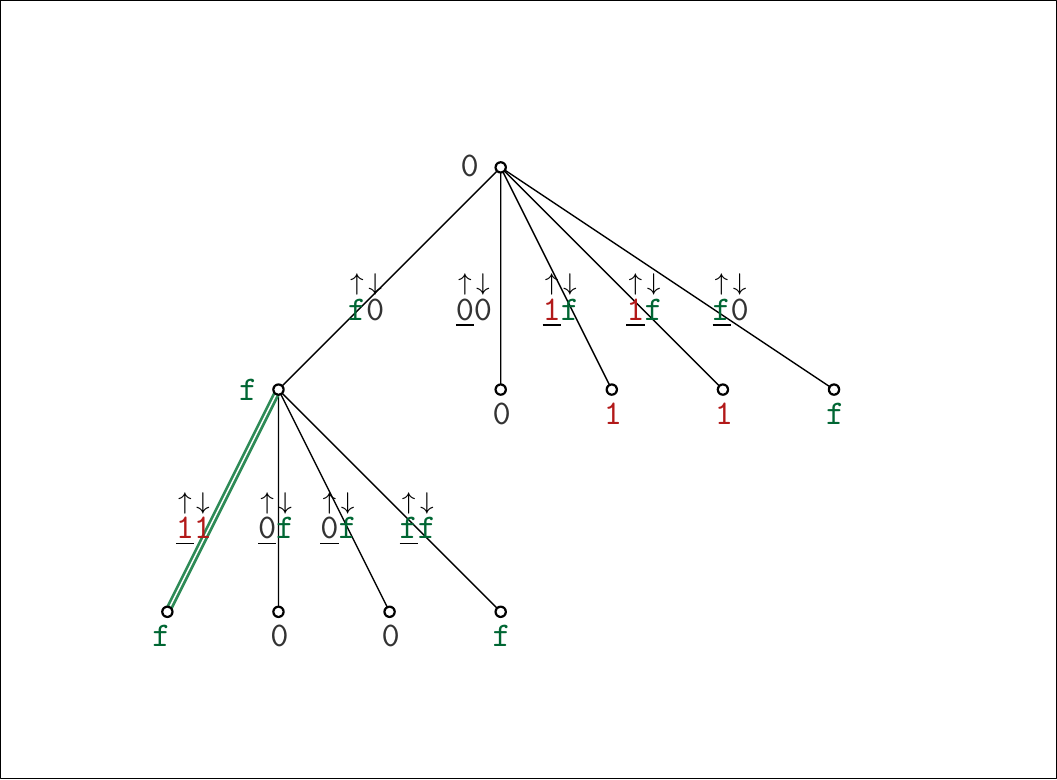}
	\caption{Root spin $\eta_\rt=\zro$}
\end{subfigure}
\quad
\begin{subfigure}[h]{0.35\textwidth}
	\centering
	\includegraphics[height=1.5in,trim=.6in .5in .8in .6in,clip]%
		{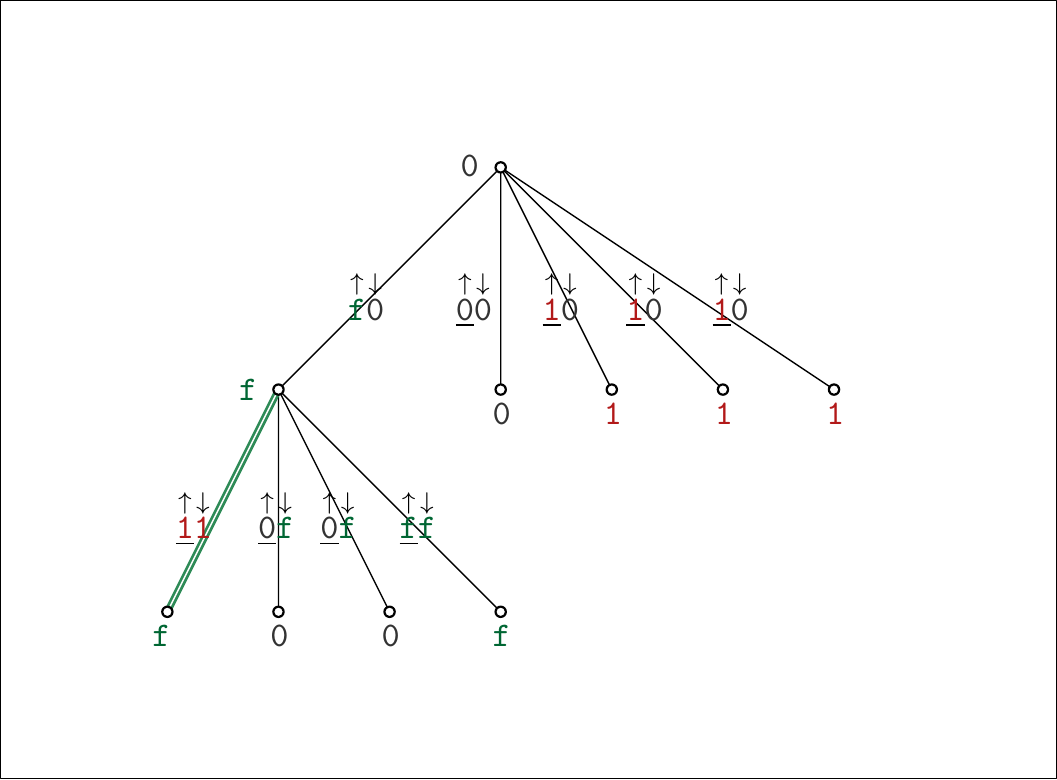}
	\caption{Root spin $\eta_\rt=\zro$}
\end{subfigure}
\caption{Completion of incoming boundary messages (underlined)}
\label{f:tree.frozen.messages}
\end{figure}

The message-passing rule for our frozen model is
$$\text{\small$\dmp_D:\set{\zro,\one,\free}^D \to\set{\zro,\one,\free},\quad
\dmp_D(\dot\ueta)
=\begin{cases}
	\one, & |\set{i:\eta_i=\one}|=0,\\
	\free, &|\set{i:\eta_i=\one}|=1,\\
	\zro, &|\set{i:\eta_i=\one}|\ge2.
	\end{cases}$}$$
See Fig.~\ref{f:tree.frozen.messages} for an illustration: in each panel of the figure, the entire configuration of messages and vertex spins is uniquely determined by the messages incoming at the boundary of the subtree depicted.
        
We then define the $\nu_t$-probability of a configuration on $\tree(t)$ to be (up to a normalizing factor) the product measure $q(\ueta^\uparrow)\equiv\prod_i q(\eta^\uparrow_i)$ tilted by the factor $\lm>1$ raised to the number of \emph{upward} (in the direction away from the boundary) $\one$-messages from vertices in $\tree(t-1)$, \emph{including} $\eta^\uparrow_\rt\equiv\dmp_d((\si_{v\to\rt})_{v\in\pd\rt})\equiv\eta_\rt$. Thus
$$
\nu_1(\eta_\rt=\one)
= \smf{\lm(1-\qo)^d}{(\lm-1)(1-\qo)^d + 1},\quad
\nu_1(\eta_\rt=\free)
= \smf{d \qo(1-\qo)^{d-2}}{(\lm-1)(1-\qo)^d + 1},$$
with the remaining probability going to $\eta_\rt=\zro$. The $\nu_t$ are consistent if and only if $q$ satisfies the \emph{frozen model recursions}
\beq\label{e:frozen.recursions}
\qo
=\smf{\lm(1-\qo)^{d-1}}{(\lm-1)(1-\qo)^{d-1} + 1},\quad
\qf=\smf{(d-1)\qo(1-\qo)^{d-2}}{(\lm-1)(1-\qo)^{d-1} + 1}
=\smf{(d-1)\qo^2}{\lm(1-\qo)},
\eeq
and $\qz=1-\qo-\qf$. (In contrast, the hard-core model with fugacity parameter $\lm$, which allows for $\zro$-vertices with any number of neighboring $\one$-vertices, has recursion $q=\lm(1-q)^{d-1}/[ \lm(1-q)^{d-1}+1 ]$ on $\tree$.)

A solution $q\equiv(\qo,\qf,\qz)$ of the frozen model recursions \eqref{e:frozen.recursions} corresponds to a root $\qo$ of the function $f(q)= (1-q)^{d-1}(\lm+q-\lm q)-q$ which, for $0\le q\le 1$ and $\lm>1$, is decreasing in $q$ and increasing in $\lm$. Therefore \eqref{e:frozen.recursions} is solved by a unique measure $q$ with $\qo$ increasing in $\lm$. Note that if $\qo=y\logdbyd$ with $y\asymp1$ then
\beq\label{e:lambda.of.q}
\lm = 
\smf{\qo [1-(1-\qo)^{d-1}]}{(1-\qo)^d}
= e^{o_d(1)} d^{y-1} y\log d.
\eeq
In the following, to emphasize the dependence on $\lm$ we sometimes write $\nu\equiv \nu^\lm$, $q\equiv q^\lm$.

\subsection{Auxiliary model and Bethe variational principle}
\label{ss:message.aux}

On the tree $\tree$, the frozen model configuration $\uz=(\ueta,\match)$ can be uniquely recovered from the configuration $\usi$ of messages on all the directed edges: each vertex spin $\eta_v$ is determined by applying $\dmp_d$ to the incoming messages, and $\match$ consists of the edges $(uv)$ with $\si_{u\to v}=\si_{v\to u}=\one$. We refer to $\usi$ as the \emph{auxiliary} configuration, and we now observe that we can define a model on auxiliary configurations on $d$-regular graphs which is in bijection with the frozen model but has the advantage of being a \emph{factor model} in a relatively simple sense.

Recall our definition of the random $d$-regular graph $\Gnd$ as given by a uniformly random matching $M$ on half-edges $H=[nd]$ incident to vertices $V=[n]$. For convenience, we now bisect each edge in $\Gnd$ by a new \emph{clause} vertex $a$, and refer to the resulting graph as the \emph{$(d,2)$-regular bipartite factor graph}: this graph has vertex set $V\cup F$ with bipartition into the set $V$ of \emph{variables} (vertices in the original graph) and the set $F$ of \emph{clauses} (edges in the original graph). The new graph has edge set $E$ where $(ia)\in E$ ($i\in V$, $a\in F$) indicates that in the original graph vertex $i$ is incident to edge $a$. These edges are \emph{labelled}, thus the new bipartite graph is simply equivalent to the original graph $G=(V,M,H)$ together with a labelling on $M$ as well as an ordering within each pair formed by $M$: this contributes a factor $(nd/2)!2^{nd/2}$ to the enumeration but clearly the problem remains unchanged, so we shall continue to use the notation $\Gnd$ for the $(d,2)$-regular bipartite factor graph.

We shall define a model on auxiliary configurations $\usi$ on the bipartite factor graph where $\si_{u\to a} = \si_{a\to v} \equiv \si_{u\to v}$ for variables $u,v\in V$ joined by clause $a\in F$ (that is, the clause acts trivially by passing on the incoming message). The spins of the auxiliary model are the bidirectional messages $\si_{(va)}\equiv\si_{(av)}\equiv(\si_{v\to a},\si_{a\to v})$, taking values in the alphabet $\msg\equiv\set{\zro,\one,\free}^2$. Write $\dsi_v$ for the $d$-tuple of spins on the edges incident to variable $v\in V$, and write $\ksi_a$ for the pair of spins on the edges incident to clause $a\in F$. We hereafter refer to the \emph{intensity} $\ii(\usi)$ of $\usi$ to mean the intensity of its corresponding frozen model configuration,
$\ii(\usi)
\equiv|\set{v:\dmp_d(\usi_v)=\one}|
+\tf12|\set{v:\dmp_d(\usi_v)=\free}|$.

The weight of configuration $\usi\in\msg^E$ under the \emph{auxiliary model} is given by
$$ 
\Psi^\lm(\usi)
\equiv
\prod_{v\in V} \dpsi^\lm(\dsi_v)
	\prod_{a\in F}\hpsi^\lm(\ksi_a)
\equiv
\lm^{\ii(\usi)}
	\prod_{v\in V}\dpsi(\dsi_v)
	\prod_{a\in F}\hpsi(\ksi_a)
\equiv
\lm^{\ii(\usi)}\Psi(\usi)
$$
for factors defined as follows: in the unweighted version, the variable factor weight $\dpsi(\dsi_v)$  is simply the indicator that each outgoing message $\si_{v\to a}$ is determined by the message-passing rule $\dmp_{d-1}$ from the incoming messages $\si_{b\to v}$, $b\in\pd v\setminus a$. We define the clause factor $\hpsi$ analogously, so (recalling that $\si$ denotes pair consisting of the variable-to-clause message followed by the clause-to-variable message) $\hpsi(\si,\acute\si)=\Ind{\si=\refl\acute\si}$ where $\refl$ is the reflection map on $\msg$, $\refl:\eta\eta'\mapsto\eta'\eta$. The weighted version gives a factor $\lm$ to vertex configurations $\dmp_d(\dsi_{\pd v\to v})=\one$ and to edge configurations $(\si,\si')=(\oo,\oo)$:
\beq\label{e:psi}
\text{\footnotesize$\DS
\begin{array}{l}
\hpsi^\lm(\si,\si')
=\hpsi(\si,\si') \lm^{\Ind{\si=\oo}}
=\Ind{\si'=\refl\si} \lm^{\Ind{\si=\oo}},\\
\dpsi^\lm(\dsi)
=\dpsi(\dsi) \lm^{\Ind{ \eta(\dsi)=\one }}
=\begin{cases}
\lm,&\dsi\in
\mathrm{Per}[(\oz^j,\of^{d-j})_{0\le j\le d}]
	\quad\text{($\one$-variable)},\\
1,&\dsi\in
\mathrm{Per}[(\oo,\fz^j,\ff^{d-1-j})_{0\le j\le d-1}]
	\quad\text{($\free$-variable)},\\
1,&\dsi\in
\mathrm{Per}[(\fo^2,\zz^j,\zf^{d-2-j})_{0\le j\le d-2}]
	\quad\text{($\zro$-variable, \emph{susceptible})},\\
1,&\dsi\in
	\mathrm{Per}[(\zo^k,\zz^j,\zf^{d-k-j})_{0\le j\le d-k,
		3\le k\le d}]
	\quad\text{($\zro$-variable, \emph{robust})},
\end{cases}
\end{array}$}
\eeq
where we have classified a variable with spin $\zro$ in the frozen model as \emph{susceptible} (hereafter~$\zros$) if it has exactly two neighboring $\one$'s, \emph{robust} otherwise (hereafter $\zror$).

\newcommand{\msgv}{{}^\mathrm{v}\hspace{-6pt}\msg}
\newcommand{\phiv}{{}^\mathrm{v}\hspace{-2pt}\phi}

\brmk\label{r:vertex.aux}
The frozen model is in exact bijection with the auxiliary model.\footnote{The correspondence remains valid even when the graph has multi-edges, provided we count neighbors with edge multiplicity --- e.g.\ if a $\zro$ neighbors a single $\one$-variable via a doubled edge, we consider it as neighboring two distinct $\one$-variables.} Indeed, let $\msgv\equiv\set{\zz,\zo,\oz,\zf,\fz,\ff,\ffp}$ and consider the mapping $\bm{v}:\msg\to\msgv$ which takes $\oo\mapsto\ffp$, $\set{\oz,\of}\mapsto\oz$, $\set{\zo,\fo}\mapsto\zo$, and acts as the identity on the remaining spins: $\bm{v}(\si)$ gives the two frozen model spins which must be incident to an edge with auxiliary model spin $\si$, with $\ffp$ indicating a matched pair of $\free$'s. 
Coordinate-wise application of $\bm{v}$ to an auxiliary configuration $\usi$ produces the configuration $\bm{v}(\usi)$ from which one can directly read the corresponding frozen configuration $\uz$. This mapping is clearly invertible, by changing the $\zo$-messages incident to $\zros$-vertices back to $\fo$. Let us note that the distribution of $\bm{v}(\usi)$ is also a factor model, with factors $\phiv$ obtained by applying $\phi$ to the pre-image under $\bm{v}$; we refer to this as the \emph{vertex-auxiliary} model.
\ermk

The primary purpose of defining the auxiliary model is that it gives us the following approach for calculating $\E\ZZ_{n\al}$. Given an auxiliary configuration $\usi$, write $m\equiv nd/2$ and consider the normalized empirical measures
\[\begin{array}{llr}
\dbh(\dsi)
	\equiv n^{-1}\sum_{v\in V} \Ind{\dsi_v=\dsi}
	& (\dsi\in\msg^d) \quad
	& \text{variable empirical measure;}\\
\hbh(\ksi) \equiv m^{-1}\sum_{a\in F}\Ind{\ksi_a=\ksi}
	& (\ksi\in\msg^2) \quad
	& \text{clause empirical measure.}
\end{array}\]
We regard $\bh\equiv(\dbh,\hbh)$ as a vector indexed by $(\supp\dpsi,\supp\hpsi)$. For $\si\in\msg$ and $\dsi\in\supp\dpsi$ let $\dH_{\si,\dsi}$ denote the number of appearances of $\si$ in $\dsi$, and similarly write $\hH_{\si,\ksi}$ for the number of appearances of $\si$ in $\ksi$. Then $\dH\in\Z^{\bar{s}\times\bm{\DOT{s}}}$, $\hH\in\Z^{\bar{s}\times\bar{s}}$
where $\bm{\DOT{s}}\equiv|\supp\dpsi|$, $\bar{s}\equiv|\supp\hpsi|=|\msg|$; and for $\bh$ to correspond to a valid configuration $\usi$,
the variable and clause empirical measures must give rise to the same edge marginals
$$\TS
\vh = d^{-1} \dH\dbh = 2^{-1} \hH\hbh,\quad
\vh(\si)\equiv (nd)^{-1}\sum_{(va)\in H}\Ind{\si_{(va)}=\si}.
$$

\bdfn\label{d:simplex}
Let $\simplex$ denote the space of probability measures $\bh\equiv(\dbh,\hbh)$ on $\supp\vph$ (i.e., $\dbh$ is a probability measure on $\supp\dpsi$ while $\hbh$ is a probability measure on $\supp\hpsi$) such that
\bnm[(i)]
\item $(\dbh,\tf{d}{2}\hbh)$ lies in the kernel of matrix
	$H_{\simplex}\equiv \bpm \dH & -\hH \epm
	\in \Z^{ \bar{s}\times( \dbs+\bar{s} ) }$, and
\item $\dbh(\dmp_d(\dsi)=\free)\le\bemax$ (cf.~\eqref{e:truncated.frozen.model}).
\enm
Let $\simplexal$ denote the subspace of measures $\bh\in\simplex$ with normalized intensity
$$\ii(\bh)
\equiv \dbh( \eta(\dsi)=\one )+ \tf12\dbh(\eta(\dsi)=\free)
= \vh(\og)+ (d/2)\,\vh(\oo)=\al.$$
We shall show
(Lem.~\ref{l:lin.bij})
that $H_{\simplex}$ is surjective, therefore $\simplex$ is an $(\bm{\DOT{s}}-1)$-dimensional space with $\simplexal$ an $(\bm{\DOT{s}}-2)$-dimensional subspace.
\edfn

The expected number of auxiliary configurations on $\Gnd$ with empirical measure $\bh$ is
$$\E\ZZ(\bh)
=
\f{\binom{n}{n\dbh}
\binom{m}{m\hbh}}
	{ \binom{nd}{ nd\vh } }
\dpsi^{n\dbh}
\hpsi^{m\hbh}
\equiv
\DSf{n! m!}{(nd)! / 
	\big\{ \prod_\si ( nd\vh(\si) )!\big\} }
\prod_{\dsi} \DSf{\dpsi(\dsi)^{n\dbh(\dsi)} }{(n\dbh(\dsi))!}
\prod_{\ksi} \DSf{\hpsi(\ksi)^{m\hbh(\ksi)} }{(m\hbh(\ksi))!}.
$$
Stirling's formula gives
$\E\ZZ(\bh)= n^{O(1)} \exp\{n\bPhi(\bh)\}$ where
\beq\label{e:bethe}
\bPhi(\bh)
\equiv
\sum_{\dsi}\dbh(\dsi)\log\DSf{\dpsi(\dsi)}{\dbh(\dsi)}
+ (d/2)
\sum_{\ksi}\hbh(\ksi)\log\DSf{\hpsi(\ksi)}{\hbh(\ksi)}
-d  \sum_\si\vh(\si)\log\DSf{1}{\vh(\si)}.\eeq
If further $\min\bh\gtrsim_d1$ as $n\to\infty$, then
\beq\label{e:poly.correction}
\E\ZZ(\bh) =
\f{
e^{O_d(n^{-1})} \mathscr{P}(\bh)}{(2\pi n)^{(\bm{\DOT{s}}-1)/2}}
 \exp\{n\,\bPhi(\bh)\},\quad
\mathscr{P}(\bh)
\equiv
\Big[
	\f{\prod_\si d\vh(\si)}
	{2 \prod_{\dsi} \DS \dbh(\dsi) \TS
	\prod_{\ksi} \DS (\tf12 d\hbh(\ksi) )
	}  \Big]^{1/2}.
\eeq
Clearly an analogous expansion holds for the expectation of the
$\lm$-weighted partition function $\ZZ^\lm(\bh) =  \lm^{n\,\ii(\bh)} \ZZ(\bh)$; we write $\bPhi^\lm(\bh) = \bPhi(\bh) + \ii(\bh)\log\lm$ for the associated rate function (with $\vph^\lm$ in place of $\vph$).

The first moment of frozen model configurations at intensity $\al$ is 
$\E\ZZ_{n\al}=\sum_{\bh\in\simplex_{\al}} \E\ZZ(\bh)$.\footnote{Though we omit it from the notation, the sum should be taken over measures $\bh$ such that $n\dbh$ and $\tf{nd}{2}\hbh$ are integer-valued.} The aim of this section is to compute the exponent $\bPhistar(\al)=\lim_n n^{-1}\log\E\ZZ_{n\al}$ by determining the maximizer $\bhstaral\equiv(\dbhstaral,\hbhstaral)$ of $\bPhi$ on $\simplexal$. Observe it is clear from the functional form of $\bPhi$ that $\dbhstaral$ and $\hbhstaral$ must be symmetric functions on $\msg^d$ and $\msg^2$ respectively.

\emph{The fugacity parameter $\lm$ serves the purpose of a Lagrange multiplier}: if $\bh\in\simplexal$ is a stationary point of $\bPhi$ restricted to $\simplexal$, then for some $\lm$ it must be a stationary point of $\bPhi^\lm$ on the \emph{unrestricted} space $\simplex$.
Strictly positive measures $\bh$ which are stationary for $\bPhi^\lm$ as a function on $\simplex$ (unrestricted) correspond to a generalization of the tree Gibbs measures considered in \S\ref{ss:frozen.tree}, where the boundary conditions are specified by a law on incoming \emph{and outgoing} messages, as we now describe. Let $\treebip$ denote the infinite tree given by bisecting each edge of $\tree$ by a new (clause) vertex; this is the local weak limit of the random $(d,2)$-regular bipartite factor graph. The vertices at level $t$ of $\treebip$ are variables for $t$ odd, clauses for $t$ even. If $\usi$ is a message configuration on the edges of $\treebip(t)$ --- including the edges $E(t-1,t)$ joining levels $t-1$ and $t$ --- then let $\Psi^\lm_t(\usi)$ denote the product of the factor weights $\dpsi^\lm(\dsi_v)$, $\hpsi^\lm(\ksi_a)$ over all $v,a\in \treebip(t-1)$. For probability measures $\hd,\hh$ on $\msg$ we define the measures
\beq\label{e:gibbs.aux}
\begin{array}{rll}
Z_t \nuaux_t(\usi)
	\hspace{-6pt}&= \Psi^\lm_t(\usi)
		\prod_{e\in E(t-1,t)} \hd_{\si_e}
	&\text{($t$ even)},\\
Z_t \nuaux_t(\usi)
	\hspace{-6pt}&= \Psi^\lm_t(\usi)
		\prod_{e\in E(t-1,t)} \hh_{\si_e}
	&\text{($t$ odd)},
\end{array}
\eeq
with $Z_t$ the normalizing constant which makes $\nuaux_t$ a probability measure.\footnote{We suppress the $\lm$-dependence from the notation except to differentiate $\vph^\lm$ from $\vph\equiv\vph^\lm|_{\lm=1}$.} The family $(\nuaux_t)_t$ is consistent if and only if $h\equiv(\hd,\hh)$ satisfies the \emph{Bethe recursions}
\beq\label{e:bp}
\begin{array}{rlr}
\dz\hd_\si
	\hspace{-6pt}&= 
	\sum_{\dsi\,:\,\si_1=\si}
	\dpsi^\lm(\dsi) \prod_{i=2}^d \hh_{\si_i},
	\quad&\text{variable Bethe recursions;}\\
\hz\hh_\si
	\hspace{-6pt}&= 
	\sum_{\ksi\,:\,\si_1=\si}
	\hpsi^\lm(\ksi) \hd_{\si_2}
	\quad&\text{clause Bethe recursions.}
\end{array}
\eeq
(with $\dz,\hz$ the normalizing constants). We show below that (as expected from our definition) these recursions are a generalization of the frozen model recursions \eqref{e:frozen.recursions}. Thus a solution $h$ of \eqref{e:bp} specifies a Gibbs measure $\nuaux\equiv\nuaux^\lm$ for the auxiliary model on $\treebip$ which generalizes the measures $\nu\equiv\nu^\lm$ described in \S\ref{ss:frozen.tree}.

The Bethe recursions read explicitly as follows. Write $\forz\equiv\set{\zro,\free}$, e.g.\ $\hh_\og\equiv\hh_\oz+\hh_\of$. The clause Bethe recursions are simply
\beq\label{e:clause.bp.explicit}
\hz\hh_\si=\hd_{\refl\si}\lm^{\Ind{\si=\oo}}
\quad\text{with }
	\hz = 1+(\lm-1)\hd_\oo.
\eeq
The variable Bethe recursions are
$$\begin{array}{rl}
\text{({\sc a})} & \DS
	\dz\hd_\oz=\dz\hd_\of=\lm(\hh_\og)^{d-1},\quad\quad
	\lm\dz\hd_\oo=\lm (\hh_\fg)^{d-1},
	\vspace{2pt}\\
& \DS \dz\hd_\fz=\dz\hd_\ff
	=(d-1)\hh_\oo (\hh_\fg)^{d-2},\quad\quad
	\dz\hd_\fo
	=(d-1)\hh_\fo (\hh_\zg)^{d-2},
	\vspace{2pt}\\
\text{({\sc b})} & \DS
	\dz\hd_\zo
	=(\hh_\zd)^{d-1}
	-(\hh_\zg)^{d-1}-(d-1)\hh_\zo(\hh_\zg)^{d-2},
	\vspace{2pt}\\
\text{({\sc c})} & \DS
	\dz\hd_\zz=\dz\hd_\zf
	=(\hh_\zd)^{d-1}-(\hh_\zg)^{d-1}-(d-1)\hh_\zo(\hh_\zg)^{d-2}\\
	&\qquad\qquad\qquad\DS +
		\tf12 (d-1)(d-2) [(\hh_\fo)^2-(\hh_\zo)^2](\hh_\zg)^{d-3}
\end{array}$$
We immediately have
$\hh_\zo=\hh_\fo$, $\hh_\zf=\hh_\ff$, and $\hh_\zz=\hh_\fz$, and then comparing  ({\sc b}) and ({\sc c}) gives $\hh_\oz=\hh_\zz=\hh_\fz$. It then follows from ({\sc a}) that the following are equivalent (with the symbol {\tiny$\checkmark$} indicating the identities we already know):
\beq\label{e:bp.sym}
\hd_\fo=\hd_\fz\stackrel{\checkmark}{=}\hd_\ff,\quad
\hh_\of=\hh_\zf\stackrel{\checkmark}{=}\hh_\ff,\quad
\lm\hd_\oo=\hd_\oz\stackrel{\checkmark}{=}\hd_\of,\quad
\hh_\oo=\hh_\zo\stackrel{\checkmark}{=}\hh_\fo.
\eeq
\emph{If \eqref{e:bp.sym} holds, then
 \eqref{e:bp} reduces to the frozen model recursions \eqref{e:frozen.recursions} with }
\beq\label{e:hh.q}
\begin{array}{rl}
3\,\hh_{\eta'\eta}
	\hspace{-6pt}&= 
	q_\eta,\quad
	\text{therefore (substituting into \eqref{e:clause.bp.explicit})}\\
3\,\hd_{\eta\eta'}
	\hspace{-6pt}&= 
	\lm^{-\Ind{\eta\eta'=\oo}} q_\eta \hz
	\quad\text{with }
	\hz=[1-\tf13 \qo(1-1/\lm)]^{-1}
\end{array}
\eeq
proving our claim that the measures $\nuaux$ generalize the measures $\nu$ of \S\ref{ss:frozen.tree}. The connection between these Gibbs measures and the rate function $\bPhi^\lm$ is given by the following variational principle:

\blem\label{l:interior.bp}
If a measure $\bh$ in the interior $\simplex^\circ$ of $\simplex$ is stationary for $\bPhi^\lm$, then $\bh$ corresponds to a solution $h\equiv h^\lm$ of the Bethe recursions \eqref{e:bp} via
\beq\label{e:bij}\TS\dbz\dbh(\dsi)
	=\dpsi^\lm(\dsi)
	\prod_{i=1}^d \hh_{\si_i},\quad
\hbz\hbh(\ksi)
	=\hpsi^\lm(\ksi)
	\prod_{i=1}^2 \hd_{\si_i},\quad
\bar{z}\vh(\si)
	=\hh_\si\hd_\si
\eeq
with $\dbz,\hbz,\bar{z}$ normalizing constants
satisfying $\bar{z}=\dbz/\dot z=\hbz/\hz$ for $\dz,\hz$ as in \eqref{e:bp}.

\bpf
Follows from \cite{dss-naesat} upon verifying that $\dH,\hH$ are surjective; we will prove a stronger condition than surjectivity in Lem.~\ref{l:lin.bij} and \eqref{e:surj}.
\epf
\elem

Recall now that our aim is to locate the global maximizer of $\bPhi$ on $\simplexal$ as a stationary point of $\bPhi^\lm$ for some value of $\lm$. 

\bthm\label{t:first.moment.exponent}
Let $\bhstarlm$ denote the unique stationary point of $\bPhi^\lm$ which corresponds via \eqref{e:bij} to the only solution $h\equiv h^\lm$ of the Bethe recursions \eqref{e:bp} which satisfies \eqref{e:bp.sym}, that is, $\hh^\lm_{\eta'\eta}=\tf13 q^\lm_\eta$. The unique maximizer $\bhstaral$ of $\bPhi$ on $\simplexal$ is given by $\bhstarlm$ for the unique value of $\lm$ at which $\bhstarlm(\dmp_d(\dsi)=\one)=\al$.
\ethm

In view of Lem.~\ref{l:interior.bp} and our preceding discussion of Gibbs measures and Lagrange multipliers, Thm.~\ref{t:first.moment.exponent} will follow by showing
\bnm[1.]
\item Any global maximizer $\bhstaral$ of $\bPhi$ on $\simplexal$ lies in the interior $\simplexalint$. Thus for some\footnote{Eventually we will find $\lm = \lm_\star = \exp\{ -(\bPhistar)'(\al) \}$.} $\lm$ it is a stationary point of $\bPhi^\lm$, and hence corresponds via \eqref{e:bij} to a solution $h\equiv h_\al$ of the Bethe recursions \eqref{e:bp}.

\item Any such Bethe solution $h\equiv h_\al$ satisfies the symmetries \eqref{e:bp.sym}, and hence is identified with the unique solution of the frozen model recursions \eqref{e:frozen.recursions}. Consequently $\lm$ is given by the unique solution of $\bhstarlm(\dmp_d(\dsi)=\one)=\al$.
\enm
Ruling out boundary maximizers for $\bPhi$ is relatively easy, so we defer the proof to \S\ref{s:second} where we will use the same argument to rule out boundary maximizers for the second-moment exponent $\bPhit$. We turn now to the more delicate task of proving the symmetries \eqref{e:bp.sym}.

\subsection{Bethe recursion symmetries}\label{ss:bp.sym}

\begin{figure}[ht]
\centering
\includegraphics[height=2.5in,trim=.7in .8in 1.4in .9in,clip]{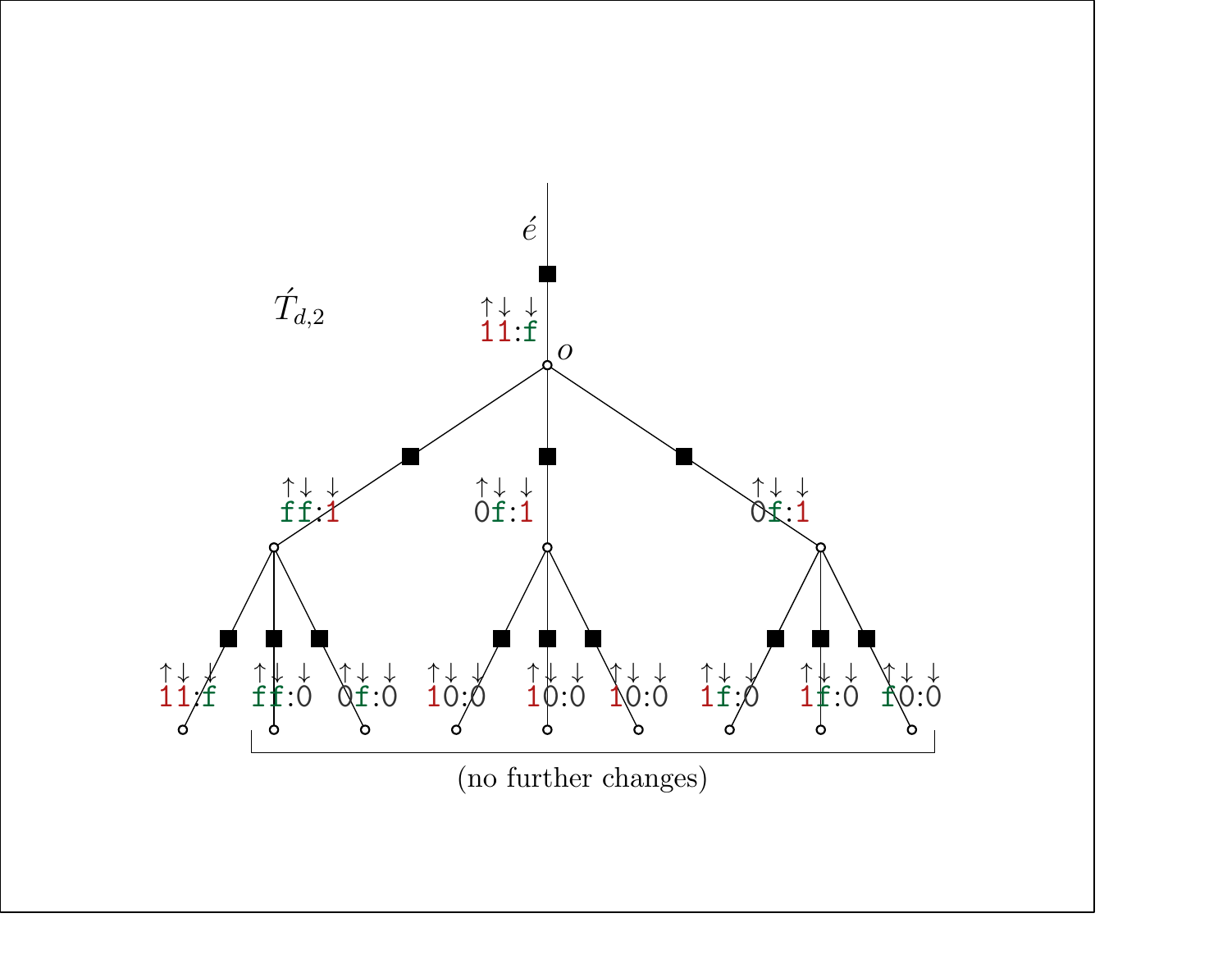}
\caption{Change of message incoming down to $\acute{e}$ is passed down
	$\etreebip$\\
	($\eta\acute\eta$:$\grave\eta$
	means message $\eta$ up, message $\acute\eta$ down in $\usi$,
	message $\grave\eta$ down in $\usi'$)
}
\label{f:bp.sym}
\end{figure}

Suppose $\bhstaral$ is an interior maximizer for $\bPhi$ on $\simplexal$, and so corresponds to a Bethe solution $h\equiv h^\lm$. Let $\etreebip$ denote $\treebip$ with a subtree incident to the root removed, such that one clause $a\in\pd\rt$ is incident to an unmatched half-edge $\acute{e}$ (Fig.~\ref{f:bp.sym}). Consider defining a Gibbs measure for $\phi^\lm$ on $\etreebip$ in the manner of \eqref{e:gibbs.aux}, with boundary law given by $h$. Then the marginal law of $\si_{\acute{e}}$ will be $\hh$, and the marginal law of the $d$-tuple of spins incident to any given vertex will be $\dbh$. Further, the Gibbs measure on $\etreebip$ can be generated in Markovian fashion, starting with spin $\si_{\acute{e}}$ distributed according to $\hh$, generating the messages on the other $d-1$ edges incident to $\rt$ according to the conditional measure $\dbh(\dsi\,|\,\si_1=\si_{\acute e})$, and continuing iteratively down the tree.

Write $\si_{\acute{e}} \equiv \bm{io}$ where $\bm{i}$ is the incoming message and $\bm{o}$ is the outgoing (in Fig.~\ref{f:bp.sym}, $\bm{o}$ is directed upwards, $\bm{i}$ downwards). Given any valid auxiliary configuration $\usi$ on the edges of $\etreebip$, changing $\bm{i}$ and passing the changed message through the tree (via $\dmp_{d-1}$) produces a new auxiliary configuration $\usi'$ (Fig.~\ref{f:bp.sym}). The symmetries \eqref{e:bp.sym} will follow by showing that for any fixed $\bm{o}$, the effect of changing $\bm{i}$ is \emph{measure-preserving} under the Gibbs measure $\nuaux_\al$ corresponding to $\bhstaral$. From our definition of the Gibbs measure via the boundary law, the measure-preserving property will follow by showing that the effect of changing $\bm{i}$ almost surely does not percolate down the tree.

Indeed, recall that we already saw directly from the Bethe recursions that $\hh_{\zro\eta}=\hh_{\free\eta}$ for any $\eta$: this corresponds to the fact that $\dmp_{d-1}$ does not differentiate between incoming messages $\zro$ and $\free$, so changing $\bm{i}$ from $\zro$ to $\free$ or vice versa has no effect at all below $\acute{e}$. We also found that $\hh_{\eta\zro}$ does not depend on $\eta$: if $\bm{o}=\zro$ then all messages outgoing from the root must be $\zro$ or $\free$, so changing $\bm{i}$ can have an effect at most one level down.

\bpf[Proof of Thm.~\ref{t:first.moment.exponent}]
We shall assume that the maximizer $\bhstaral$ of $\bPhi$ on $\simplexal$ lies in the interior $\simplexalint$, deferring the proof to \S\ref{s:second} (see Propn.~\ref{p:boundary} and Cor.~\ref{c:boundary}).

From the above discussion it remains to show the corresponding Bethe solution $h\equiv h_\al$ satisfies $\hh_\oo=\hh_\fo$: meaning that in the Gibbs measure $\nuaux_\al$ corresponding to $\bhstaral$, changing $\bm{io}$ from $\oo$ to $\fo$ or vice versa has a finite-range effect. Fig.~\ref{f:bp.sym} shows that the effect of changing $\oo$ to $\fo$ can only propagate through components of $\free$-variables, while the effect of changing $\fo$ to $\oo$ can only propagate through chains of alternating $\one$- and $\zros$-variables (cf.~\eqref{e:psi}).\footnote{Here we write $\eta$-variables to mean variables with spin $\eta$ in the corresponding frozen model configuration.}

We claim that both propagations are subcritical under $\nuaux_\al$. By \eqref{e:bij}, for any $2\le j\le d$ the density of $\zro$-variables with exactly $j$ neighboring $\one$-variables is, up to a normalizing constant, $\tbinom{d}{j}\hd_{\oz}^j\hd_{\gz}^{d-j}$ (where we have used that $\hd_\oz=\hd_\of$). Thus, among all the edges incident to $\zro$-variables, the proportion coming from $\one$-variables is
\beq\label{e:derivative.conditioned.binomial}
\f{\aone}{1-\aone-\be}
=
\f{\tf1d \sum_{j\ge2} j \binom{d}{j} u^j }
	{\sum_{j\ge2} \binom{d}{j} u^j }
\quad\text{where }u\equiv \hd_\oz/\hd_{\gz}.
\eeq
This is an increasing function of $u$ (the derivative in $u$ is the variance of a conditioned binomial random variable), so
in the regime $\albd\le\al\le\aubd$ we must have
\beq\label{e:q.of.alpha}
	\hd_\oz/\hd_{\cdot\zro}
	= \al\,[1+O(\tf{\log d}{d})].
\eeq
Comparing the number of $\one$-to-$\zros$ edges with the number of other $\one$-to-$\zro$ edges gives (applying \eqref{e:bij} and \eqref{e:q.of.alpha})
\beq\label{e:fo.zo.ratio}
\DSf{\hd_\fo}{\hd_\zo}
= \f{nd\hbh(\of,\fo)}{nd\hbh(\oz,\zo)}
= \f{2n\dbh(\eta=\zros)}{nd\al-2n\dbh(\eta=\zros)}
=
\DSf{2\binom{d}{2} u^2}
{ \sum_{j\ge3} j \binom{d}{j} u^j }
\le d^{-3/2}.
\eeq
Thus, conditioned on $\bm{io}=\fo$, the expected number of $a\in\pd\rt$ with $\si_{a\to\rt}=\free$ is less than one, from which it follows that the effect of changing $\fo$ to $\oo$ does not percolate.\footnote{Note the main estimate we required for this argument as that $\dbh(\eta=\zros)\ll d^{-1}$, which is intuitively clear and can be derived by other means. We have taken an approach which circumvents further {\it a priori} estimates.} Likewise
$$\f{\be}{1-\al-\be}
\ge
\f{\tf1d\sum_{j\ge2} (d-j) (\hd_\fz/\hd_\gz)
	\binom{d}{j} u^j}
	{\sum_{j\ge2} \binom{d}{j} u^j}
=\hd_\fz/\hd_\gz\,[1- O(\logdbyd)] 
$$
so $\hd_\fz/\hd_\zz \lesssim d^{-3/2}$ which implies that the effect of changing $\oo$ to $\fo$ does not percolate, concluding the proof.
\epf

\subsection{Explicit form of first moment exponent}\label{ss:explicit}

We conclude this section by giving the explicit form of $\bPhistar(\al)\equiv\bPhistar_d(\al)$ when $\albd\le\al\le\aubd$:

\bppn\label{p:explicit}
For $d\ge d_0$, $\albd\le\al\le\aubd$,
$\bPhistar(\al)\equiv \bPhistar_d(\al)$
is given by
\beq\label{e:explicit}
\bPhistar(\al)
=-\log[1-q(1-1/\lm)]
	-(d/2-1)\log[1-q^2(1-1/\lm)]
	-\al\log\lm
\eeq
where $\lm$ and $q$ are determined from $\al$ via
\beq\label{e:q.lambda.alpha}
\left\{
\begin{array}{l}
\lm= q\smf{1-(1-q)^{d-1}}{(1-q)^d},\text{ and }
q\equiv x\logdbyd\text{ is the unique solution of }\\
\al= q \smf{1+(d/2-1)  (1-q)^{d-1}}{ 1+q-(1-q)^{d-1} }
\text{ on the interval }
1.6\le x\le3.
\end{array}\right.
\eeq
On the interval $\albd\le\al\le\aubd$ the function $\bPhistar$ is strictly decreasing, with a unique zero $\al_\star$ in the interval's interior. The gap between $\al_\star$ and the first moment threshold $\alfm$ of the original independent set partition function \eqref{e:alpha0} is given by
\beq\label{e:threshold.gap}
\alfm-\al_\star
	= \Big(\smf{e\log d}{2 d}\Big)^2
		\Big[ 1 + O\Big( \smf{1}{\log d}\Big) \Big].
\eeq

\bpf
In the following, let $q\equiv q^\lm$ denote the solution $(\qo)^\lm$ of the recursions \eqref{e:frozen.recursions}.

\smallskip\noindent\emph{Explicit Bethe prediction.}
Substituting \eqref{e:bij} into \eqref{e:bethe} and rearranging gives 
\[
\begin{array}{rl}
\bPhistar^\lm=\bPhi^\lm(\bhstarlm)
	\hspace{-6pt}&= \log\dbz+(d/2)\log\hbz
	- d \sum_\si \vh(\si)
		\log[ \hd_\si\hh_\si
			/\vh(\si) ]\\
	&=\log\dbz
	+(d/2)\log\hbz - d\log\bar{z}.
\end{array}
\]
We use \eqref{e:hh.q} to calculate $\bar{z},\dbz,\hbz$ in terms of $q,\lm$:
\begin{align*}
\bar{z}
&= \hz \sum_\si
	\smf{\hh_\si \hh_{\refl\si}}
	{\lm^{\Ind{\si=\oo} }}
= \hz\sum_{\eta,\eta'}
	\smf{(q_\eta/3)(q_{\eta'}/3)}
	{\lm^{\Ind{\si=\oo}}}
=\smf{\hz}{9}\,
	[1-q^2(1-1/\lm)];\\
\hbz
&= \sum_\si
	\lm^{\Ind{\si=\oo}} \hd_\si \hd_{\refl\si}
= \hz^2 \sum_\si
	\smf{\hh_\si \hh_{\refl\si}}
		{\lm^{\Ind{\si=\oo}}}
= \bar{z}\hz;\\
\dbz
&= \lm (\hh_\og)^d
	+ d \hh_\oo(\hh_\fg)^{d-1}
	+ \smb{d}{2} (\hh_\fo)^2 (\hh_\zg)^{d-2}
	+ \sum_{k\ge3}
	(\hh_\zo)^k (\hh_\zg)^{d-k}\\
&= \smf{1}{3^d}
	[ 1 + (\lm-1)(1-q)^d]
= \smf{ 1-q^2(1-1/\lm) }{3^d[ 1-q(1-1/\lm) ]}
=\bar{z} \dz.
\end{align*}
Substituting back into \eqref{e:bij} gives $\dbh(\dmp_d(\dsi)=\one) =\dbz^{-1} \lm(\hh_\oz)^d =q(1-q)/[1-q^2(1-1/\lm)]$ and $\dbh(\dmp_d(\dsi)=\free) = \dbz^{-1} d\hh_\oo (\hh_\fg)^{d-1} = dq/[\lm(1-q)]\dbh(\dmp_d(\dsi)=\one)$, therefore
\[
\al=\dbh(\dmp_d(\dsi)=\one)
	+\tf12 \dbh(\dmp_d(\dsi)=\free)
	=q\smf{1-q+ dq/(2\lm)}{1-q^2(1-1/\lm)}.
\]
The recursion \eqref{e:frozen.recursions} also gives the expressions in \eqref{e:q.lambda.alpha} for $\lm$ and $\al$ solely in terms of $q$. The mapping $q\mapsto\al$ is not one-to-one on the entire interval $0\le q\le 1$, but recalling \eqref{e:q.of.alpha} we must have $q=x\logdbyd$ for $1.6 \le x\le 3$, and on this interval it is easily verified that the mapping is indeed one-to-one, with
$\lm = d^x q [ 1+O( d^{-1}(\log d)^2 ) ]$
and $D_q\,\al = 1 + O(d^{-1/2})$
where $D_q$ denotes the total derivative with respect to $q$.
This completes the verification of \eqref{e:q.lambda.alpha};
it then follows from Thm.~\ref{t:first.moment.exponent}
that $\bPhistar(\al)
=\bPhistar^\lm(\bhstarlm)-\al\log\lm$
is given by \eqref{e:explicit}. We note here that
$\bPhistar'(\al)
=-\log\lm$, therefore
$\bPhistar''(\al)
= - ( D_q\al )^{-1} D_q(\log\lm)
=-d[1+O((\log d)^{-1})]$.

\newcommand{\qis}{\acute{q}}
\newcommand{\alis}{\acute{\al}}
\newcommand{\lmis}{\acute{\lm}}
\smallskip\noindent\emph{Comparison of first-moment exponents.}
By contrast, the original independent set partition function has first-moment exponent $\Phi(\al)$ calculated in \eqref{e:first.moment}.
This exponent also has a Bethe variational characterization,
which can be expressed in terms of the fixed point
$\qis$ of the hard-core tree recursions:
\begin{align*}
\Phi(\al)
&= \log[ \lmis(1-\qis)^d+1 ]
	-(d/2) \log(1-\qis^2)-\al\log\lmis\\
&\quad \text{where }
\qis= \smf{\lmis(1-\qis)^{d-1}}{1+\lmis(1-\qis)^{d-1}}
\text{ and }
\al = \smf{\lmis(1-\qis)^d}{1+\lmis(1-\qis)^d}.
\end{align*}
This formula can be derived loosely in the same manner as \eqref{e:explicit}; its validity can be checked simply by verifying that it agrees with \eqref{e:first.moment}. 
The relation between $\al,\qis$ is 
given by $\al = \alis(\qis)= \qis/(1+\qis)$,
$\qis = \alis^{-1}(\al) = \al/(1-\al)$, so we compare $\bPhistar(\al)$ and $\Phi(\al)$ by expressing both in terms of $q$: $\bPhistar(\al)$ is given by \eqref{e:explicit} with $\al,\lm$ defined in terms of $q$ by \eqref{e:q.lambda.alpha}, while
$$\begin{array}{l}
\Phi(\al)
=-\log(1-\qis)
	-(d/2-1)\log(1-\qis^2)
	-\al\log\lmis,\\
\text{where }
	\lmis = \qis / (1-\qis)^d \text{ and }
	\qis = \alis^{-1}( \al(q) ).
\end{array}
$$
Let us emphasize that $q$ and $\qis$ are not equal but rather are related through the same $\al$; explicitly
$\qis
=q(1-q+\tf{dq}{2\lm})/(1-q-\tf{d-2}{2\lm}q^2)$. A little algebra then gives
\begin{align*}
&(\Phi-\bPhi)(\al)
=\log \smf{1-q+q/\lm}{1-q}
	+ \log \smf{1-q}{1-\qis}
	+(d/2-1)\log\smf{1-q^2+q^2/\lm}{1-\qis^2}
	-\al\log\smf{\lmis}{\lm}\\
&= \log \smf{1-q+q/\lm}{1-q}
	+ \log\smf{\DS (1-q)(1-q-\tf{d-2}{2\lm}q^2)}
		{\DS (1-q)^2-\tf{d-1}{\lm}q^2 }
	+(d/2-1)
	\log\smf{ \DS (1-q)^2( 1-\tf{d-2}{2\lm} \tf{q^2}{1-q} )^2 }
		{ \DS (1-q)^2-\tf{d-1}{\lm} q^2 }
	-\al\log\smf{\lmis}{\lm}\\
&= \log\Big( 1 + \smf{q}{\lm(1-q)}\Big)
-\log\Big( 1-\smf{d-2}{2\lm}\smf{q^2}{1-q} \Big)
+d\,\log\f{
	\DS  1-\tf{d-2}{2\lm} \tf{q^2}{1-q}  }
	{ \DS [ 1-\tf{d-1}{\lm}
		\tf{q^2}{(1-q)^2} ]^{1/2} }
	-\al\log\smf{\lmis}{\lm}.
\end{align*}
Taylor expanding (recalling
$q =x \logdbyd$, $\lm = d^x q [ 1+O( d^{-1}(\log d)^2 ) ]$)
gives
\begin{align*}
\al\log\smf{\lmis}{\lm}
&=\al\log\smf{\qis}{q}
+d\al\log\smf{1-q}{1-\qis}
-\overbrace{\al\log[1-(1-q)^d]}^{ O(q^2/\lm) }\\
&=\al\log
	\Big( 1 + \smf{dq}{2\lm}
		\smf{\DS 1 + (1-\tf{2}{d})q }
		{\DS 1-q-\tf{d-2}{2\lm}q^2 } \Big)
+d\al\log
	\Big(
	1+ \smf{dq^2}{2\lm}
		\smf{\DS 1 + (1-\tf2d)q }
		{\DS (1-q)^2 (1-  \tf{(d-1)q^2}{\lm(1-q)^2}) }
		\Big)
		-O\Big(\smf{q^2}{\lm}\Big)\\
&= \smf{dq^2}{2\lm}(dq+1)
	+O\Big( \smf{d^2 q^4}{\lm}
		+  \smf{d^3 q^4}{\lm^2}  \Big)
= \smf{dq^2}{2\lm}(dq+1)
	+O\Big(
	\smf{(\log d)^3}{d^{x+1}}
	+ \smf{(\log d)^2}{d^{2x-1}}
	\Big).
\end{align*}
Substituting into the above expression for
$(\Phi-\bPhi)(\al)$
and expanding the other terms gives
\begin{align}\nonumber
&(\Phi-\bPhi)(\al)
= \smf{q}{\lm}
	+ \smf{dq^2}{2\lm}(dq+2)
	+O\Big( \smf{dq^3}{\lm}
		+\smf{d^2q^3}{\lm^2} \Big)
	-\al\log\smf{\lmis}{\lm}\\
\label{e:exponent.gap}
&=\smf{q}{\lm}
	+ \smf{dq^2}{2\lm}
	+O\Big( \smf{d^2q^4}{\lm}
		+\smf{d^3q^4}{\lm^2} \Big)
= \smf{q}{2\lm}(dq+2)
	+O\Big(
	\smf{(\log d)^3}{d^{x+1}}
	+ \smf{(\log d)^2}{d^{2x-1}}
	\Big)
\asymp d^{-x}\log d.
\end{align}

\newcommand{\qfm}{\qis_\square}
\newcommand{\qsq}{q_\square}
\newcommand{\lmsq}{\lm_\square}

\smallskip\noindent\emph{Comparison of first-moment thresholds.}
It is clear from the above that $\bPhistar$ has a unique zero
$\albd<\al_\star<\alfm$
with $\alfm$ the first moment threshold of the original independent partition function, determined in Lem.~\ref{l:is.first.moment}.
Let $\qsq,\lmsq$ denote the solution of \eqref{e:q.lambda.alpha} for $\al=\alfm$, and note $\qsq=\alfm[1+O( d^{-1}(\log d)^2)]$.
Consider
$0\le\de\le d^{-2}(\log d)^3$:
applying \eqref{e:exponent.gap} with the above estimate of $\bPhistar''(\al)$ gives
\[
\begin{array}{rl}
\bPhistar(\alfm-\de)
\hspace{-6pt}&=\bPhistar(\alfm)+\de\log\lmsq +O(d\de^2)\\
&=-(2\lmsq)^{-1}
	\alfm (d\alfm+2)[ 1+O(d^{-1}(\log d)^3) ]
	+\de\log\lmsq,
\end{array}\]
so the gap between the threshold is given by
$$\alfm-\al_\star
	=[ 1+O(d^{-1}(\log d)^3) ]
	\smf{\alfm(d\alfm+2)}{2\lmsq\log\lmsq}
	\asymp d^{-2}(\log d)^2.$$
For a more precise estimate,
recall from Lem.~\ref{l:is.first.moment}
that $\alfm=\alfmt[1+O(d^{-1}\log d)]$
where
$\alfmt$ is the zero of $F(\al)\equiv\log(e/\al) -
	(d+1)\al/2$.
Then
\[\begin{array}{rl}
\log\lmsq
\hspace{-6pt}&=F(\qsq)
	+(d\qsq/2)+O(d^{-1}(\log d)^2)\\
&=(d\qsq/2)+O( d^{-1}(\log d)^3 )
=(d\alfmt/2)\,[1+O( d^{-1}(\log )^2 )].\end{array}\]
Substituting into the above gives \eqref{e:threshold.gap},
concluding the proof.
\epf
\eppn

\section{Second moment of frozen model}\label{s:second}

In this section we compute the exponential growth rate $\bPhitstar(\al)=\lim_{n\to\infty} n^{-1}\log\E[\ZZ_{n\al}^2]$ of the second moment of the (truncated) frozen model partition function \eqref{e:truncated.frozen.model}. This is done in the same framework as in introduced in \S\ref{s:first}, regarding the second moment as the first moment of the partition function of \emph{pair} frozen model configurations $(\uz^1,\uz^2)$ on the same underlying graph. The corresponding model of pair auxiliary configurations $\uta\equiv(\usi^1,\usi^2)$ has factors $\dpsit\equiv\dpsi\otimes\dpsi$ and $\hpsit\equiv\hpsi\otimes\hpsi$, and rate function $\bPhit$ on the space $\simplextal$ of empirical measures $\bht$ with both marginals in $\simplexal$.

\bthm\label{t:second.moment}
The rate function $\bPhit$ on $\simplextal$ attains its maximum only at the product measure $\bhtstaral\equiv\bhstaral\otimes\bhstaral$
or at the measure $\bhidal$ with marginals $\bhstaral$ which is supported on pair configurations $\uta\equiv(\usi,\usi)$.
\ethm

\subsection{Intermediate overlap regime}

Write $\ZZt_{n\al}(n\rho)$ for the contribution to $(\ZZ_{n\al})^2$ from \emph{pair frozen model} configurations $(\uz^1,\uz^2)$ on $G$ in which exactly $n\rho$ vertices take spin $\one$ in both configurations. In this subsection we show that the value of $\rho$ maximizing $\E[\ZZt_{n\al}(n\rho)]$ must either be very small or very close to $\al$. This will be done by a rather crude comparison with the contribution $\Zt_{n\al}(n\rho)$ to $(Z_{n\al})^2$ from \emph{pair independent set} configurations $(\ux^1,\ux^2)$ on $G$ with overlap $n\rho$. At the same time we shall prove that the global maximizer of $\bPhi$ on $\simplexal$ lies in the interior $\simplexalint$ (as required in the proof of Thm.~\ref{t:first.moment.exponent}), and that any global maximizer of $\bPhit$ on $\simplextal$ \emph{away from} $\bhidal$ lies in the interior $\simplextalint$. More refined estimates for the regimes near $\bhtstaral$ and $\bhidal$ are reserved for \S\ref{ss:indep}~and~\ref{ss:rigid}.

Given a pair frozen configuration $\uz^i\equiv(\ueta^i,\match^i)$ ($i=1,2$) on $G$, define
\beq\label{e:frozen.pair.configuration}
\uom\equiv\uom(\uz^1,\uz^2)\in\pairs^V,\quad
\pairs\equiv\set{\zro,\one,\free}^2\cup\set{\ffne}
\eeq
by setting $\om_v=(\eta^1_v,\eta^2_v)$ \emph{unless} $(\eta^1_v,\eta^2_v)=\ff$ and $v$ is matched to different vertices under $\match^1$ and $\match^2$, in which case we set $\om_v=\ffne$. For $\eta\in\set{\zro,\one,\free}$ write $\pairs_\eta$ for the spins $\om\in\pairs$ with either coordinate equal to $\eta$. For $S\subseteq V$ write
\beq\label{e:subset.spin}
\begin{array}{rrlrl}
\text{(for $\om\in\pairs$)}
	& S_\om
	\hspace{-6pt}&\equiv
	\set{v\in S:\om_v=\om},
	& n_\om
	\hspace{-6pt}&\equiv |V_\om|;\\
\text{(for $\eta\in\set{\zro,\one,\free}$)}
	&
	S_\eta
	\hspace{-6pt}&\equiv
	\set{v\in S:\om_v\in\pairs_\eta},
	& n_\eta
	\hspace{-6pt}&\equiv
	|V_\eta|.
\end{array}
\eeq
Let $\Vumzf$ denote the set of $\zf$-vertices whose matched partner does not have spin $\of$, and $\numzf\equiv |\Vumzf| = n_\zf-n_\of$. Define symmetrically $\Vumfz$ and $\numfz\equiv|\Vumfz|$.
We decompose
$$
\E[\Zt_{n\al}(n\rho)]
=\sum_{\pi\in M(n\al,n\rho)}
	\E[\Zt(\pi)]
\quad\text{and}\quad
\E[\ZZt_{n\al}(n\rho)]
=\sum_{\zpi\in\cM(n\al,n\rho)}
	\E[\ZZt(\zpi)]$$
where $M(n\al,n\rho)$ is the appropriate space of empirical measures on $\set{\zro,\one}^2$ and $\cM(n\al,n\rho)$ is the appropriate space of empirical measures on $\pairs$.

\bppn\label{p:ap.compare.frozen.is}
Suppose $\albd\le\al\le\aubd$ and $0 \le \rho\le\al-\erho\logdbyd$ for $\erho$ a small positive constant uniform in $d$. For $d\ge d_0(\erho)$, $n\ge n_0(d)$,
$$\E[\ZZt_{n\al}(n\rho)] =\E[\Zt_{n\al}(n\rho)] \exp\{ O[ nd^{-1}(d^{-0.49}+d^{-3\erho/4} ) ] \},$$
(recalling that $\ZZ$ refers to the frozen model while $Z$ refers to the independent set model).

\bpf
Let $\pi\in\cM(n\al,n\rho)$;
we decompose the expected 
contribution from
$\pi$ to the pair frozen model partition function as
\beq\label{e:def.J.pi.spin}
\E[\ZZt(\pi)]= \mathbf{c}(\pi)\,\E[\mathbf{J}^\pi_\free \mathbf{J}^\pi_\one \mathbf{J}^\pi_\zro]
\eeq
where $\mathbf{c}(\pi)$ denotes the multinomial coefficient $\tbinom{n}{n\pi}$, and the $\mathbf{J}^\pi_\eta$ are defined with respect to a \emph{fixed} spin configuration $\uom\in\pairs^V$ with empirical measure $\pi$, as follows: $\mathbf{J}^\pi_\free$ is the number of matchings on the $\free$-vertices; $\mathbf{J}^\pi_\one$ is the indicator that $\one$-vertices neighbor only $\zro$-vertices and further that every vertex in $V_\zf\cup V_\fz$ is forced; and $\mathbf{J}^\pi_\zro$ is the indicator that every vertex in $V_\zz\cup V_\zo\cup V_\oz$ is ``forced,'' that is, has at least two $\one$-neighbors in each coordinate in which it takes spin $\zro$.

\noindent\emph{Matchings on $\free$-vertices.}\\
Let $\jfree(\pi)\equiv\E[\mathbf{J}^\pi_\free]$ denote the expected number of matchings on the $n_\free$ free vertices in the spin configuration $\uom$. We claim that for any valid $\pi$,
\beq\label{e:expected.matchings.f}
\begin{array}{l}
(2n)^{-1} \kfree(\pi) \le \jfree(\pi)\le\kfree(\pi)
	\quad\text{where }\\
\DS\kfree(\pi)\equiv
\f{ d^{n_\free+n_\ffne} }
	{ [d/(d-1)]^{n_\ffne} }
\f{
	(n_\ff-1)!!
	(n_\zf)_{n_\of} (n_\fz)_{n_\fo}
	(\numzf+n_\ffne-1)!!
	(\numfz+n_\ffne-1)!!
	}
	{ \fdf{ nd }{ (n_\free+n_\ffne)/2 }}.
\end{array}
\eeq
(For each vertex in $V_\free\setminus V_{\ffne}$ we distinguish one half-edge to participate in the matching. For each $\ffne$-vertex we distinguish an ordered pair of half-edges, the first to participate in $\match^1$ and the second in $\match^2$. There are $n_\free+n_\ffne$ distinguished half-edges in total.) If $n_\free<2$ then $\jfree(\pi)=\kfree(\pi)$; in general, $\kfree$ upper bounds $\jfree$ because it counts matchings without enforcing the constraint that two $\ffne$-vertices cannot be matched in both $\match^1$ and $\match^2$.

For the lower bound when $n_{\ffne}\ge2$, suppose without loss that $\numfz\ge\numzf$, and note this implies $a^\bullet \equiv \numfz+n_\ffne$ must be at least $4$ --- otherwise $n_\ffne=2$ while $\numzf=\numfz=0$ so there is no valid matching on the $\ffne$-vertices. Say that all the distinguished half-edges have been assigned except for the $a^\bullet$ half-edges incident to $\Vumfz\cup V_\ffne$ which were chosen to participate in $\match^1$. Now match these remaining half-edges one pair at a time, but avoid forming pairs already present in $\match^2$. The number of choices for the first $(a^\bullet/2)-1$ pairs is $\ge (a^\bullet-2)(a^\bullet-4)\cdots2 \ge (a^\bullet-3)!!
	\ge (a^\bullet-1)!! /n $.
	
\begin{figure}[h]
\centering
\begin{subfigure}[h]{\textwidth}\centering
	\includegraphics[width=.5\textwidth]{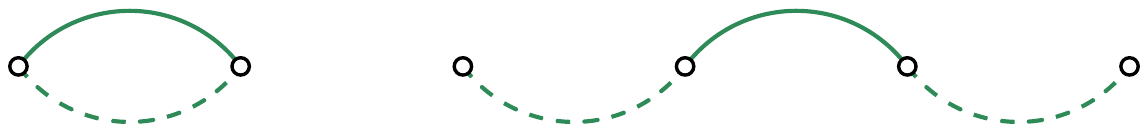}
	\caption{Matching of final $a^\bullet$ half-edges
	fails because final (leftmost) pair already present in $\match^2$}
	\label{f:bad}
\end{subfigure}\\ \medskip
\begin{subfigure}[h]{0.45\textwidth}\centering
	\includegraphics[width=\textwidth]{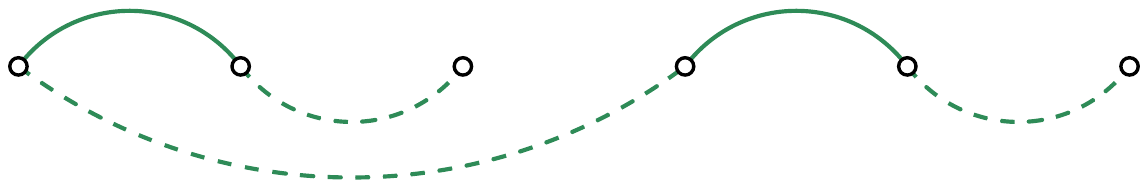}
	\caption{One valid switching}
	\label{f:good.first}
\end{subfigure}\quad
\begin{subfigure}[h]{0.45\textwidth}\centering
	\includegraphics[width=\textwidth]{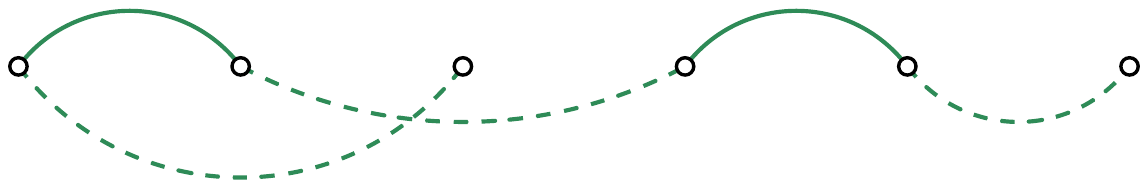}
	\caption{A distinct valid switching}
	\label{f:good.second}
\end{subfigure}
\caption{Two valid switchings of a failed matching
	(dashed lines $\match^1$, solid lines $\match^2$).}
\label{f:fmatch.switching}
\end{figure}

The procedure succeeds if and only if the final pair remaining is not already present in $\match^2$. To bound the probability that it fails, note that if given a failed matching in which the final pair is already present in $\match^2$, we can choose any of the first $(a^\bullet/2)-1$ pairs, and switch the half-edges in one of two ways to produce a valid matching (Fig.~\ref{f:fmatch.switching}). Thus each failed matching maps to $a^\bullet-2$ valid matchings.

In the reverse direction, if the final pair is given then the failed matching can be uniquely recovered from the valid matching, so each valid matching has at most $n_\ffne/2$ preimages.
Thus the ratio of failed to valid matchings is
(recalling $a^\bullet\ge4$) at most
$(n_\ffne/2)/(a^\bullet-2)
\le (2-4/a^\bullet)^{-1}\le1$.
This proves that the matching procedure succeeds with probability
at least $\slf12$, and the lower bound in \eqref{e:expected.matchings.f} follows. Combining with $\mathbf{c}(\pi)$ gives
\[
n^{O(1)}\,\mathbf{c}(\pi)\,\jfree(\pi)
= \f{ d^{(n_\free+n_\ffne)/2}\,
	(1-\slf1d)^{n_\ffne}\,
	2^{-n_\ff/2}\,
	n!\,
	(\numzf+n_\ffne-1)!!\,
	(\numfz+n_\ffne-1)!!}
{ 
(\prod_{\om\in\set{\zro,\one}^2}n_\om!)\,
(\slf{n_\ff}{2})!
	\,n_\of!
	\,n_\fo!
	\,\numzf!
	\,\numfz!
	\,\prod_{i=0}^{(n_{\free}+n_\ffne)/2-1}
		(n- \slf1d - \slf{2i}{d})
	}.\]

\medskip\noindent\emph{Edges from $\one$-vertices.} \\
Fix any matching $\match$ on the $\free$-vertices, and let $\ue[\om]$ denote the number of unmatched half-edges incident to $V_\om$ given the matching on the $\free$-vertices; the total is $\ue=\sum_\om\ue[\om] = nd-n_\free-n_\ffne$. Consider the assignment of the remaining $\ue[\one] = n_\one d -n_\of-n_\fo$ unmatched half-edges from vertices with spin $\one$ in either coordinate. Recalling $\forz\equiv\set{\zro,\free}$, define $\chi_{K,F}$ to be the indicator that there are $K$ edges between $V_\oz$ and $V_\zo$, $F_1$ edges between $V_\og$ and $V_\zf$, and $F_2$ edges between $V_\go$ and $V_\fz$; and let $\mathbf{j}_{\one,K,F}(\pi,\match)\equiv \E[\mathbf{J}^\pi_\one \chi_{K,F}\giv\match]$. Then, writing $F_\bullet \equiv F_1+F_2$,
\begin{align*}
\mathbf{j}_{\one,K,F}(\pi,\match)
&= \f{(\ue[\og])_{F_1} (\ue[\zf])_{F_1}}
	{(F_1)!}
\f{(\ue[\go])_{F_2} (\ue[\fz])_{F_2}}
	{(F_2)!}
\f{(\ue[\oz])_K(\ue[\zo])_K}{K!}
\f{(\ue[\zz])_{\ue[\one]-2K-F_\bullet}}
	{ \fdf{\ue}{ \ue[\one]-K } }\\
&= 
\exp\{ O(F_\bullet \tf{\log d}{d}) \}
\underbrace{
	\f{(\ue[\og])_{F_1} (\ue[\zf])_{F_1}
	(\ue[\go])_{F_2} (\ue[\fz])_{F_2}}
	{(F_1)!(F_2)!
	(\ue[\zz]-\ue[\one])^{F_\bullet}
	}
	}_{\mathbf{b}_{\one,F}(\pi,\match)}
\times
	\underbrace{
	\f{(\ue[\oz])_K(\ue[\zo])_K
	(\ue[\zz])_{\ue[\one]-2K}}
	{K! \fdf{\ue}{\ue[\one]-K} }
	}_{\wt{\mathbf{a}}_{\one,K}(\pi,\match)}.
\end{align*}
Let $\mathbf{a}_{\one,K}(\pi,\match)\equiv\wt{\mathbf{a}}_{\one,K}(\pi,\match)\,/\,\gone(\ue,\ue[\one],\ue[\zz])$ where $\gone$ is as in \eqref{e:unweighted.one}. We then estimate
\beq\label{e:edges.from.ones.K}
\mathbf{a}_{\one,K}(\pi,\match)
\bigg[
\f{(\ue[\oz]\ue[\zo])^K}
	{K!(\ue[\zz]-\ue[\one])^K}
\bigg]^{-1}
\begin{cases}
\le \exp\{O(K\bemax)\}
	&\text{for all $K$;}\\
\ge  \exp\{ O(K \tf{\log d}{d})\}
	&\text{for
	$K\lesssim K^\star
	\equiv \ue[\oz]\ue[\zo]/\ue$.}
\end{cases}.
\eeq
Similarly, writing $F^\star\equiv(F_1,F_2)^\star\equiv(\ue[\og]\ue[\zf]/\ue,\ue[\go]\ue[\fz]/\ue)$, we have
$$\mathbf{b}_{\one,F}(\pi,\match)
\bigg[
\f{
(\ue[\og]\ue[\zf])^{F_1}
(\ue[\go]\ue[\fz])^{F_2}
}
{(F_1)!(F_2)!(\ue[\zz]-\ue[\one])^{F_\bullet}}	
\bigg]^{-1}
\begin{cases}
\le 1 &\text{for all $F$}, \\
\ge\exp\{ O(F_\bullet\tf{\log d}{d}) \}
	&\text{for }
	F \lesssim F^\star.
\end{cases}$$

\medskip\noindent\emph{Forcing of $\zro$-vertices.} \\
Recall \eqref{e:unweighted.zof} that $\gzro(n,nd\target[\od])$ denotes the probability, with respect to a uniformly random assignment of $nd\target[\od]$ half-edges to $n$ degree-$d$ vertices, that each vertex is forced (that is, receives at least two of the incoming half-edges). Similarly, define the following:
\bnm[1.]
\item For $n^\bullet\le n$, suppose each variable $1\le i\le n^\bullet$ has degree $d$, while each variable $n^\bullet<i\le n$ has already received one forcing half-edge, and now has degree $d-1$. Let $\gzroum(n,n^\bullet,nd\target[\od])$ denote the probability, with respect to a uniformly random assignment of $nd\target[\od]$ half-edges to the $n$ vertices, that each vertex is forced (thus, each vertex $n^\bullet<i\le n$ is required only to receive at least one incoming half-edge).

\item Writing
$\target\equiv
(\target[\oo],\target[\oz],\target[\zo])$,
let 
$\gzz(n,nd\target)$
denote the probability,
with respect to a uniformly random assignment
of
$nd\target$
half-edges
to $n$ degree-$d$ vertices,
that each vertex
is forced in both coordinates (that is, receives at least two of the incoming half-edges both from $\set{\oo,\oz}$ and from $\set{\oo,\zo}$).
\enm
Then
$\mathbf{j}_{\zro,K,F}(\pi,\match)
\equiv
\E[\mathbf{J}^\pi_\zro \giv \match, \mathbf{J}^\pi_\one\chi_{K,F} >0 ]$
can be written as
\begin{align*}
\mathbf{j}_{\zro,K,F}(\pi,\match)
&=
\overbrace{\gzro(n_\oz,K)\,
	\gzro(n_\zo,K)}^{\text{(\sc a)}}
\times
\overbrace{
\gzroum(n_\zf,\numzf,F_1)\,
	\gzroum(n_\fz,\numfz,F_2)
}^{\text{(\sc b)}}\\
&\qquad\times
\underbrace{\gzz(n_\zz,
	(
	\ue[\oo],\ue[\og]-K-F_1,
	\ue[\go]-K-F_2
	))}_{\text{(\sc c)}}.
\end{align*}
The functions $\gzro$ and $\gzz$ were estimated in Propn.~\ref{p:forcing}. It is clear that $\gzroum$ satisfies the same estimates as $\gzro$, so we conclude
$$
\begin{array}{rll}
\text{\sc (a)}
	\hspace{-6pt}&\text{is }
	\exp\{ O( (n_\oz+n_\zo)
	 d^{-4\erho/5} )\}
&\text{for }
	K\ge
	\slf{9}{10}\cdot
	\slf{\ue[\oz] \ue[\zo]}{\ue};\\
\text{\sc (b)}
	\hspace{-6pt}&\text{is }
		\exp\{ O( ( n_\zf+n_\fz )
		d^{-4\erho/5} )\}
&\text{for }
	F_1\ge\slf{9}{10}
		\cdot \slf{\ue[\og]\ue[\zf]}{\ue}, \
	F_2\ge\slf{9}{10}
		\cdot \slf{\ue[\go]\ue[\fz]}{\ue};
	\\
\text{\sc (c)}
\hspace{-6pt}&\text{is }
	\exp\{O( nd^{-1.65} )\}
&\text{for }
	K \lesssim n_\oz \wedge n_\zo.
\end{array}
$$

\medskip\noindent\emph{Comparison of frozen model with independent set.} \\
From the preceding estimates, $\argmax_K \mathbf{a}_{\one,K}(\pi,\match)$ must be close to $K^\star$. If $K$ deviates from $K^\star$ by more than a $(1\pm\ep)$-factor for $\ep$ a small constant uniform in $d$, then the resulting decrease in $\mathbf{a}_{\one,K}(\pi,\match)$ will be much larger than any possible gain in $\mathbf{j}_{\zro,K,F}(\pi,\match)$:
$$
\begin{array}{rl}
\mathbf{a}_{\one,K}(\pi,\match)
	\,/\,\mathbf{a}_{\one,K^\star}(\pi,\match)
	\hspace{-6pt}&\le
	\exp\{ - \tf23 \ep^2 \cdot K^\star \}
	\le \exp\{ -\tf12\ep^2 \cdot n (\erho\log d)^2/d \},\\
\mathbf{j}_{\zro,K,F}(\pi,\match)
	\,/\,
	\mathbf{j}_{\zro,K^\star,F}(\pi,\match)
	\hspace{-6pt}&\le
	1 \,/\,
	\mathbf{j}_{\zro,K^\star,L}(\pi,\match)
	\le
	\exp\{ O( n(\erho\log d)/d ) \}.
\end{array}
$$
Therefore we may restrict consideration to $\slf{9}{10}\cdot K^\star \equiv \underline{K} \le K \le \overline{K} \equiv \slf{11}{10}\cdot K^\star$, since the contribution to $\E[\ZZt(\pi)]$ from $K$ outside this regime must be an exponentially small fraction of the total. By the same reasoning we may similarly restrict $F\equiv(F_1,F_2)$ to lie in the regime $\slf{9}{10}\cdot F^\star\equiv\underline{F} \le F \le \overline{F}\equiv \slf{11}{10}\cdot F^\star$, therefore
\beq\label{e:restrict.K.F}
\E[\ZZt(\pi)]
\asymp \mathbf{c}(\pi)\,\jfree(\pi)\,\E_{\match}
	\bigg[ \sum_{\underline{K}\le K\le \overline{K}}
	\sum_{\underline{F}\le F\le \overline{F}}
	\mathbf{j}_{\one,K,F}(\pi,\match)\,
	\mathbf{j}_{\zro,K,F}(\pi,\match)
	\bigg]
\eeq
where $\E_{\match}$ denotes the expectation over the random matching $\match$, conditioned on $\mathbf{J}^\pi_\free$ being positive. For any $\pi\in\cM(n\al,n\rho)$, consider the measure $\zpi\in M(n\al,n\rho)$ defined by
\[\zpi_\oo=\pi_\oo,\quad
	\zpi_\oz=\pi_\oz+\pi_\of+\tf12\pi_\fd,\quad
	\zpi_\zo=\pi_\zo+\pi_\fo+\tf12\pi_\dotf,\]
with the remaining probability going to $\zpi_\zz$. Then
\[\f{\E[\ZZt(\pi)]}{\E[\Zt(\zpi)]}
	=\f{\mathbf{c}(\pi)\,\jfree(\pi)}{\mathbf{c}(\zpi)}
	\f{\E_{\match}[ \sum_{K,F}
		\mathbf{j}_{\one,K,F}(\pi,\match)
		\mathbf{j}_{\zro,K,F}(\pi,\match) ] }
	{\sum_K \mathbf{j}_{\one,K,0}(\zpi)},\]
and it follows from the above estimates that this is $d^{O(n\bemax)}\exp\{ O(nd^{-1.65} + n_\one d^{-3\erho/4})\}$, implying the result.
\epf
\eppn

\bcor\label{c:rule.out.intermediate}
Suppose $\albd\le\al\le\aubd$ and $\erho>0$ a small positive constant uniform in $d$. For $d\ge d_0(\erho)$, $\bPhit$ can only attain its global maximum on $\simplextal$ either in the near-independence regime $\simplextalind\subseteq\simplextal$ of measures with overlap $\rho\le d^{-1}$, or in the near-identical regime $\simplextalid\subseteq\simplextal$ of measures with overlap $\al-\erho\logdbyd\le \rho\le \al$.

\bpf
We shall prove that for $d\ge d_0(\erho)$ and $n\ge n_0(d)$, the ratio of $\E[\ZZt_{n\al}(n\rho)]$ to $\max_r\E[\ZZt_{n\al}(nr)]$ is exponentially small in $n$ for all $d^{-1.45}\le\rho\le\al-\erho\logdbyd$. We begin with an analogous calculation in the independent set model: there is no matching on the $\free$-vertices and no forcing constraint on the $\zro$-vertices, so
$$\f{\E[\Zt_{n\al}(n\rho)]}{\E[Z_{n\al}]}
=
\sum_{\pi\in M(n\al,n\rho)}\f{
	\mathbf{c}(\pi)
	\E[\mathbf{J}^\pi_\one]}
{\binom{n}{n\al} \gone(nd,nd\al,0)}
=\sum_{\pi\in M(n\al,n\rho)}\f{
	\mathbf{c}(\pi)
	\gone(nd,n_\one d,0)
	\sum_K \mathbf{a}_{\one,K}(\pi)
	}
{\binom{n}{n\al} \gone(nd,nd\al,0)},$$
where $\mathbf{a}_{\one,K}(\pi)\equiv\mathbf{a}_{\one,K}(\pi,\emptyset)$ for $\mathbf{a}_{\one,K}(\pi,\match)$ as in \eqref{e:edges.from.ones.K}. For $n_\one\lesssim n \logdbyd$ we have
$$\gone(nd,n_\one d,0)=
\exp\{ -\tf12 n_\one^2 d/n
+ O(n d^{-2}(\log d)^3 ) \}.$$
The estimates of Propn.~\ref{p:ap.compare.frozen.is} imply
$$
\sum_K \mathbf{a}_{\one,K}(\pi)
=\exp\{ O(K\tf{\log d}{d}) \}
\sum_K 
\f{( n_\oz n_\zo d)^K}
	{K!(n_\zz-n_\one)^K}
= \exp\{ n_\oz n_\zo d/n
	+O(n d^{-2}(\log d)^3 )\}.
$$
Combining gives
\begin{align*}
\E[\Zt_{n\al}(n\rho)]
&=\E[Z_{n\al}]\,
\exp\{ ng_\al(\rho)
	+ O( n d^{-2}(\log d)^3 )
	 \}
\quad\text{where}\\
g_\al(\rho)
&\equiv \al H( \tf{\rho}{\al})
	+(1-\al) H(\tf{\al-\rho}{1-\al})
	-\tf{d}{2} n\al^2
	+ \tf{d}{2} n\rho^2.
\end{align*}
The function $g_\al$ has first derivative $(g_\al)'(\rho) = 2\log(\al-\rho)-\log\rho-\log(1-2\al+\rho)+d\rho$; differentiating again gives $(g_\al)''(\rho)=d-2(\al-\rho)^{-1}-\rho^{-1}-(1-2\al+\rho)^{-1}$, 
so we see that $g_\al$ is strictly convex on the interval
$2/d\le\rho\le\al-3/d$.
From the expression for $(g_\al)'$ we see that the (unique) minimizer $\rho_\circ$ of $g_\al$ on this interval must lie near $\logdbyd$, and so applying Propn.~\ref{p:ap.compare.frozen.is} gives
\begin{align}
\label{e:sup.moments.ratio.left}
\f{\sup_{2d^{-1.48}\le\rho\le\rho_\circ}
\E[\ZZt_{n\al}(n\rho)]}
	{\E[\ZZt_{n\al}(nd^{-1.48} )]}
&=
\exp\{ O(nd^{-1.49})\}
\f{\sup_{2d^{-1.48}\le\rho\le\rho_\circ}
\E[\Zt_{n\al}(n\rho)]}
	{\E[\Zt_{n\al}(nd^{-1.48} )]
	};\\
\label{e:sup.moments.ratio.right}
\f{\sup_{\rho_\circ\le\rho\le\al-\erho \logdbyd}
	\E[\ZZt_{n\al}(n\rho)]}
{
\E[\ZZt_{n\al}(n\al
	- n\erho(\log d)/(2d)
	)]}
&=\exp\{ O(n/d) \}
\f{\sup_{\rho_\circ\le\rho\le\al-\erho \logdbyd}
	\E[\Zt_{n\al}(n\rho)]}
{
\E[\Zt_{n\al}(n\al- n\erho(\log d)/(2d)
)]
}.
\end{align}
We estimate
$(g_\al)'(\rho) \le -(\log d)/10$ for $d^{-1.85}\rho \le 100/d$, 
and similarly $(g_\al)'(\rho) \ge (\log d)/10$ for
$1.6(\log d)/d\le \rho \le \al-d^{-1.25}$.
 Thus
\begin{align*}
\eqref{e:sup.moments.ratio.left}
\text{ is }
&\le
 \exp\{ n\,
 [g_\al(2d^{-1.48}) - g_\al( d^{-1.48} )]
 +O(nd^{-1.49}) \}
\le \exp\{ - n d^{-1.48}\},\\
\eqref{e:sup.moments.ratio.right}
\text{ is }
&\le
\exp\{ n
[g_\al(n\al-n\erho\logdbyd)
-g_\al(n\al-n\erho(\log d)/(2d) )
]
+O(n/d)\}\\
&\le
\exp\{ -n\erho\logdbyd \}.
\end{align*}
These estimates cover the entire interval $d^{-1.45}\le \rho \le \al-\erho\logdbyd$, implying the result.
\epf
\ecor

\subsection{Boundary estimates}

\bppn\label{p:ap.intermediate}
For $\albd\le\al\le\aubd$ the following hold:
\bnm[(a)]
\item\label{p:ap.intermediate.a}
For $d\ge d_0$, $n\ge n_0(d)$,
$\argmax_{0\le\be\le\bemax}
\E[\ZZ_{n(\al-\be/2),n\be}]$ is positive.

\item \label{p:ap.intermediate.full.supp}
Suppose $0 \le \rho\le\al-\erho\logdbyd$ for $\erho$ a small positive constant uniform in $d$. Then for $d\ge d_0(\erho)$, $n\ge n_0(d)$, any measure $\pi\in\cM(n\al,n\rho)$ maximizing $\E[\ZZt(\pi)]$ has full support on $\pairs$; further $\pumzf\equiv\pi_\zf-\pi_\of$ and $\pumfz\equiv\pi_\fz-\pi_\fo$ must be positive.
\enm

\bpf
Suppose that $\pi_\om\equiv n_\om/n$ is the maximizer of $\E[\ZZt(\pi)]$ on $\cM(n\al,n\rho)$. We prove \eqref{p:ap.intermediate.full.supp} by a series of comparison estimates. Write $\de\equiv n/e^{d^2}$.
\bnm[1.]
\item Suppose $\pi_\zf=0$, and consider the measure
\beq\label{e:pi.add.zf}
\zpi
=\pi + (\de/n)\,(2 \cdot\Ind{\zf}-\Ind{\zz,\zo})
\in\cM(n\al,n\rho).\eeq
Let $\match$ denote any matching on the $\free$-vertices for the $\pi$-configuration; clearly $\ue[\zf]=F_1=0$. Likewise let $\zmatch$ denote any matching on the $\free$-vertices for the $\zpi$-configuration; note that in this case $0\le F_1\le\overline{\ue}_\zf=2\de d-O(\de)$. We estimate
\begin{align*}
&\f{\mathbf{j}_{\one,K,(F_1,F_2)}(\zpi,\zmatch)}
	{\mathbf{j}_{\one,K,(0,F_2)}(\pi,\match)}
= 
\f{(\ue[\og])_{F_1}(\overline{\ue}_\zf)_{F_1}}{  F_1 !\, \ue^{F_1} \, d^{O(\de)} }
	\f{ (\ue[\go]-\de d)_{F_2} }{ (\ue[\go])_{F_2} }
	\f{ (\ue[\zo]-\de d)_K }{ (\ue[\zo])_K }
\stackrel{\star}{=}
d^{O(\de)},\\
&\f{\mathbf{j}_{\zro,K,(F_1,F_2)}(\zpi,\zmatch)}
	{\mathbf{j}_{\zro,K,(0,F_2)}(\pi,\match)}
= e^{O(\de)}\,\gzroum(\zn_\zf,\zn_\zf,F_1)
= e^{O(\de)}\,\gzro(\de,F_1)
\stackrel{\star}{=} e^{O(\de)}
\end{align*}
where $\star$ holds for $K,F$ appearing in the sum \eqref{e:restrict.K.F} with respect to the $\zpi$-configuration.
Thus
$$
\f{\E[\ZZt(\zpi)]}{\E[\ZZt(\pi)]}
=
\f{\mathbf{c}(\zpi)\,\jfree(\zpi)}
{\mathbf{c}(\pi)\,\jfree(\pi)}
\underbrace{
\f{
\E_{\zmatch}[\sum_{K,F}
	\mathbf{j}_{\one,K,F}(\zpi,\zmatch)\,
	\mathbf{j}_{\zro,K,F}(\zpi,\zmatch)
	]}
	{ \E_{\match}[\sum_{K,F_2}
	\mathbf{j}_{\one,K,(0,F_2)}(\pi,\match)\,
	\mathbf{j}_{\zro,K,(0,F_2)}(\pi,\match)
	]}
	}_{ \ge d^{O(\de)} }
	\ge \f{(n/\de)^\de}{d^{O(\de)}}
	\gg1,
$$
contradicting the assumption that $\pi$ was a maximizer $\E[\ZZt(\pi)]$ on $\cM(n\al,n\rho)$.

\item 
Since the preceding estimate shows that $\pi_\zf$ and $\pi_\fz$ must be positive, we must have $\de'\equiv (n_\zf \wedge n_\fz)/e^{d^2}\gg0$ scaling linearly with $n$. Suppose $\pumzf=0$: if we let $\zpi$ be defined by \eqref{e:pi.add.zf} with $\de'$ in place of $\de$, then a simple term-by-term comparison gives
$$
\f{\E[\ZZt(\zpi)]}{\E[\ZZt(\pi)]}
= e^{O(\de')}
\f{\mathbf{c}(\zpi)\,\jfree(\zpi)}
{\mathbf{c}(\pi)\,\jfree(\pi)}
\ge \f{(n/\de')^{\de'}}{d^{O(\de')}}
	\gg1,
$$
again a contradiction. Similarly, if $\pi_\ff=0$, then
$\E[\ZZt(\pi)]\ll\E[\ZZt(\zpi)]$
$$\zpi
=\pi + (\de'/n)\,(2 \cdot \Ind{\ff}- \Ind{\oz,\zo})\in\cM(n\al,n\rho).$$
This proves that $\pumzf$, $\pumfz$, and $\pi_\ff$ must be positive.

\item Let $\de''\equiv (\numzf\wedge\numfz\wedge n_\ff)/e^{d^2}$; by the preceding estimates we have $\de''\gg0$ scaling linearly with $n$. By similar calculations as before,
if $\pi_\ffne=0$ then
$\E[\ZZt(\pi)]
\ll \E[\ZZt(\zpi)]$ for
$$\zpi=\pi+(\de''/n)\,(
	2\cdot\Ind{\ffne}
	-\Ind{\zf,\fz} )
	\in\cM(n\al,n\rho).$$
If $\pi_\of=0$, then
$\E[\ZZt(\pi)]
\ll \E[\ZZt(\zpi)]$ for
$$\zpi=\pi+(\de''/n)\,(\Ind{\of,\zf}
	-2\cdot\Ind{\ff})
	\in\cM(n\al,n\rho).$$
\enm
This concludes the proof of \eqref{p:ap.intermediate.full.supp}. The first-moment estimate \eqref{p:ap.intermediate.a} follows simply by restricting consideration to $\pi$ supported on $\set{\zz,\oo,\ff}$.
\epf
\eppn

\bppn\label{p:boundary}
For $\albd\le\al\le\aubd$ the following hold:
\bnm[(a)]
\item\label{p:boundary.a}
For $d\ge d_0$, any global maximizer of $\bPhi$ on $\simplexal$ must be strictly positive on $\supp\phi$.
\item\label{p:boundary.b}
Fix any small constant $\erho>0$ uniform in $d$. For $d\ge d_0(\erho)$, any global maximizer of $\bPhit$ on $\simplextal$ which lies outside of $\simplextalid$ must be strictly positive on $\supp\psit$.
\enm

\newcommand{\bhv}{{}^\mathrm{v}\hspace{-1pt}\bh}
\newcommand{\vhv}{{}^\mathrm{v}\hspace{-1pt}\vh}
\newcommand{\dbhv}{{}^\mathrm{v}\hspace{-1pt}\dbh}
\newcommand{\psitv}{{}^\mathrm{v}\hspace{-1pt}\psit}
\newcommand{\dpsitv}{{}^\mathrm{v}\hspace{-1pt}\dpsit}
\newcommand{\bPhitv}{{}^\mathrm{v}_2\hspace{-1pt}\bm{\Phi}}
\newcommand{\dbdev}{{}^\mathrm{v}\hspace{-1pt}\dbde}
\newcommand{\bdev}{{}^\mathrm{v}\hspace{-1pt}\bde}
\newcommand{\vdev}{{}^\mathrm{v}\hspace{-1pt}\vde}
\newcommand{\nuauxv}{{}^\mathrm{v}\hspace{-1pt}\nuaux}

\bpf
We shall prove \eqref{p:boundary.b}; the proof of \eqref{p:boundary.a} is similar but simpler.

\medskip\noindent\emph{Boundary derivative of rate function.} \\
As we have noted before, the functional form of $\bPhit$ implies that the optimal $\bh$ in $\simplextal$ must be symmetric, with
$$\TS\bPhit(\bh)
=\sum_{\dta}\dbh(\dta)\log\dpsit(\dta)
+ H(\dbh)
- \tf{d}{2}H(\vh)$$
For $\bde\equiv(\dbde,\hbde)$ such that $\bh+t\bde$ is also symmetric and lies in $\simplextal$ for $t\ge0$ small, consider
$$
\Tlog(\bPhit)(\bh;\bde)\equiv
\lim_{t\downarrow0}\f{\bPhit(\bh+t\bde)-\bPhit(\bh)}{t\log(1/t)}
=
\dbde[(\supp\dbh)^c]
-(d/2)\, \vde[(\supp\vh)^c].
$$
To show that $\bh$ is not a maximizer it suffices to exhibit $\Tlog(\bPhit)(\bh;\bde)>0$ for some $\bde$. In particular, it follows by convexity that for any $\bh\in\simplextal$,
$\bh+t(\bhtstaral-\bh)\in\simplextalint$ for $t>0$ small
and $\bhtstaral$ as in the statement of Thm.~\ref{t:second.moment}. Therefore, if $\bh$ is a maximizer such that the edge marginal has full support $\supp\vh=\msg^2$, then necessarily $\supp\bh=\supp\psit$, since otherwise $\Tlog(\bPhit)(\bh;\bhtstaral-\bh)>0$.

\medskip\noindent\emph{Positive maximizer for vertex-auxiliary model.} \\ Clearly the preceding calculation applies equally well to the rate function $\bPhitv$ of the pair version of the \emph{vertex-auxiliary} model introduced in Rmk.~\ref{r:vertex.aux}, with factors $\psitv\equiv\phiv\otimes\phiv$. We now show that the maximizer $\bhv$ of $\bPhitv$ in the regime $\rho\le\al-\erho\logdbyd$ satisfies $\supp\vhv=\msgv^2$, hence (from the above) $\supp\bhv=\supp\psitv$.

By Propn.~\ref{p:ap.intermediate}\ref{p:ap.intermediate.full.supp},
$\bhv$ must satisfy $\vhv(\spint{r}{s})>0$ for $r,s$ any of the sets in the partition $\set{\oz}$, $\set{\zo}$, $\set{\ffp}$, $\zg=\set{\zz,\zf}$, $\fg=\set{\fz,\ff}$ (e.g.\ $\vhv(\spint{\ffp}{\fg})>0$ due to the presence of $\ffne$-variables. If $\tau\equiv(\si^1,\si^2)\notin\supp\vhv$, then $\acute\tau\equiv(\acute\si^1,\acute\si^2)\in\supp\vhv$ for some choice of $\acute\si^i$ belonging to the same subset as $\si^i$ in this partition. Therefore $(\acute\tau,\vec\tau)\in\supp\dbhv$ for some $\vec\tau\in\msgv^{2(d-1)}$, so $(\tau,\vec\tau)\in\supp\dpsitv\setminus\supp\dbhv$. 
By the symmetry $\vhv(\si)=\vhv(\refl\si)$, there
must exist $\vec\tau'\in\msgv^{2(d-1)}$ such that $(\refl\acute\tau,\vec\tau')\in\supp\dbhv$. If $\tau$ belongs to
$\forz\forz\times(\forz\forz\cup\set{\ffp,\oz,\zo})$,
then $(\refl\tau,\vec\tau')$ must be in $\supp\dpsitv\setminus\supp\dbhv$. Taking
$$\dbdev
=
\I_{(\tau,\vec\tau)}-\I_{(\acute\tau,\vec\tau)}
+\I_{(\refl\tau,\vec\tau')}-\I_{(\refl\acute\tau,\vec\tau')},\quad
\vdev
=\tf1d (\I_\tau - \I_{\acute\tau}
	+ \I_{\refl\tau} - \I_{\refl\acute\tau})
$$
gives $\Tlog(\bPhitv)(\bhv;\bdev) = 2-1>0$,
proving our claim that the vertex-auxiliary maximizer must have
$\supp\vhv=\msgv^2$, hence $\supp\bhv=\supp\psitv$.

\begin{figure}[ht]
\centering
\includegraphics[height=1in,trim=.7in .8in .7in .7in,clip]{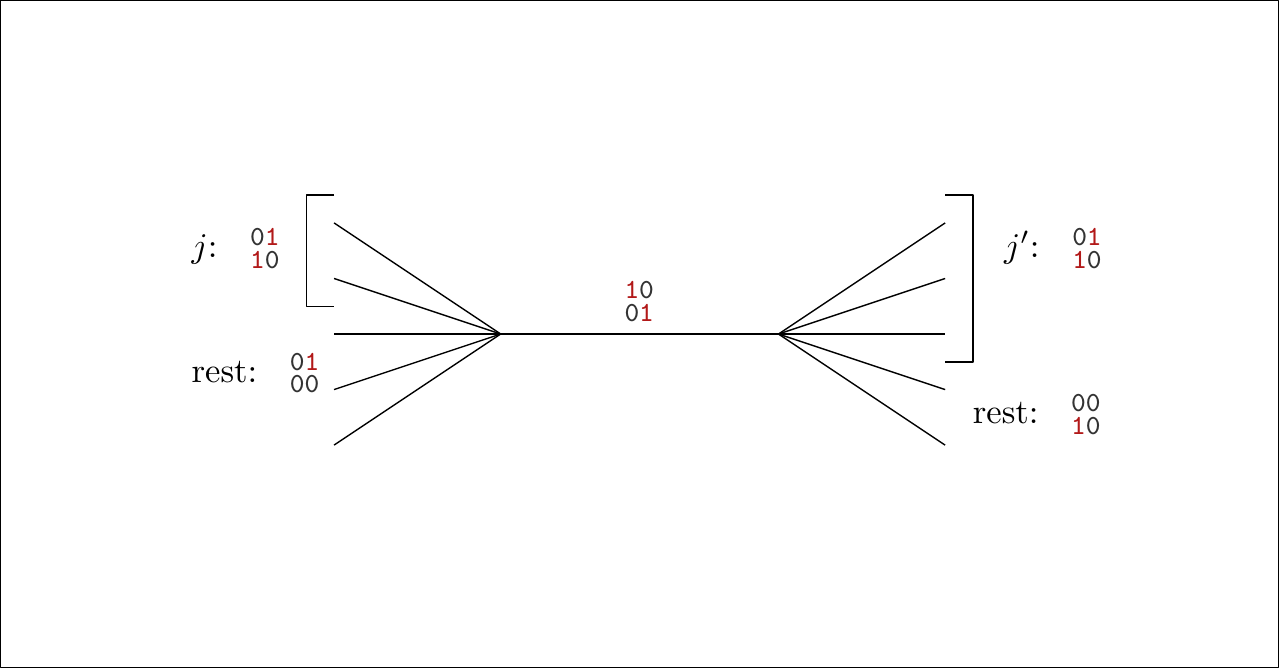}
\caption{Any combination of $j,j'$ must appear in vertex-auxiliary model}
\label{f:two.factors}
\end{figure}

\medskip\noindent\emph{Positive maximizer for auxiliary model.} \\
It directly follows from $\supp\vhv=\msgv^2$ that the maximizer $\bh$ for the original auxiliary model on $\simplextal\setminus\simplextalid$ satisfies $\vh(\spint{r}{s})>0$ for $r,s$ any of the sets in the partition $\set{\zz}$, $\set{\zf}$, $\set{\fz}$, $\set{\ff}$, $\set{\oo}$, $\og$, $\go$. Further, since $\supp\dbhv=\supp\dpsitv$ we may replace $\go$ in the above list by $\set{\zo}$, $\set{\fo}$. By symmetry, the only remaining possibility for $\tau\notin\supp\vh$ is that either $\tau$ or $\refl\tau$ belongs to $\og\times\go$. Consider the Gibbs measure $\nuauxv$ associated to $\bhv$ on the subtree shown in Fig.~\ref{f:two.factors}: clearly, all combinations of $j,j'$ appear with positive $\nuauxv$-probability. Thus $\og\times\go$, and symmetrically $\go\times\og$, must belong in the support of $\vh$, concluding our proof that $\supp\bh=\supp\psit$.
\epf
\eppn

Recall \eqref{e:truncated.frozen.model} and 
Defn.~\ref{d:simplex} that we have truncated the frozen model and the spaces $\simplexal,\simplextal$ by restricting the density of $\free$-variables, so in order to establish $\bPhi$ and $\bPhit$ have \emph{interior} global maximizers in $\simplexal$ and $\simplextal$ respectively, we must verify in addition to Propn.~\ref{p:boundary} that the maximizer does not occur near the boundary with density $\bemax$ of $\free$-variables; this will be done in \S\ref{ss:indep}.

\subsection{Near-independence regime}\label{ss:indep}

In this subsection we complete our analysis of the near-independence regime $\simplextalind$ to prove
\bppn\label{p:indep.maximizer}
The unique global maximizer of the restriction of $\bPhit$ to $\simplextalind$ is $\bhtstaral$.
\eppn

Recall that in the pair frozen model, $\cM(n\al,n\rho)$ denotes the space of empirical measures $\pi$ on $\pairs$ which contribute to $\E[\ZZ_{n\al}^2]$ and have $\pi_\oo=\rho$.

\blem\label{l:ap.near.indep}
For $\albd\le \al\le\aubd$ let $y$ be defined by $\al\equiv y\logdbyd$. If $\ep$ is any small positive constant uniform in $d$, then the following hold for $d\ge d_0(\ep)$, $n\ge n_0(d)$:
\bnm[(a)]
\item \label{l:ap.near.indep.a}
The maximizer of $\E[\ZZ_{\al-\be/2,\be}]$ over $\be\le\bemax$ must satisfy $d^{y-1}\be/\al \le d^\ep$.
\item \label{l:ap.near.indep.b} For $\rho \le d^{-1}$, the maximizer of $\E[\ZZt(\pi)]$ over $\pi\in\cM(n\al,n\rho)$ must satisfy
$$
\f{\pi_\of}{\pi_\zf\pi_\oz},
\f{\pi_\fo}{\pi_\fz\pi_\zo},
\f{d \pi_\ff}{ \pi_\zf \pi_\fz},
\f{\pi_\ffne}{ \pi_\zf \pi_\fz},
\f{d^{y-1}\pi_\zf}{ \pi_\oz}
\le d^\ep.
$$
\enm
\vspace{4pt}

\bpf
Suppose the measure
$\pi_\om\equiv n_\om/n$ maximizes $\E[\ZZt(\pi)]$ over $\cM(n\al,n\rho)$. As in the proof of Propn.~\ref{p:ap.intermediate}, we shall estimate $\pi$ by comparing it with various nearby measures $\zpi_\om\equiv \zn_\om/n$. Note that if we fix any small positive constant $\ep$ which is uniform in $d$, then \eqref{e:restrict.K.F} continues to hold if the sum is restricted to
\beq\label{e:restrict.K.F.eps}
1-\ep\le
\set{K/K^\star,F/F^\star} \le 1+\ep.
\eeq
Let $\de\equiv e^{-d^2}\min\set{n_\om:\om\in\pairs}$; then Propn.~\ref{p:ap.intermediate} implies that $\de$ must scale linearly with $n$.

\bnm[1.]
\item \emph{Comparison of $\of$ to $\zf$.} \\
Let $\zpi=\pi+ (\de/n)\,(\Ind{ \oz,\zf } - \Ind{ \of,\zz })$. Let $\match$ (resp.\ $\zmatch$) be any matching of the $\free$-vertices for the $\pi$-configuration (resp.\ $\zpi$-configuration). For $\rho\lesssim d^{-1}$ and $K,F$ satisfying \eqref{e:restrict.K.F.eps},
\begin{align*}
\f{\mathbf{j}_{\one,K,F}(\zpi,\zmatch)}
{\mathbf{j}_{\one,K,F}(\pi,\match)}
&= e^{O(\de)}
\f{ (\ue[\zf]+\de d)_{F_1} }{ (\ue[\zf])_{F_1} }
\f{(\ue[\oz]+\de d)_K}{(\ue[\oz])_K}
\f{(\ue[\zz]-\de d)_{\ue[\one]-2K-F_\bullet}}
	{(\ue[\zz])_{\ue[\one]-2K-F_\bullet}}
\ge \f{e^{O(\de)}}{d^{O(\de\ep)}},
\end{align*}
where we used that
$(A+c)_b/(A)_b
= \exp\{ (bc/A) [1+O(x)] \}$
for $|b|,|c| \le xA$. Thus
$$
\f{\E[\ZZt(\zpi)]}{\E[\ZZt(\pi)]}
\ge 
\f{\mathbf{c}(\zpi)\,\jfree(\zpi)}
	{\mathbf{c}(\pi)\,\jfree(\pi)}
\f{e^{O(\de)}}{d^{O(\de\ep)}}
\ge
\Big(
\f{\pi_\of}
{ \pi_\oz\pi_\zf  }
\f{e^{O(1)}}{ d^{O(\ep)} }
\Big)^\de,
$$
so if $\pi$ is the maximizer then it must satisfy $\pi_\of \le d^{O(\ep)}\, \pi_\zf \pi_\oz$.

\item \emph{Comparison of $\ff,\ffne$ to $\zf,\fz$.} \\
Let $\zpi=\pi + (\de/n)\,(\Ind{\zf,\fz}-\Ind{\zz,\ff})$. By a similar calculation as above we find
$$\f{\E[\ZZt(\zpi)]}{\E[\ZZt(\pi)]}
\ge 
\f{\mathbf{c}(\zpi)\,\jfree(\zpi)}
	{\mathbf{c}(\pi)\,\jfree(\pi)}
\f{e^{O(\de)}}{d^{O(\de\ep)}}
\ge
\Big(
\f{d \pi_\ff}{ \pi_\zf\pi_\fz}
	\f{e^{O(1)}}{d^{O(\ep)}}
\Big)^{\de/2},
$$
proving that the maximizer $\pi$ must satisfy
$d\pi_\ff \le d^{O(\ep)}\, \pi_\zf\pi_\fz $. The calculation with $\ffne$ in place of $\ff$ is similar and shows that the maximizer $\pi$ must satisfy
$\pi_\ffne \le d^{O(\ep)}\,\pi_\zf\pi_\fz$.

\item \emph{Comparison of $\zf$ to $\zo$.}\\
Let $\zpi\equiv \pi+(\de/n)\,( \Ind{\zz,\zo}-2 \cdot \Ind{\zf} )$. Then
\begin{align*}
&\f{\mathbf{j}_{\one,K,F}(\zpi,\zmatch)}
{\mathbf{j}_{\one,K,F}(\pi,\match)}
=\f{(\ue[\zf]-2\de d)_{F_1}}{(\ue[\zf])_{F_1}}
	\f{ (\ue[\go]+\de d)_{F_2} }{ (\ue[\go])_{F_2} }
	\f{ (\ue[\zo]+\de d)_K }{ (\ue[\zo])_K }
		\f{ (\ue[\zz])^{\de d} \, e^{O(\de)}}{ 
			( \ue[\zz]-\ue[\one] )^{\de d} } 
\ge \f{d^{\de y} e^{O(\de)} } { d^{O(\de\ep)} }.
\end{align*}
From the preceding estimates,
$\pi_\of + \pi_\ffne \lesssim
\pi_\zf \logdbyd$, and so
$$\f{\E[\ZZt(\zpi)]}{\E[\ZZt(\pi)]}
\ge 
\f{\mathbf{c}(\zpi)\,\jfree(\zpi)}
	{\mathbf{c}(\pi)\,\jfree(\pi)}
\f{d^{\de y} e^{O(\de)} } { d^{O(\de\ep)} }
\ge
\Big(
\f{n\,\numzf }
	{n_\zz \, n_\zo}
\f{\numzf}{\numzf+n_\ffne}
\f{d^y \, e^{O(1)}}{ d^{1+O(\ep)} }
\Big)^\de
\ge \Big(
	\f{\pi_\zf}{\pi_\zo}
	\f{d^y e^{O(1)}}
	{ d^{1+O(\ep)} }
	\Big)^\de,
$$
implying that the maximizer $\pi$ must satisfy
$d^{y-1}\pi_\zf \le d^{O(\ep)} \,\pi_\zo$.
\enm
Adjusting $\ep$ as needed, this concludes the proof of \eqref{l:ap.near.indep.b}. To prove \eqref{l:ap.near.indep.a}, consider $\pi$ supported on $\set{\zz,\oo,\ff}$ with $\pi_\one=\al=y\logdbyd$ and $\pi_\free=\be$, and let 
$\zpi=\pi
+(\de/n)\,(\Ind{\zz,\oo}
	-2\cdot\Ind{\ff})$.
Clearly we must always have $K=F_1=F_2=0$,
and we estimate
$$\f{\E[\ZZt(\zpi)]}{\E[\ZZt(\pi)]}
\ge 
\f{\mathbf{c}(\zpi)\,\jfree(\zpi)}
	{\mathbf{c}(\pi)\,\jfree(\pi)}
e^{O(\de)} d^{\de y}
\ge \Big(
	\f{\pi_\free}{\pi_\one}
	\f{ d^{y-1} }{ e^{O(1)} }
	\Big)^\de,
$$
so that the maximizer must satisfy
$d^{y-1} \pi_\free \le
	d^{O(\ep)} \pi_\one$,
implying the result.
\epf
\elem

\bcor\label{c:boundary}
Let $\erho>0$ be a small constant uniform in $d$. For $d\ge d_0(\erho)$,
the following hold for all $\albd\le\al\le\aubd$:
\bnm[(a)]
\item \label{c:boundary.a} Any global maximizer of $\bPhi$ on $\simplexal$ lies in the interior $\simplexalint$.
\item \label{c:boundary.b} Any global maximizer of $\bPhit$ on $\simplextalind$ must be an interior stationary point.
\enm

\bpf
Propn.~\ref{p:boundary}\ref{p:boundary.a} and Lem.~\ref{l:ap.near.indep}\ref{l:ap.near.indep.a} combine to give \eqref{c:boundary.a}, while \eqref{c:boundary.b} follows by combining Cor.~\ref{c:rule.out.intermediate}, Propn.~\ref{p:boundary}\ref{p:boundary.b}, and Lem.~\ref{l:ap.near.indep}\ref{l:ap.near.indep.b}.
\epf
\ecor

Cor.~\ref{c:boundary} was required in the proof of Thm.~\ref{t:first.moment.exponent}; it also implies (with Lem.~\ref{l:interior.bp}) that any maximizer $\bh$ on $\bPhit$ on $\simplextalind$ corresponds to a solution $h$ of the pair Bethe recursions for some $\lm\equiv (\lm_1,\lm_2)$ (\eqref{e:bij} and \eqref{e:bp} with $\smash{\psit^\lm\equiv\phi^{\lm_1}\otimes\phi^{\lm_2}}$ in place of $\phi^\lm$). It remains to identify this Bethe solution with the one corresponding to $\bhtstaral=\bhstaral\otimes\bhstaral$.

We first show an analogous result in the simpler pair frozen model recursions, as follows. The unique fixed point $q=q^\lm$ of \eqref{e:frozen.recursions} defines a Gibbs measure $\nu$ on frozen model configurations $\uz$ on the infinite $d$-regular tree $\tree(t)$. We shall define now a Gibbs measure ${}_2\nu$ on \emph{pair} frozen configurations $(\uz^1,\uz^2)$: the message-passing rule is given simply by applying $\dmp_{d-1}$ separately in each coordinate. We weight according to $\lm_1$ raised to the number of upward messages in $\tree(t-1)$ which are $\one$ in the first coordinate, times $\lm_2$ raised to the number of upward messages which are $\one$ in the second coordinate. However, we will allow for more general boundary laws in which the incoming messages are not necessarily independent across the two copies: we will take the boundary messages $\tau_{v\to w}\equiv (\si^1_{v\to w},\si^2_{v\to w})$ to be distributed according to a measure $\tq$ on $\set{\zro,\one,\free}^2$ (still i.i.d.\ across the different boundary edges $v\to w$). Recall $\forz\equiv\set{\zro,\free}$, e.g.\ $\qgg\equiv \tq(\set{\zro,\free}^2) =1-2\qo+\qoo$, and $\qog = q_\one-\qoo = \qgo$. The reader can check that for $q$ to define a consistent family $({}_2\nu_t)_t$ of finite-dimensional distributions, its marginals $\tq_{\eta\cdot}\equiv\sum_{\eta'} \tq_{\eta\eta'}$ and $\tq_{\cdot\eta}\equiv\sum_{\eta'}\tq_{\eta'\eta}$ must satisfy the single-copy recursions \eqref{e:frozen.recursions} (with respect to $\lm_1$ and $\lm_2$ respectively), and further we must have
\beq\label{e:pair.frozen.recursions}
\begin{array}{l}
z\qoo = \lm_1\lm_2(\qgg)^{d-1},\quad
z\qof = \lm_1(d-1) \qgo (\qgg)^{d-2},\quad
z\qfo = \lm_2(d-1) \qog (\qgg)^{d-2}\\
z\qff = (d-1) \qoo (\qgg)^{d-2} + (d-1)(d-2) \qog\qgo (\qgg)^{d-3},
\end{array}
\eeq
with the rest being determined by margin constraints. Clearly $\smash{q^{\lm_1} \otimes q^{\lm_2}}$ solves \eqref{e:pair.frozen.recursions}, and the following lemma identifies a regime in which it is the unique solution:

\blem\label{l:pair.frozen.recursions}
For $d^{1/10} \le \set{\lm_1,\lm_2} \le d^2$, the recursions \eqref{e:pair.frozen.recursions} have the unique solution $\smash{q^{\lm_1}\otimes q^{\lm_2}}$ in the regime
$\qoo \le \slf{9}{10}\cdot \logdbyd$.

\bpf
Write ${}_i q$ for the solutions of \eqref{e:frozen.recursions} with respect to $\lm_i$ ($i=1,2$). A measure $\tq$ solving \eqref{e:pair.frozen.recursions} corresponds to a root $\qoo$ of the function
$$\begin{array}{rl}
f(x)
\hspace{-6pt}&\equiv
x \,[ 1+(\lm_1-1)(1-{}_1\qo)^{d-1}+(\lm_2-1)(1-{}_2\qo)^{d-1} ]\\
&\qquad
- w(x)^{d-1} \,[ \lm_1\lm_2-x (\lm_1-1)(\lm_2-1) ]
\end{array}
$$
where $w(x)\equiv 1-{}_1\qo-{}_2\qo+x$, so that $w(\qoo)=\qgg$. We calculate
$$\begin{array}{rl}
f'(x)
\hspace{-6pt}&\equiv
 1+(\lm_1-1)(1-{}_1\qo)^{d-1}+(\lm_2-1)(1-{}_2\qo)^{d-1}\\
&\qquad- w(x)^{d-2}\,[
	(d-1)( \lm_1\lm_2-x(\lm_1-1)(\lm_2-1) )-w(x)\,(\lm_1-1)(\lm_2-1)
	].
\end{array}$$
Recall \eqref{e:lambda.of.q} that $\lm$ is monotone in $\qo\equiv\qo^\lm$, so in the regime $d^{1/10}\le\lm_i\le d^2$ we must have ${}_i\qo = y_i\logdbyd$ with $\slf{21}{20}\le y_i \le 3$, and consequently
$$(\lm_i-1)(1-{}_i\qo)^{d-1}
= {}_i\qo/(1-{}_i\qo) - (1-{}_i\qo)^{d-2}
= O( \logdbyd )
\quad\text{for $i=1,2$.}$$
If $x=\rho\logdbyd$ then $w(x)^{d-2}\asymp d^{\rho-y_1-y_2}$, so for $\rho\le\slf{9}{10}$ we find $f'(x)=1-O(d^{-1/20})$, which shows clearly that $f$ has the unique root $\qoo=({}_1\qo)({}_2\qo)$ in the stated regime.
\epf
\elem

\bpf[Proof of Propn.~\ref{p:indep.maximizer}]
As noted above, Cor.~\ref{c:boundary}\ref{c:boundary.b} implies that any maximizer $\bh$ of $\bPhit$ on $\simplextalind$ corresponds to a solution $h$ of the pair Bethe recursions \eqref{e:bp} with respect to some $\smash{\phi^{\lm_1}\otimes\phi^{\lm_2}}$. We now show that $h$ must satisfy the Bethe symmetries $\hh(\bm{io})=\hh(\bm{i}'\bm{o})$, where $\bm{i}$ or $\bm{i}'$ now indicates the incoming \emph{pair} of variable-to-clause messages, and $\bm{o}$ the outgoing pair of clause-to-variable messages.

We apply the same argument used to show the first-moment symmetries $\hh_{\eta'\eta}=\hh_{\eta''\eta}$ in the proof of Thm.~\ref{t:first.moment.exponent}: that is, we shall argue that the effect of changing the (pair) message $\bm{i}$ incoming to $\etreebip$ does not percolate down the tree. Again, if $\bm{i}$ and $\bm{i}'$ differ only by changing $\zro$ to $\free$ or vice versa in either coordinate, the effect does not propagate at all so we immediately find $\hh_{\bm{io}}=\hh_{\bm{i}'\bm{o}}$. Further, whenever the outgoing message is $\zro$ the value of $\hh$ does not depend on the incoming message. Analogously to the equivalences \eqref{e:bp.sym}, all remaining symmetries will follow by showing $\hh(\spint{\si}{\oo})=\hh(\spint{\si}{\fo})$ for any $\si$: as in the proof of Thm.~\ref{t:first.moment.exponent}, this follows by showing
$$d\,\hd(\spint{\si}{\fo})
\ll \hd(\spint{\si}{\zo})
\quad\text{and}\quad d\, \hd(\spint{\si}{\ff}) \ll \hd(\spint{\si}{\zf}) \quad\text{for all }\si\in\msg.$$
\bnm[1.]
\item The estimates $d\,\hd(\spint{\si}{\fo})\ll \hd(\spint{\si}{\zo})$ are all proved in essentially the same manner which is similar to the proof of $\hd_\fo/\hd_\zo\ll d^{-1}$ (see~\eqref{e:fo.zo.ratio}), so we sketch the case $\si=\fz$ (which is slightly delicate) and leave the rest to the reader.

By \eqref{e:bij} and the $\zro$-$\free$ symmetry noted above, $\hd(\spint{\si}{\fo}) / \hd(\spint{\si}{\zo})=\hbh(\spint{\si&\refl\si}{\fo&\of})/\hbh(\spint{\si&\refl\si}{\zo&\oz})$, the ratio between the number of $\free\zros$--$\zo$ edges to the number of $\free\zror$--$\zo$ edges. Let $\mu=n_\fo/n_\fz$ denote the proportion of $\fz$-variables which are matched to $\fo$-variables, and as before let $F_2$ denote the number of non-matching edges between $V_\od$ and $V_\fz$. Again using the $\zro$-$\free$ symmetry we calculate
\beq\label{e:fz.edges.from.od}
\begin{array}{l}
F_2 / n_\fz
= \mu\,\E[Y\giv Y\ge1]+(1-\mu)\,\E[Y\giv Y\ge2], \vspace{2pt}
\\
Y\sim\Bin(d-1,p)
\text{ where }
p\equiv \hd(\spint{\gf}{\oz}) / \hd(\spint{\gf}{\dotz}).
\end{array}
\eeq
Recall \eqref{e:restrict.K.F} that $F_2/n_\fz$
 must be close to $d\pi_\zo$; a slight variation of the argument following \eqref{e:derivative.conditioned.binomial} then gives $p\approx \pi_\zo$. We therefore obtain the desired bound
$$
\f{\hd(\spint{\fz}{\fo})}{\hd(\spint{\fz}{\zo})}
=\f{\mu\,\E[Y\I_{\set{Y=1}}\giv Y\ge1]
	+(1-\mu)\,\E[Y\I_{\set{Y=2}}\giv Y\ge2]}
{\mu\,\E[Y\I_{\set{Y>1}}\giv Y\ge1]
	+(1-\mu)\,\E[Y\I_{\set{Y>2}}\giv Y\ge2]}
\le d^{-3/2}.
$$

\item The estimates $d\,\hd(\spint{\si}{\ff})\ll\hd(\spint{\si}{\zf})$ can be proved by estimating the  ratios of edge densities $\hbh(\spint{\si&\refl\si}{\ff&\ff})/\hbh(\spint{\si&\refl\si}{\zf&\fz})$. These are relatively straightforward from the estimates of Lem.~\ref{l:ap.near.indep} and are left to the reader.
\enm
This proves the symmetries $\hh(\bm{io})=\hh(\bm{i}'\bm{o})$, so we conclude that $h$ must in fact correspond to a solution $\tq$ of the pair frozen model recursions \eqref{e:pair.frozen.recursions} via $9\,\hh(\bm{io})=\tq_{\bm{o}}$ (cf.~\eqref{e:hh.q}). It remains to verify that $\tq$ falls within the regime of Lem.~\ref{l:pair.frozen.recursions}. As before, let $K$ denote the number of $\zo$--$\oz$ edges: 
from \eqref{e:bij} and the $h$-to-$\tq$ correspondence,
\beq\label{e:oz.edges.from.zo}
K / n_\oz
	= \f{ \sum_{j\ge2} j \tbinom{d}{j}
		(\qgo  / \qgg)^j }
	{ \sum_{j\ge2} \tbinom{d}{j}
		(\qgo  / \qgg)^j },\quad
\text{so \eqref{e:restrict.K.F} implies 
$\qgo\approx \pi_\zo \qgg$
and $\qog\approx \pi_\oz \qgg$.}
\eeq
The parameter $p$ in \eqref{e:fz.edges.from.od} can be expressed as $\qgo/\qfzdot$, so \eqref{e:oz.edges.from.zo} implies $p\approx \pi_\oz$, therefore $\P(Y\ge2)\asymp1$ for $Y\sim\Bin(d-1,p)$. Another application of \eqref{e:bij} gives
$$
\f{n_\fo}{\numfz}
=\f{\mu}{1-\mu}
= \f{\qoo \, \P(Y\ge1)}{\qog \,\P(Y\ge2)}
\asymp \f{\qoo}{\qog},
$$
so Lem.~\ref{l:ap.near.indep}\ref{l:ap.near.indep.b} implies (crudely) that $\qoo/\qog \le d^{-1/2}$. From these estimates we may apply Lem.~\ref{l:pair.frozen.recursions} to conclude that $\smash{\tq=q^{\lm_1}\otimes q^{\lm_2}}$. By the one-to-one correspondence (Propn.~\ref{p:explicit}) between $\al$ and $q$ in this regime we conclude $\bh =\bhstaral$ as claimed.
\epf

\subsection{\emph{A priori} rigidity estimate}\label{ss:rigid}

In this section we analyze near-identical frozen model configurations to prove

\bppn\label{p:second.moment.id}
The contribution to $\E[\ZZ_{n\al}^2]$ from configurations $\bh\in\simplextalid$ is $\asymp_d\E\ZZ_{n\al}$.
\eppn

\begin{figure}[ht]
\includegraphics[trim=.8in .6in .4in .9in,clip]{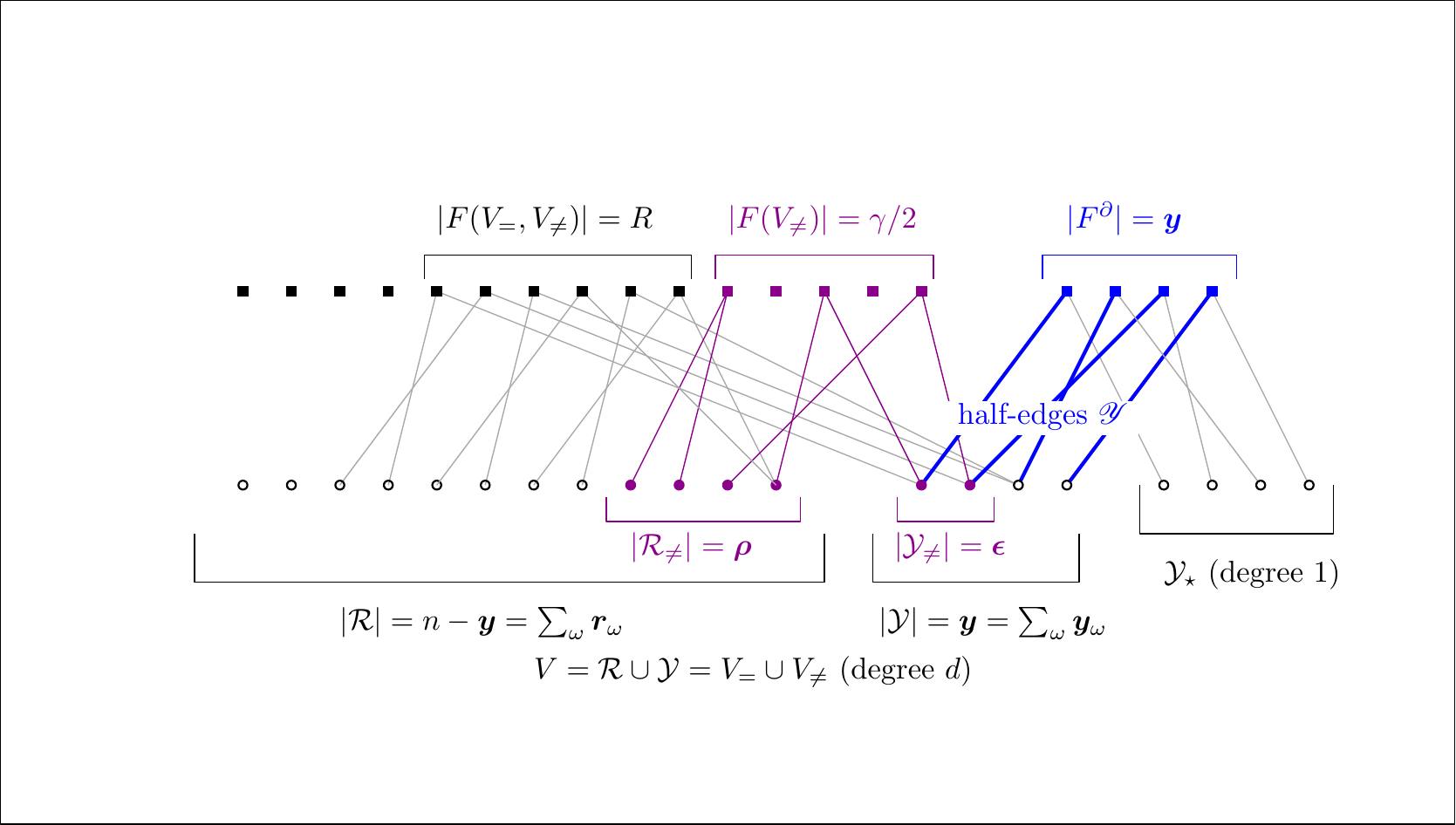}
\caption{Graph $\Gmin$ with dangling edges}
\label{f:dangling}
\end{figure}

The proof of Propn.~\ref{p:second.moment.id} is based on an {\it a~priori} estimate showing that frozen model configurations are sufficiently rigid that one typically does not find a large cluster of configurations near a given one. For application in our proof of the tightness of $\mis_n$ we shall prove this estimate for graphs drawn from the following slight generalization of the configuration model which allows for some unmatched edges (Fig.~\ref{f:dangling}). Let $V$ be a set of $n$ vertices, each incident to $d$ half-edges. Let $\extY$ be a disjoint set of $\by\lesssim_d\log n$ vertices, each incident to a single half-edge. Finally let $F$ be a set of $(nd+\by)/2$ clauses, each incident to $2$ half-edges. Let $\Gmin$ be the graph formed by taking a random matching between the half-edges incident to $\bdV\equiv V\cup \extY$ with the half-edges incident to $F$. We shall write $\bdF \equiv\pd \extY\subseteq F$, $\intY\equiv V\cap\pd\bdF$, and $\intR\equiv V\setminus \intY$.

\bdfn\label{d:frozen.dangling}
An auxiliary model configuration on $\Gmin$ is a message configuration $\usi^\pd$ on the edges of $\Gmin$ (including variable-to-clause messages from $\extY$ to $\bdF$) such that the configuration of messages incident to any variable or clause within $\Gmin$ is valid.
A frozen model configuration $\uz$ on $\Gmin$ is a spin configuration $\smash{\ueta\in\set{\zro,\one,\free}^{\bdV}}$ together with a subset $\match\subseteq F$ such that every $v\in V$ satisfies properties \eqref{d:frozen.matched.i}-\eqref{d:frozen.matched.iii} of Defn.~\ref{d:frozen.matched}\footnote{The properties need not be satisfied by vertices on $\extY$.}, and further the density of $\free$-variables in each coordinate is $\le\bemax$.
\edfn

Arguing as in Rmk.~\ref{r:vertex.aux}, there is a one-to-one correspondence between the frozen and auxiliary models on $\Gmin$. As in \eqref{e:frozen.pair.configuration}, we associate to a pair frozen configuration $\uz^i\equiv(\ueta^i,\match^i)$ ($i=1,2$) on $\Gmin$ a spin configuration $\smash{\uom\equiv\uom(\uz^1,\uz^2)\in\pairs^{\bdV}}$. As in \eqref{e:subset.spin}, for $\om\in\pairs$ and $S\subseteq\bdV$ we write $S_\om$ for the subset of vertices in $S$ with spin $\om$; in particular, write $S_{\ne}$ for the subset of vertices in $S$ with spin in $\smash{\pairsne\equiv\pairs\setminus\set{\zz,\oo,\ff}}$, and $\smash{S_=\equiv S\setminus S_{\ne}}$. Let
$$
\begin{array}{rlrlrl}
\intN_\om
	\hspace{-6pt}&\equiv |\intR_\om|,\quad
	&\extN_\om
	\hspace{-6pt}&\equiv|\intY_\om|,\quad
	&n_\om
	\hspace{-6pt}&\equiv |V_\om|=\intN_\om+\extN_\om,\\
\intNE
	\hspace{-6pt}&\equiv\intN_{\ne},\quad
	&\extNE
	\hspace{-6pt}&\equiv\extN_{\ne},\quad
	&\De
	\hspace{-6pt}&\equiv
	n_{\ne} =\intNE+\extNE.
\end{array}
$$
For $S,S'\subseteq\bdV$ write
$\edges(S,S')$ for the set of clauses $a\in F$ joining a vertex in $S$ to a vertex in $S'$; write $\edges(S)\equiv\edges(S,S)$ for the clauses internal to $S$.

\newcommand{\Ot}{{}_2\bm{O}}
\newcommand{\tauY}{\uta_\YY}

Write $\YY$ for the half-edges joining $\intY$ to $\bdF$. In the following, we fix boundary conditions given by a auxiliary pair configuration $\tauY\equiv(\usi^1_\YY,\usi^2_\YY)$ on $\YY$. Let $\ZZt^\pd[\pi|\tauY]$ denote the cardinality of the set $\Ot^\pd[\pi|\tauY]$ of pair frozen model configurations on $\Gmin$ which are consistent with $\tauY$ and have empirical measure
$\pi_\om=\intN_\om/\intN$ when restricted to $\intR$. We decompose
\beq\label{e:Z.boundary.decompose.A}
\ZZt^\pd[\pi|\tauY]
= \sum_{A\ge0} \ZZt^\pd_A[\pi|\tauY]
= \sum_{A\ge0}
	| \Ot^\pd_A[\pi|\tauY] |
\eeq
where $\Ot^\pd_A[\pi|\tauY]\subseteq \Ot^\pd[\pi|\tauY]$ denotes the subset of configurations which have $\gm/2\equiv|\edges(V_{\ne})|=\intNE+A$ internal edges among the unequal spins.\footnote{It is straightforward to see that $A\ge0$; we show a stronger inequality in Lem.~\ref{l:trees.w}.}
 We shall compare the expectation of $\ZZt^\pd_A[\pi|\tauY]$ with that of $\ZZ^\pd[\pi^1|\usi^1_\YY]$ ---
the number of frozen model configurations on $\Gmin$ which are consistent with $\usi^1_\YY$ and have empirical measure given by the projection $\pi^1$ of $\pi$ onto the first coordinate:
\beq\label{e:diff.int.notation}
|\intR|\pi^1_\eta
\equiv |\intR|\pi_{\eta\cdot}
= \intN_{\eta\cdot}
\equiv \intN_{\eta\eta} + \intNE_\eta
\quad(\eta\in\set{\zro,\one,\free}).
\eeq
Note $\ZZ^\pd[\pi^1|\usi^1_\YY]=\proj^1(\Ot^\pd[\pi|\tauY])$ where $\proj^1$ is the projection mapping $(\uz^1,\uz^2)\mapsto\uz^1$.

\bppn\label{p:rigidity}
There exists a small constant $\erho>0$ (uniform in $d\ge d_0$) such that for any empirical measure $\pi$ on $\pairs$ whose projections $\pi^1,\pi^2$ have intensities $\albd\le\al^1,\al^2\le\aubd$,
\[ \E[\ZZt^\pd[\pi|\tauY]]
\le
\f{ \E[\ZZ^\pd[\pi^1|\usi^1_\YY]]
	+ \E[\ZZ^\pd[\pi^2|\usi^2_\YY]] }{d^{(\intNE-\extNE)/10}}\]
provided $\De\equiv\intNE+\extNE\le n \erho \logdbyd$.

\bpf
Fix $\smash{\uom\in\pairs^{\bdV}}$ which is consistent with $\tauY$ and whose restriction to $\intR$ has empirical measure $\pi$, and write $\smash{\ueta^1\equiv\ueta^1(\uom)\in\set{\zro,\one,\free}^{\bdV}}$ for its projection onto the first coordinate. Write $\aone\equiv\pi_\od$, $\be\equiv\pi_\fd$. Write $\ue[\om]$ for the number of half-edges incident to $V_\om$ and not matched to a half-edge from $\extY$ (for example $\ue[\ne]=\intNE d+\extNE(d-1)$); the total number of such half-edges is $\ue\equiv\sum_\om \ue[\om]=nd-\by$. Let $R\equiv|\edges(V_=,V_{\ne})|=\ue[\ne]-\gm$.

\medskip\noindent\emph{Assignment of half-edges from $V_{\ne}$.} \\
We first estimate the probability of the event $F(V_{\ne},V_\oo)=\emptyset$: to this end, define
\[\DD_{A,B}(\uom)
\equiv \set{|\edges(V_{\ne})|=\gm/2\equiv\intNE+A,
|\edges(V_{\ne},V_\ff)|=B,\edges(V_{\ne},V_\oo)=\emptyset}.\]
Then, writing $\bin_{d,p}(j) \equiv \tbinom{d}{j} p^j(1-p)^{d-j}$ and recalling the notation \eqref{e:falling.factorial}, we have
\[\P[\DD_{A,B}(\uom)]
\le \underbrace{
	\f{ (\ue[\ne])_\gm\,(\gm-1)!! }{ \gm!\, \fdf{\ue}{\gm/2}  }
	}_{ \lesssim e^{O(\gm)} ( \De^2 d/(n\gm) )^{\gm/2} }
\times \underbrace{
	\f{ (\ue[\ff]+\ue[\zz])_R }{ \fdf{ \ue-\gm }{R} }
	}_{ \lesssim d^{ O(\De\erho) } e^{-\ue[\oo]R/\ue} }
\times \underbrace{
	\f{\tbinom{R}{B} (\ue[\ff])_B ( \ue[\zz] )_{R-B}}
		{ (\ue[\ff]+\ue[\zz])_R }
	}_{ \lesssim \bin_{R,2\be}(B) },\]
so altogether $\P[\DD_{A,B}(\uom)] \lesssim  e^{O(A)} d^{O(\De\erho)} ( \De^2 d/(n\gm) )^{\gm/2}e^{-\ue[\oo]R/\ue}\bin_{R,2\be}(B)$.

\medskip\noindent\emph{Assignment of remaining half-edges.} \\
Let $R'\equiv \ue[\ne]+R=2\ue[\ne]-\gm\le 2\De d$
denote the number of half-edges which are either incident to $V_{\ne}$ or matched with a half-edge from $V_{\ne}$. Let $G_=$ denote the induced subgraph on $V_=$, and let $\cJ_=(\uom)$ be the number of frozen model configurations on $G_=$ which are consistent with $\uom$. Meanwhile let $\cJ(\ueta^1)$ be the number of frozen model configurations on $\Gmin$ which are consistent with $\ueta^1$. We compare $\E[\cJ_=(\uom)\giv \DD_{A,B}(\uom)]$ to $\E[\cJ(\ueta^1)]$ by modifying the proof of Propn.~\ref{p:ap.compare.frozen.is}: the boundary $\extY$ causes some minor complication, but otherwise the calculation is much simplified by the restriction to vertices of equal spin.

Let $\Vucff$ denote the subset of vertices in $V_\ff$ not matched to $\extY$ under $\tauY$. Since there are by definition no edges between $\intR$ and $\extY$, clearly $\Vucff\cap\intR=\intR_\ff$,
so we write
$$\intYucff\equiv\Vucff\cap\intY,\quad
	\extNucff\equiv|\intYucff|,\quad
\nucff\equiv|\Vucff|\equiv\intN_\ff+\extNucff.$$
Similarly, let $\Vucfd$ denote the subset of vertices in $V_\fd$ not matched to $\extY$ under $\usi^1_\YY$: then
$\Vucfd\cap\intR=\intR_\fd$,
$\Vucfd\cap\intR_{=}=\intR_\ff$, and
$\Vucfd\cap\intY_= = \intYucff$, so we write
\[
\begin{array}{l}
\intYucfd\equiv\Vucfd\cap\intY,\quad
	\extNucfd\equiv|\intYucfd|,\quad
	\nucfd
	\equiv |\Vucfd| \equiv \intN_\fd+\extNucfd,\\
\extNEucfd\equiv|\Vucfd \cap \intY_{\ne}|
		= \extNucfd-\extNucff,\quad
	\Deucfd \equiv |\Vucfd \cap V_{\ne}|
	= \intNE_\free+\extNEucfd
	= \nucfd-\nucff.
\end{array}
\]
Consider the expected number of matchings on $\Vucff$
conditioned on $\DD_{A,B}(\uom)$, divided by the expected number of matchings on $\Vucfd$: this ratio is given by
$$
\f{(\prod_{v\in \Vucff}
	|\pd v\setminus (\extY\cup V_{\ne} )| )
	\,(\nucff-1)!! \,/\, \fdf{ \ue-R' }{ \nucff/2} }
{ (\prod_{v\in \Vucfd}
	|\pd v\setminus \extY|)
	\,(\nucfd-1)!!\,/\, \fdf{\ue }{ \nucfd/2}
}
\le \JJ_\free
\equiv\f{ \gfree(\nucff,\ue-R') }
	{ \gfree( \nucfd,\ue ) }\,
	[d/(d-1)]^{\extNEucfd},
$$
where $\gfree(n_\free,E)\equiv d^{n_\free}(n_\free-1)!! / \fdf{E}{n_\free/2}$ (cf.~\eqref{e:expected.matchings.f}), and the factor $[d/(d-1)]^{\extNEucfd}$ accounts for the vertices in $v\in\intYucfd\setminus\Vucff$ which contribute factors $d-1$ in the denominator on the left-hand side. We estimate $\JJ_\free \le e^{O(\De)}(d\be)^{-\Deucfd/2}$.

Because of the restriction to $V_=$, in the current setting the method of Propn.~\ref{p:ap.compare.frozen.is} reduces to a first-moment calculation, yielding
\[\f{\E[\cJ_=(\uom)\giv \DD_{A,B}(\uom)]}{\E[\cJ(\ueta^1)]}
\le \JJ_\free\times
	\underbrace{
	\f{ \gone( \ue-R'-\nucff,\ue[\oo],\ue[\ff]-B-\nucff ) }
	{ \gone( \ue-\nucfd,\ue[\od],\ue[\fd]-\nucfd) }}_{  \JJ_\one(B) }
	\times\JJ_\zro\]
where $\JJ_\zro$ is the (conditional) probability of all vertices in $V_\zz$ being forced, divided by the (conditional) probability of all vertices in $V_\zd$ being forced: explicitly, for $v\in V$ set $\ell^v=1$ if $v$ has a neighbor $w\in\extY$ with $\si^1_w=\one$, otherwise set $\ell^v=0$. 
Fix $0<\vth<1$ and define independent random variables $X^v\sim\Bin(|\pd v\setminus\extY|,\vth)$, with joint law $\P_\vth$: then
$$
\JJ_\zro
\le\f{ \P_\vth(X^i + \ell^i\ge2 \,\forall i\in V_\zz \setminus \pd V_{\ne}
	\giv \sum_{i\in V_\zz\setminus\pd V_{\ne}}
	X^i=\ue[\oo] ) }
	{ \P_\vth(X^i + \ell^i\ge2 \,\forall i\in V_\zd
	\giv \sum_{i\in V_\zd} X^i=\ue[\od] ) },
$$
where we have ignored the forcing constraint on the vertices in $V_\zz$ neighboring $V_{\ne}$. A trivial variation on the proof of Propn.~\ref{p:forcing} gives $\JJ_\zro=e^{O(\De)}$.

Recall $R'\equiv\ue[\ne]+R$; there remain $\ue[\zz]-R+B$ available half-edges from $V_\zz$. Therefore
$$
\JJ_\one(B)
=\f{ (\ue[\zz]-R+B)_{\ue[\oo]}
		\times\fdf{ \ue-\nucfd }{ \ue[\od] } }
	{\fdf{\ue[=]-R-\nucff }{ \ue[\oo] }
		\times ( \ue[\zd] )_{\ue[\od]}}
= d^{O(\De\erho)}
	\f{ \exp\{ \ue[\oo] (B+\ue[\ne])/\ue \} }
	{ \exp\{ \ue[\oo] ( \ue[\zd]-\ue[\zz] )/\ue \} }.
$$
Writing $\De_\zro\equiv n_\zro-n_\zz$ and recalling $\ue[\ne]-R=\gm=2(\intNE+A)$ then gives
\[
\sum_B \f{\bin_{R,2\be}(B) \JJ_\one(B)}{ \exp\{ \ue[\oo] R/\ue \} }
=d^{O(\De\erho)} e^{O(A)}
	\sum_B 
	\f{ \bin_{R,2\be}(B)   \exp\{ \ue[\oo] B/\ue \}}
	{\exp\{ \ue[\oo]( \ue[\zd]-\ue[\zz] )/\ue \}}
=\f{d^{O(\De\erho)} e^{O(A)}}
	{\exp\{\aone d \De_\zro \}}.
\]

\medskip\noindent
\emph{Combinatorial factors.} \\
The preceding estimates were for a fixed configuration $\uom\in\pairs^{\bdV}$ consistent with $\tauY,\pi$. Accounting for the permutations of $\uom$ and $\ueta^1$ gives (recalling $\intN=|\intR|$)
$$
R^1_A[\pi|\tauY]
\equiv
\f{ \E[\ZZt_A^\pd[ \pi|\tauY ]] }
	{ \E[\ZZ^\pd[ \pi^1 |\usi^1_\YY ]] }
\le e^{O( \De+A )}
\f{ \binom{\intN}{\intNE,\intN_\zz,\intN_\oo,\intN_\ff} }
	{ \binom{\intN}
		{\intN_\zd,\intN_\od,\intN_\fd} }
	\f{\sum_B
	\E[\cJ_=(\uom); \DD_{A,B}(\uom)]}{\E[\cJ(\ueta^1)]}.
$$ 
The ratio of multinomial coefficients is 
$\lesssim e^{O(\intN)} (n/\intNE)^{\intNE}
	(\aone)^{\intNE_\one}\be^{\intNE_\free}$.
Combining with the preceding estimates gives
\beq\label{e:near.id.pi.A}
R^1_A[\pi|\tauY]
\le d^{O(\De\erho)} e^{O(A)}
	(\intNE/\gm)^{\intNE}
	\Big( \smf{\De^2 d}{n\gm}\Big)^A
	\f{ d^{\intNE} (\aone)^{\intNE_\one}\be^{\intNE_\free} }
	{ \exp\{\aone d\De_\zro \} (d\be)^{\Deucfd/2} },
\eeq
where $\gm/2=\intNE+A$ and $\De\le \intNE e^{\extNE/\intNE}$ were used to bound $( \De^2/(\intNE\gm) )^{\intNE} \le e^{O(\intNE)}(\intNE/\gm)^{\intNE}$. Using $\aone \ge (3/2)\,\logdbyd$ and $\be\le\bemax$ then gives
$$R^1_A[\pi|\tauY]
\le
\f{ d^{O(\De\erho)} (\log d)^{O(\De)}
	(\intNE/\gm)^{\intNE}
( \f{\De^2 d}{n\gm})^A}
	{
	\exp\{ [\tf32\De_\zro-\intNE_\zro
		+\tf14(\intNE_\free-\extNEucfd)]\log d \}
	}.$$
Recall $\intN_\of \le n_\zf$ and $\intN_\fo\le n_\fz$; by symmetry we shall suppose $n_\zo+n_\zf \ge n_\oz+n_\fz$. Then
\[\tf32\De_\zro-\intNE_\zro\ge \tf12\De_\zro
\ge \tf14
	( n_\zo+n_\zf
	+ n_\oz+n_\fz )
\ge \tf18 ( \intNE-\intN_\ffne ),\]
so combining with $\intNE_\free-\extNEucfd\ge \intNE_\free-\extNE$ gives, for some constant $C$ uniform in $d$,
$$
\sum_{A\ge0}R^1_A[\pi|\tauY]
\le 
	\f{ d^{O(\De\erho)} (\log d)^{O(\De)}
	 }{ d^{(\intNE-\extNE)/8} }
\sum_{A\ge0} (\intNE/\gm)^{\intNE}
	\Big( \smf{C \De^2 d}{n\gm}\Big)^A
\le\f{  d^{O(\De\erho)} (\log d)^{O(\De)} }
		{ d^{(\intNE-\extNE)/8} }.$$
To see the last inequality, note the sum over $A\ge0$ is clearly $\lesssim1$ if $C\De^2 d\le 2n\gm$. If $C\De^2 d>2n\gm$, then recalling $\intNE \le \gm/2$ and optimizing over $\gm$ gives $(\intNE/\gm)^{\intNE} (C\De^2 d/(n\gm))^A\le d^{ O(\De\erho)}$. It follows that there exists a small constant $\erho$ (uniform in $d$) such that
$$
\sum_{A\ge0}R^1_A[\pi|\tauY]
\le d^{-(\intNE-\extNE)/10}
\quad\text{for }
	\De=\intNE+\extNE\le n\erho \logdbyd,$$
implying the result.
\epf
\eppn

\bpf[Proof of Propn.~\ref{p:second.moment.id}]
Follows from Propn.~\ref{p:rigidity} applied to our original random graph $\Gnd$ with no unmatched edges.
\epf

\bpf[Proof of Propn.~\ref{t:second.moment}]
Follows by combining Propns.~\ref{p:indep.maximizer}~and~\ref{p:second.moment.id}.
\epf

\section{Negative-definiteness of free energy Hessians}\label{s:nd}

In this section we prove Thm.~\ref{t:frozen.first.moment}
as well as
\bthm\label{t:frozen.second.moment}
It holds uniformly over $\albd\le\al\le\aubd$ that
$\E[\ZZ_{n\al}^2]
\lesssim_d
(\E\ZZ_{n\al})^2+\E\ZZ_{n\al}.
$
\ethm

\bppn\label{p:nd}
For $d\ge d_0$, the Hessians
$H\bPhi(\bhstaral)$
and $H(\bPhit)(\bhtstaral)$
as functions on $\simplexal$ and $\simplextal$ respectively
are negative-definite.
\eppn

The calculation of this section is similar to that of \cite[\S7]{dss-naesat}. Let $\bh\in\simplexalint$ with $\dbh$ and $\hbh$ both symmetric, and let $\bde$ be any signed measure on $\supp\vph$ (not necessarily symmetric) with $\bh+s\bde\in\simplexalint$ for sufficiently small $\abs{s}$. Then
\beq\label{e:rate.fn.second.derivative}
\pd_s^2\bPhi(\bh+s\bde)|_{s=0}
=-\angl{(\dbde/\dbh)^2}_{\dbh}
-(d/2)\, \angl{(\hbde/\hbh)^2}_{\hbh}
+d\, \angl{(\vde/\vh)^2}_{\vh},
\eeq
where $a/b$ denotes the vector given by coordinate-wise division of $a$ by $b$, and $\angl{\cdot}_h$ denotes integration with respect to measure $h$, e.g.\ $\angl{(\vde/\vh)^2}_{\vh}=\sum_\si \vde(\si)^2/\vh(\si)$. Consider maximizing \eqref{e:rate.fn.second.derivative} over $\bde$ subject to fixed marginals $\vde$: we find that the optimal $\hbde$ will be symmetric, with $\hbde(\si,\refl\si)=\vde(\si)=\vde(\refl\si)$. The optimal $\dbde$ will be of form $d\,\dbde(\dsi)=\dbh(\dsi)\sum_{i=1}^d \dchi_{\si_i}$ with $\dchi$ chosen to satisfy the margin constraint --- which, after a little algebra, becomes the system of equations
\beq\label{e:phi.delta}
\vH^{-1}\vde = d^{-1}\,[I+(d-1)\dM]\,\dchi
\eeq
where $\vH\equiv\diag(\vh)$ and $\dM$ denotes the stochastic matrix with entries
\beq\label{e:markov}
\dM_{\si,\si'}
\equiv
\f{1}{\vh(\si)}
	\sum_{\dsi} \dbh(\dsi)
		\Ind{(\si_1,\si_2)=(\si,\si')}.
\eeq
If such $\dchi$ exists, then the minimal value of $\angl{(\dbde/\dbh)^2}_{\dbh}$ subject to marginals $\vde$ is $\ip{\vde}{\dchi}$
(which clearly remains invariant under translations of $\dchi$ by vectors in the kernel of $I+(d-1)\dM$).

Throughout the following we take $\bh$ to be $\bhstaral$ (first moment) or $\bhtstaral$ (second moment).

\blem
The eigenvalues of $\dM$ counted with geometric multiplicity are
$$\eig(\dM)=(1,1,1,0,0,0,-(d-1)^{-1},\lm_1,\lm_2)$$
where $|\lm_1|\le d^{-1.9}$ and $0<|\lm_2-(d-1)^{-1}|\le d^{-1.2}$.
The right eigenvector $\bar{x}$ corresponding to eigenvalue $-(d-1)^{-1}$ is given by $\bar{x}_\si = (d-1)\Ind{\si=\oo}-\Ind{\si=\fz\text{ or }\ff}$.

\bpf\label{l:dM.evals}
The matrix $\dM$ is $\vh$-reversible and block diagonal with blocks
$\dM_\one$ acting on $\set{\oz,\of}$,
$\dM_\free$ acting on $\set{\oo,\fz,\ff}$, and
$\dM_\zro$ acting on $\set{\fo,\zo,\zz,\zf}$. We compute
$$\text{\footnotesize
$\dM_\one = 
\kbordermatrix{
	    & \oz                  & \of\\
	\oz & \tfrac{\qz}{1-\qo} & \tfrac{\qf}{1-\qo} \\
	\of & \tfrac{\qz}{1-\qo} & \tfrac{\qf}{1-\qo} \\
},\quad
\dM_\free
= \kbordermatrix{
    & \oo        & \fz  & \ff \\
\oo &          0 & \tfrac{\qz}{1-\qo}
                 & \tfrac{\qf}{1-\qo} \\
\fz & \tfrac{1}{d-1} & \tfrac{d-2}{d-1} \tfrac{\qz}{1-\qo}
                 & \tfrac{d-2}{d-1} \tfrac{\qf}{1-\qo} \\
\ff & \tfrac{1}{d-1} & \tfrac{d-2}{d-1} \tfrac{\qz}{1-\qo}
                 & \tfrac{d-2}{d-1} \tfrac{\qf}{1-\qo} \\
}$}$$
so $\eig(\dM_\one)=\set{1,0}$ and $\eig(\dM_\free)=\set{1,0,-(d-1)^{-1}}$; the right eigenvector of $\dM_\free$ corresponding to eigenvalue $-(d-1)^{-1}$ is $(d-1,-1,-1)$. We also compute
$$\text{\footnotesize$\dM_\zro
= \kbordermatrix{
    & \fo & \zo & \zz & \zf \\
\fo & \tfrac{1}{d-1} & 0 & \tfrac{d-2}{d-1} \tfrac{\qz}{1-\qo}
                     & \tfrac{d-2}{d-1} \tfrac{\qf}{1-\qo}\\
\zo & 0 & \ep+\qo(1-\ep) & \qz(1-\ep) & \qf(1-\ep)\\
\zz & \ep & \qo(1-\ep) & \qz(1-\ep) & \qf(1-\ep)\\
\zf & \ep & \qo(1-\ep) & \qz(1-\ep) & \qf(1-\ep)\\
}$}$$
where (using $\vh$-reversibility of $\dM$, or alternatively the frozen model recursion)
$$\ep
= \smf{(d-2)\qo^2 (1-\qo)^{d-3}}
	{ 1-[1+(d-2)\qo](1-\qo)^{d-2} }
= \smf{d-2}{d-1} \f{\qo \qf}{(1-\qo) \qz}
\asymp \qo\qf.$$
Write $\vH_\zro$ for the restriction of $\vH$ to $\set{\fo,\zo,\zz,\zf}$; then the symmetric matrix $\vH_\zro^{1/2}\dM_\zro\vH_\zro^{-1/2}$ has spectral decomposition $\sum_{i=1}^4 \lm_i e_i e_i^t$
with $\set{\lm_i}=\eig(\dM_\zro)$.
We can take
\[
\lm_3=1, \ e_3=\vh^{1/2}\quad\text{and}\quad
\lm_4=0, \ e_4 = (0,0, \qf(1-\qo)^{-1}, -\qz(1-\qo)^{-1})^{1/2}.\]
For the other two eigenvalues, consider the following ``almost'' eigenvalue equations:
\[
v_1 = ( (d-2)\qo, 1-\qo,-\qo,-\qo ),\quad
\dM_\zro v_1 = \ep( 0,1-\qo,(d-2)\qo,(d-2)\qo ),
\]
so $v_1$ is almost in the kernel of $\dM_\zro$; and
\[v_2=(1,0,0,0), \quad
	(\dM_\zro- (d-1)^{-1} I)v_2
	= \ep(0,0,1,1),
	\]
so $v_2$ is almost an eigenvector of $\dM_\zro$ with eigenvalue $(d-1)^{-1}$. We have $\ip{\vH_\zro^{1/2} v_1}{e_4}=0$ so $\vH_\zro^{1/2} v_1 = \sum_i a_i e_i$ with $a_4=0$. We then calculate
$$
\f{\sum_i \lm_i^2 a_i^2}{\sum_i a_i^2}
= \f{\|\vH_\zro^{1/2} \dM_\zro v_1\|^2}{ \|\vH_\zro^{1/2} v_1\|^2 }
=\f{\ep^2 \qz\qo(1-\qo) [ \qo + \tfrac{1-\qo}{(d-2)^2} ]}
	{\qo[ \qf\qo^2 + \tfrac{\qz(1-\qo)}{(d-2)^2} ]}
\asymp \qo (d\ep)^2
\asymp \qo (\qf\log d)^2.
$$
From Propn.~\ref{p:explicit} we have $\lm\ge d^{0.59}$,
and substituting into \eqref{e:frozen.recursions} gives $\qf \asymp 
d\qo^2/\lm \le d^{-1.58}$. It then follows from the above that at least one of the other two eigenvalues, say $\lm_1$, must have absolute value $\le d^{-2}$. Similarly, representing $\smash{\vH_\zro^{1/2}v_2=\sum_i b_i e_i}$ gives
$$\f{\sum_i (\lm_i-(d-1)^{-1})^2 b_i^2}{\sum_i b_i^2}
=\f{\|\vH_\zro^{1/2} (\dM_\zro-(d-1)^{-1}\cdot I) v_2\|^2}{ \|\vH_\zro^{1/2} v_2\|^2 }
=\f{\ep^2\qz(1-\qo)}{\qo\qf} \asymp \qo\qf \le d^{-2.5},
$$
so the last eigenvalue $\lm_2$ must satisfy $|\lm_2-(d-1)^{-1}| \le d^{-1.2}$. Note however that
$$\det [\dM_\zro -(d-1)^{-1}\cdot I]
=  \f{(d-2)\ep}{(d-1)^3} 
	[ d(\qo+\ep-\qo\ep)-(1+\qo+\ep-\qo\ep) ]\asymp \ep d\qo$$
so $\lm_2$ does not exactly equal $(d-1)^{-1}$.
\epf
\elem

\bpf[Proof of Propn.~\ref{p:nd}]
Consider the first moment Hessian $H\bPhi(\bhstaral)$ on $\simplexal$.
For convenience, rewrite \eqref{e:phi.delta} as $\vH^{-1/2}\vde = \dS\vH^{1/2}\dchi$ with $\dS$ the symmetric matrix $\dS\equiv d^{-1} \vH^{1/2}[I+(d-1)\dM]\vH^{-1/2}$. By the requirement that $\bh,\bh+\eta\bde$ both lie in $\simplexalint$ for $|\eta|$ small, the permissible vectors $\vde$ span a linear subspace $\ol{W}$.
In particular, necessarily $\ip{\vde}{\bar x}=0$ for the vector $\bar{x}$ which was determined in Lem.~\ref{l:dM.evals} to span the kernel of $\dL\equiv d^{-1}[I+(d-1)\dM]$. For such $\vde$, \eqref{e:phi.delta} is easily solved by considering the invariant action of $\dS$ on the subspace $(\vH^{1/2}\bar{x})^\perp$: if $U\in\R^{9\times 8}$ with columns giving an orthonormal basis for this subspace, then for all $\vde$ orthogonal to $\bar{x}$ (in particular, for all $\vde\in\ol{W}$) we may define
$\vH^{1/2}\dchi=U(U^t\dS U)^{-1}U^t \vH^{-1/2}\vde$,
consequently the maximal value of $\pd_\eta^2\bPhi(\bh+\eta\bde)|_{\eta=0}$ subject to marginals $\vde$ is
$$-\ip{\vde}{\dchi}+(d/2)\,\angl{(\vde/\vh)^2}_{\vh}
=-(U^t\vH^{-1/2}\vde)^t
\DOT{Q}
(U^t\vH^{-1/2}\vde),\quad
\DOT{Q}\equiv (U^t\dS U)^{-1}-(d/2)\,I$$
The eigenvalues of $\DOT{Q}$ are given by
$$\smf{d}{1+(d-1)\lm}-\smf{d}{2}\quad
	\text{for }
	\lm\in\eig(\dM) \setminus \set{ -(d-1)^{-1} },$$
so $\DOT{Q}$ is non-singular since we saw in Lem.~\ref{l:dM.evals} that $(d-1)^{-1}\notin\eig(\dM)$. Since we proved in Thm.~\ref{t:first.moment.exponent} that $\bhstaral\in\simplexalint$ is the global maximizer of $\bPhi$ on $\simplexal$, the restriction of the above quadratic form to the space of permissible $\vde$ (formally, to $U^t\vH^{-1/2}\ol{W}$) must be negative semi-definite, so by non-singularity we see that it is in fact negative-definite.

The proof for the second moment Hessian $H(\bPhit)(\bhtstaral)$ on $\simplextal$ is similar; note $\bhtstaral=\bhstaral\otimes\bhstaral$ implies $\dM_2=\dM\otimes\dM$, with eigenvalues $\set{\lm\lm' :\lm,\lm'\in\eig(\dM)}$. The kernel of $d^{-1}[I+(d-1)\dM_2]$ is spanned by vectors $\bar{x}\otimes\bar{y}$ or $\bar{y}\otimes\bar{x}$ with $\bar{x}$ as before and $\bar{y}$ any right eigenvector of $\dM_2$ with eigenvalue $1$, and again the permissible measures $\vde$ must be orthogonal to the kernel. Negative-definiteness then follows as above from the observation that $(d-1)^{-1}\notin\eig(\dM_2)$.
\epf

\blem\label{l:lin.bij}
For any $\si,\si'\in\msg$ there exists a signed integer measure $\bde\equiv\bde_{\si-\si'}=(\dbde,\hbde)$ with $\supp\bde\subseteq\supp\vph$ such that
$$0=\ip{\dbde}{1}=\ip{\hbde}{1}=\ii(\bde)\quad\text{and}
\quad \dH\dbde-\hH\hbde=\I_\si-\I_{\si'}.$$
The analogous condition holds for the support of the second-moment factors $\psit$.

\bpf
Define a graph on $\msg^2$ by placing an edge between $\tau$ and $\tau'$ if and only if there exists $\bde_{\tau-\tau'}$ satisfying the requirements in the statement of the lemma. We shall show this graph is connected (hence complete). By considering $\supp\dpsit=(\supp\dpsi)^2$ (see \eqref{e:psi}) it is clear that
$\set{\spint{\oz}{\si},\spint{\of}{\si}}$,
$\set{\spint{\fz}{\si},\spint{\ff}{\si}}$, and
$\set{\spint{\zz}{\si},\spint{\zo}{\si},\spint{\zf}{\si}}$
are connected subsets of $\msg^2$ for any $\si\in\msg$ --- for example, $\spint{\oz}{\si}$ and $\spint{\of}{\si}$ can be connected via
$$\bde_{\spint{\oz}{\si}-\spint{\of}{\si}}
=(\dbde,\hbde)
=\Big(
	\Ind{\dspint{\oz & \oz^{d-1}}{\si &\usi}}
	-\Ind{\dspint{\of & \oz^{d-1}}{\si & \usi}}
	,0\Big)$$
for any choice of $\usi\in\msg^{d-1}$ with $(\si,\usi)\in\supp\dpsi$. Then $\spint{\zz}{\si}$ and $\spint{\oz}{\si}$ can be connected via
$$\bde_{\spint{\zz}{\si}-\spint{\oz}{\si}}
=\bde_{\spint{\zo}{\refl\si}-\spint{\zz}{\refl\si}}
+\Big(0,\Ind{\dspint{\oz&\zo}{\si&\refl\si}}
	-\Ind{\dspint{\zz&\zz}{\si&\refl\si}}
	\Big);$$
the spins $\spint{\zz}{\si}$ and $\spint{\fz}{\si}$ can be similarly connected. We can then connect $\spint{\fo}{\si}$ and $\spint{\fz}{\si}$ via
\[\begin{array}{rl}
\bde_{\spint{\fo}{\si}-\spint{\zo}{\si}}
\hspace{-6pt}&=
	\Big(
	\Ind{\dspint{\fo&\fo&\zz&\zz^{d-3}}{\si&\si_2&\si_3&\vec\si}}
	-\Ind{\dspint{\zo&\zo&\zo&\zz^{d-3}}{\si&\si_2&\si_3&\vec\si}},
	\Ind{\dspint{\fo&\of}{\si_2&\refl\si_2}}
	-\Ind{\dspint{\zo&\oz}{\si_2&\refl\si_2}}
	\Big)\\
&\qquad{}+\bde_{\spint{\of}{\refl\si_2}-\spint{\oz}{\refl\si_2}}
+\bde_{\spint{\zo}{\si_3}-\spint{\zz}{\si_3}}
\end{array}
\]
for any choice of $\si_2$, $\si_3$, and $\usi$ such that $(\si,\si_2,\si_3,\usi)\in\supp\dpsi$. Lastly, to connect $\spint{\oo}{\si}$ and $\spint{\fo}{\si}$,
find $\si_2$ and $\usi\equiv(\si_3,\ldots,\si_d)$
with $(\si,\si_2,\usi)\in\supp\dpsi$, and take
\[\TS
\bde_{\spint{\oo}{\si}-\spint{\fo}{\si}}
=\Big(
\Ind{\dspint{\oo&\fz&\fz^{d-2}}{\si&\si_2&\usi}}
-\Ind{\dspint{\fo&\fo&\zz^{d-2}}{\si&\si_2&\usi}}
,0
\Big)
+\bde_{ \spint{\fo}{\si_2}-\spint{\fz}{\si_2} }
+\sum_{j=3}^d
	\bde_{ \spint{\zz}{\si_j}-\spint{\fz}{\si_j} }\]
This shows that for any $s,\si\in\msg$, $\spint{s}{\si}$ is connected to $\spint{\zz}{\si}$, which in turn by symmetry is connected to $\spint{\zz}{\zz}$. This proves connectivity of $\msg^2$;
connectivity of $\msg$
trivially follows.
\epf
\elem

\bpf[Proof of Thms.~\ref{t:frozen.first.moment}~and~\ref{t:frozen.second.moment}]
Recall (Defn.~\ref{d:simplex}) that $\E\ZZ_{n\al}$ is the sum of $\E\ZZ(\bh)$ over probability measures $\bh\equiv(\dbh,\hbh)$ on $\supp\psi$ such that $\bg\equiv(\dbg,\hbg)\equiv(n\dbh,\tf{nd}{2}\hbh)$ is integer-valued, and lies in the kernel of $H_{\simplex}\equiv\bpm \dH & -\hH \epm$ with intensity $\ii(\bg)=n\al$. Let $\E[\starZZ_{n\al}]$ denote the restriction of this sum to measures $\bg$ within a euclidean ball of radius $n^{1/2}\log n$ centered at $\bgstaral$; then Propn.~\ref{p:nd} implies $\smash{\E\ZZ_{n\al}=e^{o(1)}\,\E[\starZZ_{n\al}]}$.

Lem.~\ref{l:lin.bij} shows that the integer matrix $H_{\simplex}$ defines a surjection from $L'\equiv\set{\bde\in\R^{\supp\psi}:\ip{\dbde}{1}=\ip{\hbde}{1}=\ii(\bde)=0}$ to $\set{\vde\in\R^\msg:\ip{\vde}{1}=0}$, so we see that $L\equiv L'\cap (\ker H_{\simplex})\cap\Z^{\supp\psi}$ is an $(\bm{\DOT{s}}-2)$-dimensional lattice with spacings $\asymp_d1$. The measures $\bg$ contributing to $\E[\starZZ_{n\al}]$ are given by the intersection of the euclidean ball $\{ \|\bg-\bgstaral\|\le n^{1/2}\log n \}$ with an affine translation of $L$. The expansion \eqref{e:poly.correction} then shows that $n^{1/2}\E[\starZZ_{n\al}]$ defines a convergent Riemann sum, implying Thm.~\ref{t:frozen.first.moment}.

The same argument implies that the contribution
$\E[\ZZauxind_{n\al}]$ to the second moment
from measures in $\simplextalind$ is
$\asymp_d n^{-1}\exp\{ n\,[\bPhit(\bhtstaral)] \} \asymp_d (\E\ZZ_{n\al})^2$, so combining with Thm.~\ref{t:second.moment} and Propn.~\ref{p:second.moment.id} gives $\E[\ZZ_{n\al}^2]\asymp_d(\E\ZZ_{n\al})^2+\E\ZZ_{n\al}$.
\epf

\section{Constant fluctuations}

In this section we prove Thm.~\ref{t:lbd} by a variance decomposition of $\ZZ_{n\al}$. An application of our method in a simpler setting appears in \cite[\S8]{dss-naesat}, where the decomposition is introduced in greater detail. In what follows we write $\ZZ_{n\al}$ for the partition function of auxiliary model configurations $\usi$ with intensity $\ii(\usi)=n\al$ and empirical measure $\bh$ within distance $n^{-1/2}\log n$ of the first moment maximizer $\bhstaral$. In this section we restrict attention to the frozen model at intensities $n\albd\le n\al\le \misn$, and we write $f=O_d(g)$ or $f\lesssim_d g$ to indicate $f\le C(d)\,g$ where $C(d)$ can be chosen uniformly over $\al$ in this regime.

\bppn\label{p:var.log}
It holds uniformly over $\albd\le\al\le \misn/n$ that
$$\Var \log(\ZZ_{n\al}+\ep\,\E\ZZ_{n\al})
	\lesssim_d
	1 + 
	c_{d,\ep} (\E\ZZ_{n\al})^{-1}$$
\eppn

\bpf[Proof of Thm.~\ref{t:lbd}]
Write $\Lep_{n\al}\equiv\log(\ZZ_{n\al}+\ep\,\E\ZZ_{n\al})$. Since 
Thm.~\ref{t:frozen.first.moment}~and~Thm.~\ref{t:frozen.second.moment}
together
imply $\E[(\ZZ_{n\al})^2]\asymp_d(\E\ZZ_{n\al})^2$ uniformly over $\albd\le\al\le \misn/n$, there must exist a constant $\de\equiv\de(d)>0$ for which $\P(\ZZ_{n\al}\ge\de\,\E\ZZ_{n\al})\ge\de$, 
hence $\E\Lep_{n\al} \ge \de\log\de+(1-\de)\log\ep+\log\E\ZZ_{n\al}$. On the other hand, $\ZZ_{n\al}\le\ep\,\E\ZZ_{n\al}$ if and only if $\Lep_{n\al}\le\log(2\ep\,\E\ZZ_{n\al})$, which for small $\ep>0$ is much less than the lower bound on $\E\Lep_{n\al}$. Applying Chebychev's inequality and Propn.~\ref{p:var.log} therefore gives
$$\P(\ZZ_{n\al}\le\ep\,\E\ZZ_{n\al})
\le \f{\Var\Lep_{n\al}}
	{ ( \E\Lep_{n\al} - \log (2\ep\,\E\ZZ_{n\al}) )^2 }
\lesssim_d
\f{1 + c_{d,\ep} (\E\ZZ_{n\al})^{-1}}{ (\de\log(\de/\ep) -\log 2 )^2 }.
$$
With $c_\star$ as above (see~Cor.~\ref{c:first.moment.ubd}),
define $C_{d,\ep} \equiv 4c_\star[(\log c_{d,\ep})\vee0]$:
it then follows from Thm.~\ref{t:frozen.first.moment} that for all $\albd\le\al\le\misn-C_{d,\ep}$ that
$$c_{d,\ep} (\E\ZZ_{n\al})^{-1}\lesssim_d 1,\quad
\text{therefore }
\P(\ZZ_{n\al}=0)\lesssim ( \de\log(\de/\ep)-\log2 )^{-2}.
$$
The right-hand side tends to zero as $\ep$ decreases to zero, proving the theorem.
\epf

We prove Propn.~\ref{p:var.log} by controlling the increments of the Doob martingale of the random variable $\Lep_{n\al}\equiv\log(\ZZ_{n\al}+\ep\,\E\ZZ_{n\al})$  with respect to the edge-revealing filtration $(\filt_i)_{1\le i\le m}$, $m\equiv (nd/2)$, for the graph $G_n\sim\Gnd$. For $i\le m-2$ let $G^\circ$ denote the graph with edges $a_1,\ldots,a_{i-1}$. Randomly partition the remaining unmatched half-edges into disjoint sets $\KK$ and $\WW$ with $|\KK|=4$. We showed in \cite[\S8]{dss-naesat} by a coupling argument that for any $m'\le m-2$,
\beq\label{e:naesat.coupling.ubd}
\begin{array}{rl}
\Var\Lep_{n\al}
\hspace{-6pt}&\le
m'\max_{\Ak,\Aprime}\E[|\Vep_{n\al,m'}|^2]
	+\sum_{i=m'+1}^{m-1}\max_{\Ak,\Aprime}\E[|\Vep_{n\al,i}|^2]\\
&\qquad
\text{for }\Vep_{n\al,i}
	\equiv\E_{i-1}[\log(\ZZ_{n\al}+\ep\,\E\ZZ_{n\al})
	-\log(\acute\ZZ_{n\al}+\ep\,\E\ZZ_{n\al}) ]
\end{array}
\eeq
where $\E_{i-1}$ is expectation conditioned on the graph $G^\circ$; and $\ZZ_{n\al}\equiv\ZZ_{n\al}(\Ak),\acute\ZZ_{n\al}\equiv \ZZ_{n\al}(\Aprime)$ are the partition functions for coupled completions $G,\Gprime$ of $G^\circ$ to $d$-regular graphs: $G$ is formed by placing random edges $\Aw$ on $\WW$ and the given edges $\Ak$ on $\KK$, while $\Gprime$ is formed by placing the \emph{same} edges $\Aw$ on $\WW$ and the given edges $\Aprime$ on $\KK$. (Note that $\Vep_{n\al,i}$ depends on $\smash{\Ak,\Aprime}$ though we suppress this from the notation.)

Assume throughout that $m'\le m-2$ with $m-m'\lesssim_d\log n$. For such $m'$, we showed in \cite[\S8]{dss-naesat} that the sum appearing in \eqref{e:naesat.coupling.ubd} has two dominant components: the first is an ``independent-copies contribution'' coming from pair configurations with empirical measure near $\bhtstaral\equiv\bhstaral\otimes\bhstaral$; we showed that this contribution can be controlled under some abstract conditions (Lem.~\ref{l:ind.copies.conditions} and Cor.~\ref{c:indep} below). The other component is an ``identical-copies contribution'' coming from closely correlated pair configurations: this was controlled in \cite[\S8]{dss-naesat} under the assumption that the first moment is exponentially large, and the main work of this section is to control the identical-copies contribution assuming only that the first moment is bounded below by a large constant. Before turning to this we first show in \S\ref{ss:fluctuations.ind} that the independent-copies contribution is controlled by a straightforward application of the method of \cite[\S8]{dss-naesat}.

\subsection{Independent-copies contribution}
\label{ss:fluctuations.ind}

Let $B^\circ_\ell(\KK)$ denote the ball of graph distance $\ell$ about $\KK$ in the graph $G^\circ$; the leaves of $B^\circ_\ell(\KK)$ are half-edges incident to variables ($\ell$ odd) or clauses ($\ell$ even). Set
\beq\label{e:TT}
\TT=B^\circ_{2t'}(\KK),\quad
t' = t\wedge \min\set{\ell: B^\circ_{2\ell} \text{ has fewer than $|\KK|$ connected components} }.
\eeq
$\WW$ can only intersect $\TT$ in its leaves; and we shall let $\UU$ denote the leaves of $\TT$ without $\WW$.

Regarding $G$ as a $(d,2)$-regular graph of variables $V$ (degree $d$) and clauses $F$ (degree $2$), for $S\subseteq V\cup F$ and any configuration $\usi_S$ on the edges incident to the vertices $s\in S$, let
$$\Psi^\lm_{S,l}(\usi_S)
\equiv
\Ind{\ii(\usi_S)=l}
\prod_{v\in S\cap V} \dpsi^\lm(\dsi_v)
\prod_{a\in S\cap F} \hpsi^\lm(\ksi_a)
\in\set{0,\lm^l};
$$
we write $\Psi_{S,l}$ for the unweighted version given by taking $\lm=1$. Then, with $\YY\equiv\UU\cup\WW$,
$$\ZZ_{n\al}
=\sum_{r,j} \sum_{\usi_\YY} \Psi_{\Aw,j}(\usi_\WW)
	\ka_{r-j}(\usi_\UU)
	\ZZ^\pd_{n\al-r}[\usi_\YY]
=
\lm^{-n\al}
\sum_{r,j}
\sum_{\usi_\UU}
	\Psi^\lm_{\Aw,j}(\usi_\WW)
\ka^\lm_{r-j}(\usi_\UU)
\,
\ZZ^{\pd,\lm}_{n\al-r}[\usi_\YY]$$
where $\ka_l$ (resp.\ $\kaprime_l$) is the number of configurations on $\TT\cup\Ak$ (resp.\ $\TT\cup\Aprime$) with intensity $l$, with $\lm$-weighted versions denoted by $\ka^\lm_l,\kaprime^\lm_l$. Though we suppress it from the notation, let us emphasize that $\ka_l$ depends on $\TT$ and also on $\usi_\WW$, since $\TT$ may intersect $\WW$. Writing $\usi_\TT\sim\usi_\YY$ to indicate that $\usi_\TT$ and $\usi_\YY$ agree on the leaves of $\TT$,
$$\ka_l(\usi_\UU)
\equiv\ka_{\TT,l}(\usi_\UU|\usi_\WW)
\equiv\sum_{\usi_T\sim\usi_\YY} \Psi_{\TT\cup\Ak,l}(\usi_\TT).$$
Let $\ZZ^{\pd,\lm}_{n\al}$ denote the sum of $\ZZ^{\pd,\lm}_{n\al}[\usi_\YY]$ over all $\usi_\YY$, and $\ZZ^{\pd,\lm}_{n\al}[\cdot,\usi_\WW]$ the sum over the $\UU$-coordinates only. Write $\ETT$ for expectation over the law of the graph $\Gmin$ conditioned on $\TT$. As shown in \cite[\S8]{dss-naesat}, Lem.~\ref{l:lin.bij}
can be used to prove
\beq\label{e:unmatched.measure}
\ETT[
\ZZ^{\pd,\lm}_{n\al-r}[\usi_\UU,\usi_\WW]]
\sim
\bp^\lm(\usi_\UU)\,
	\ETT[\ZZ^{\pd,\lm}_{n\al-r}[\cdot,\usi_\WW]]
\sim
\bp^\lm(\usi_\YY)\,
	\ETT[\ZZ^{\pd,\lm}_{n\al-r}]
\eeq
where $\bp^\lm$ indicates the product measure with marginals $\hd^\lm$ coming from the Bethe solution corresponding to $\bhstarlm=\bhstaral$.
We then decompose $\ZZ_{n\al}$ according to a Fourier basis for $L^2(\msg^\UU,\bp^\lm)$:
take $(\fourier^\lm_1,\ldots,\fourier^\lm_{|\msg|})$ to be an orthonormal basis for $L^2(\msg,\hd^\lm)$ with $\fourier^\lm_1\equiv1$. Then the functions
$\fourier^\lm_{\uss}(\usi_\UU)
\equiv\prod_{u\in\UU} \fourier^\lm_{s(u)}(\si_u)$
($\uss\in[|\msg|]^\UU$)
form an orthonormal basis for $L^2(\msg^\UU,\bp^\lm)$. By Plancherel's identity
\beq\label{e:plancherel}
\ZZ_{n\al}
=\lm^{-n\al}
	\sum_{j,l}
	\sum_{\usi_\WW} \Psi^\lm_{\Aw,j}(\usi_\WW)
	\sum_{\uss}
	(\ka^\lm_l)^\wedge_{\uss}
	(\mathbf{F}^\lm_{n\al-j-l})^\wedge_{\uss}
\eeq
where $^\wedge$ indicates the Fourier transform with respect to the basis $\fourier^\lm$, and
$$
\mathbf{F}^\lm_{n\al-r}(\usi_\UU)
\equiv
\mathbf{F}^\lm_{n\al-r}(\usi_\UU|\usi_\WW)
\equiv
\f{\ZZ^{\pd,\lm}_{n\al-r}[\usi_\UU,\usi_\WW]}{\bp^\lm(\usi_\UU)}.
$$
By the method of \cite[\S8]{dss-naesat} applied to the decomposition \eqref{e:plancherel}, the independent-copies contribution to $\Var\Lep_{n\al}$ is controlled subject to a few conditions which are proved in the following lemma. For $\uss\in[|\msg|]^\UU$ write $|\uss|\equiv|\set{u\in\UU:s_u\ne1}|$, and write $\emptyset$ for the identically-$1$ vector in $[|\msg|]^\UU$. 

\newcommand{\emptykalm}[1]{(\overline{\ka}^\lm_{#1})^\wedge_\emptyset}

\blem\label{l:ind.copies.conditions}
The auxiliary model satisfies the following:
\bnm[(a)]
\item \label{l:ind.copies.conditions.b}
(Symmetries) On the event $\mathbf{T}$ that $\TT$ consists of $4$ disjoint tree components with $\TT\cap\WW=\emptyset$, we have $(\ka^\lm_l)^\wedge_{\uss}=(\kaprime^\lm_l)^\wedge_{\uss}$ for all $|\uss|\le1$; and the zeroth order Fourier coefficient $(\ka^\lm_l)^\wedge_\emptyset = \sum_{\usi_\UU} \ka^\lm_l(\usi_\UU) \bp^\lm(\usi_\UU)$ restricted to the event $\mathbf{T}$ takes a constant value $\emptykalm{l}$ which does not depend on the edges $\Ak$. On the event $\mathbf{C}^\circ$ that $\TT$ either contains a single cycle or has a single intersection with $\WW$ (but not both), we have $(\ka^\lm_l)^\wedge_\emptyset=(\kaprime^\lm_l)^\wedge_\emptyset$.

\item \label{l:ind.copies.conditions.c}
(Correlation decay) For $\al\ge\albd$, the Gibbs measure $\nu$ corresponding to $\bhstaral$ on the infinite $d$-regular tree has correlation decay at rate faster than square root of the branching rate: if edges $e,e'$ are separated by $t$ intermediate vertices then
$|\Cov_\nu(f(\si_e),g(\si_{e'}))|
	\lesssim_d
	\|fg\|_\infty d^{-(3/2)t}
	\ll \|fg\|_\infty (d-1)^{-t/2}$.
\enm

\bpf
\eqref{l:ind.copies.conditions.b}
Follows by the argument of \cite[\S8]{dss-naesat} (simpler in the current setting since random literals are not involved).

\medskip\noindent
\eqref{l:ind.copies.conditions.c}
Let $N$ denote the event that for some intermediate vertex $v$ along the path $\gm$ joining $e,e'$ in the tree, there exist at least two vertices $w\in\pd v\setminus\gm$ with $\si_{w\to v}=\one$. 
The spins $\si_e,\si_{e'}$ are exactly independent on the event $N$, which has probability $\nu(N)\ge1-d^{-(3/2)t}$.
\epf
\elem

As shown in \cite[\S8]{dss-naesat}, the conditions of Lems.~\ref{l:lin.bij}~and~\ref{l:ind.copies.conditions} taken together give control over the independent-copies component of $\Var\Lep_{n\al}$. Explicitly, let
\[\wt\ka_l\equiv\ka_l+\kaprime_l+ (\ka_l+\kaprime_l)^\wedge_\emptyset
=\ka_l+\kaprime_l+ \lm^{-l} (\ka^\lm_l+\kaprime^\lm_l)^\wedge_\emptyset,\]
and let $\ol\psi_j(\usi_\WW)$ denote the expectation of $\Psi_{\Aw,j}(\usi_\WW)$ over the random edges $\Aw$. Then define
\beq\label{e:vid}
\vid^\TT_{n\al,i}
\equiv
\sum_{r^1,r^2}
\sum_{j^1,j^2}
\sum_{\tauY}
\ol\psi_{j^1}(\usi^1_\WW)
\ol\psi_{j^2}(\usi^2_\WW)
\wt\ka_{r^1-j^1}(\usi^1_\UU)
\wt\ka_{r^2-j^2}(\usi^2_\UU)
\idZZ^\pd_{
	n\al-(r^1,r^2)}[\tauY]
\eeq
where $\smash{\idZZ^\pd_{n\al-(r^1,r^2)}}$ refers to the contribution to the pair partition function on $\Gmin\equiv G^\circ\setminus\TT$ from empirical measures within distance $n^{-1/5}$ of $\bhidal$, with intensity $n\al-r^i$ in the $i$-th coordinate for $i=1,2$. Then \cite[\S8]{dss-naesat} together with Lems.~\ref{l:lin.bij}~and~\ref{l:ind.copies.conditions} implies the following

\bcor\label{c:indep}
For $\albd\le\al\le\misn/n$ and
$T$ as in \eqref{e:TT} with $t\equiv t(d,\ep)= 4\log_d(1/\ep)$,
$$\Var\Lep_{n\al}
\lesssim_d
1+o_\ep(1)
+ \underbrace{
	\f{m'\E\vid^\TT_{n\al,m'} 
	+ \sum_{i=m'+1}^{m-2} \E\vid^\TT_{n\al,i} }
	{(\ep\,\E\ZZ_{n\al})^2}}_{\text{identical-copies component}}$$
for any $m'$ with $m-m'\lesssim_d\log n$.
\ecor

The remainder of this section is devoted to the handling of the identical-copies component. Let us first note that Lem.~\ref{l:lin.bij} can also be used to see how \eqref{e:unmatched.measure} varies with $n\al-r$: it is clear from the lemma that we can find signed integer measures $\bde^\mathrm{v}$, $\bde^\mathrm{c}$, $\bde^\ii$, all on the support $\supp\vph$ of the first-moment auxiliary model factors \eqref{e:psi}, such that
\beq\label{e:surj}
\begin{array}{rlrlrlrl} 
\ip{\dbde^\mathrm{v}}{1}
	\hspace{-6pt}&=1, &
\ip{\hbde^\mathrm{v}}{1}
	\hspace{-6pt}&=0, &
\ii(\hbde^\mathrm{v})
	\hspace{-6pt}&=0, &
\dH\dbde^\mathrm{v}-\hH\hbde^\mathrm{v}
	\hspace{-6pt}&=d\,\I_\zz\\
\ip{\dbde^\mathrm{c}}{1}
	\hspace{-6pt}&=0, &
\ip{\dbde^\mathrm{c}}{1}
	\hspace{-6pt}&=1, &
\ii(\hbde^\mathrm{c})
	\hspace{-6pt}&=0, &
\dH\dbde^\mathrm{c}-\hH\hbde^\mathrm{c}
	\hspace{-6pt}&=2\,\I_\zz\\
\ip{\dbde^\ii}{1}
	\hspace{-6pt}&=0, &
\ip{\hbde^\ii}{1}
	\hspace{-6pt}&=0, &	
\ii(\hbde^\ii)
	\hspace{-6pt}&=1, &	
\dH\dbde^\mathrm{c}-\hH\hbde^\ii
	\hspace{-6pt}&=0.
\end{array}\eeq
For example, taking $\bde_{\zo-\oz}$ as given by Lem.~\ref{l:lin.bij}, we can take $\bde^\ii$ to be
\[\bde^\ii
=( \Ind{ (\oz^d) }
	-\Ind{ (\zo^d)} ,0) + 
	d\,\bde_{ \zo-\oz }.\]
If $\Gmin$ has $n^\pd$ variables and $m^\pd$ clauses, then set
\[\bde^{\usi_\YY}_r
=(n-n^\pd)\bde^\mathrm{v}+ (m-m^\pd)\bde^\mathrm{c}
	+ \sum_{y\in\YY} \bde_{\si_y-\zz}
	+ r \bde^\ii.\]
Away from the simplex boundary, empirical measures contributing to $\ZZ^\pd_{n\al-r}[\usi_\YY]$ can be parametrized as $\bg-\bde^{\usi_\YY}_r$ where $\bg$ runs over empirical measures contributing to $\ZZ_{n\al}$. Write $\usi_\YY\equiv(\outm,\inm)$ where $\outm\in\set{\zro,\one,\free}^\YY$ are the outgoing messages of $\usi^i_\YY$, and $\inm\in\set{\zro,\one,\free}^\YY$ the incoming. Then
\beq\label{e:z.dangling.ones.comparison}
\f{\E\ZZ^\pd_{n\al-r}[\usi_\YY]}
	{\lm^r \E\ZZ_{n\al}}
=\f{\E\ZZ^{\pd,\lm}_{n\al-r}[\usi_\YY]}
	{\E\ZZ^\lm_{n\al}}
\sim
\bm{c}^{-1}\,\bp^\lm(\usi_\YY)
=\f{\bm{c}^{-1}\,\bq(\outm) \lm^{-|\set{y:\si_y=\oo}|}}
	{ [ 3-(1-\lm^{-1})\qo ]^{|\YY|} }.
\eeq
for $\bm{c}$ a proportionality constant depending only on $n-n^\pd$ and $m-m^\pd$, and $\bq$ the product measure with marginals $q\equiv (\qz,\qo,\qf)$.

\subsection{Refined rigidity estimate}

\newcommand{\trw}{\mathfrak{T}}

In this subsection we prove a refined version of Propn.~\ref{p:rigidity} for small $\De$.

\bdfn\label{d:trw}
Given a pair frozen configuration $(\uz^1,\uz^2)$ on $\Gmin$, let $\trw_\of$ denote the collection of tree components of $\bdV_{\ne}$ which intersect $\intY$ in a single $\of$-vertex with matched partner in $\extY$, and consist
entirely of $\of$-$\zf$ matched pairs. Let $\trw_\zf$ denote the collection of tree components of $\bdV_{\ne}$ which intersect $\intY$ in a single $\zf$-vertex with matched partner in $\intR$, and consist entirely of $\of$-$\zf$ matched pairs except for one additional vertex in $\extY$. Define symmetrically $\trw_\fo$ and $\trw_\fz$, and let $\trw\equiv\trw(\uz^1,\uz^2)$ denote the union of $\trw_\of$, $\trw_\zf$, $\trw_\fo$, $\trw_\fz$.
\edfn

\blem
\label{l:trees.w}
Let $(\uz^1,\uz^2)$ be a frozen model pair configuration on $\Gmin$.
\bnm[(a)]
\item \label{l:trees.w.a}
Any tree component $T$ of $\bdV_{\ne}$ with $|T\cap\extY|\le1$ must be in $\trw$.

\item \label{l:trees.w.b}
 With $|\edges(V_{\ne})|=\gm/2\equiv\intNE+A$ as before, we have
\beq\label{e:Aprime}
A''\equiv A-\tf12|\intY_{\ne}\setminus\trw|\ge0.
\eeq
Further, let $\numzf$ denote the
cardinality of the set $\Vumzf$
of vertices in $V_\zf$ whose matched partner does not have spin $\of$; symmetrically define $\numfz=|\Vumfz|$.Then
\beq\label{e:Aprime.lbd}
2A''\ge
[|\edges(V_\zo\cup V_\oz)|-2\intN_\zo-\extN_\zo]
+[|\edges(V_\zo\cup V_\oz)|-2\intN_{\oz}-\extN_{\oz}]
+\numzf+\numfz
\eeq
with both $|\edges(V_\zo\cup V_\oz)|-2\intN_\zo-\extN_\zo$
 and $|\edges(V_\zo\cup V_\oz)|-2\intN_{\oz}-\extN_{\oz}$ non-negative.
\enm

\bpf
Recall Defn.~\ref{d:frozen.dangling} that $\uz\equiv(\match,\ueta)$ is a frozen model configuration on $\Gmin$ if and only if every
vertex in $V=\bdV\setminus\extY$ satisfies properties \eqref{d:frozen.matched.i}-\eqref{d:frozen.matched.iii} of Defn.~\ref{d:frozen.matched}; the properties need not be satisfied on $\extY$. 

\newcommand{\cmin}{\mathfrak{C}}
\newcommand{\rmin}{\mathfrak{R}}

\medskip\noindent\eqref{l:trees.w.a}
Let $\cmin$ be a component of $\bdV_{\ne}$, and write $\chi\equiv\chi_\cmin\equiv|\edges(\cmin)|-|\cmin|$. Assume $\cmin\setminus\extY\ne\emptyset$ since the result is otherwise trivial. Since $\bdV_{\ne}$ can have no neighbors in $\bdV_\oo$, every $v\in\cmin\cap V$ must satisfy properties \eqref{d:frozen.matched.i}-\eqref{d:frozen.matched.iii} \emph{within $\cmin$ alone} (in both coordinates), i.e.\ in the absence of the rest of the graph. Write $d_v(\cmin)\equiv|\edges(v,\cmin)|$, and note that if $v_0$ is a leaf vertex of $\cmin$ (meaning $d_v(\cmin)=1$) then $v_0\in V_\of\cup V_\fo \cup \extY$.

Suppose then that $v_0\in V_\of$ is a leaf vertex of $\cmin$; it must be paired under $\match^2$ to a partner $v_1\in\cmin$ with spin $\zf$. If $v_1\notin\extY$, then (by property \eqref{d:frozen.matched.ii}) it must have at least one other neighbor $v_2\in\cmin$ with spin $\of$ or $\oz$. Since $v_0$ is a leaf and $v_1$ is a $\zf$-vertex matched to $v_0$, the presence of these two vertices in $\cmin$ is not needed to ensure that the other vertices satisfy the properties \eqref{d:frozen.matched.i}-\eqref{d:frozen.matched.iii}. Therefore, in the case that $v_1$ has \emph{exactly two} neighbors in $\cmin$ with each joined to $v_1$ by a single edge, we may choose to remove both $v_0,v_1$ from $\cmin$ causing no change in $\chi$. We iterate this procedure on leaf vertices in $V_\of$ or $V_\fo$ to produce a reduced graph $\rmin\equiv\rmin(\cmin)\subseteq\cmin$ which is one of the following:
\bnm[{\sc a.}]
\item \label{reduc.a}
$\rmin=\set{v_0,v_1}$ where $v_0\in\intY$ has spin $\of$ or $\fo$, and $v_1\in\extY$ is its matched partner.
\item \label{reduc.b}
$\rmin=\set{v_0}$ where $v_0\in\extY$ has spin $\oz$, $\of$, $\zo$, or $\fo$.
\item \label{reduc.c}
$\rmin\setminus\extY\ne\emptyset$; and
any $v\in V$ which is a leaf vertex of $\rmin$ is either a $\of$-vertex matched to a $\zf$-vertex $u$ with $d_u(\rmin)\ge3$, or symmetrically an $\fo$-vertex matched to an $\fz$-vertex $w$ with $d_w(\rmin)\ge3$.

\enm
In cases {\sc \ref{reduc.a}} and {\sc \ref{reduc.b}}, $\chi_{\cmin}=\chi_{\rmin}=-1$ so $\cmin$ is a tree consisting entirely of $\of$-$\zf$ and $\fo$-$\fz$ matched pairs except in case {\sc \ref{reduc.b}} for the vertex on $\extY$. We claim further that all the matched pairs must be of the same type: indeed, since the only edges between $\bdV_\of\cup \bdV_\zf$ and $\bdV_\fo\cup \bdV_\fz$ must occur between $\bdV_\zf$ and $\bdV_\fz$, every $v\in\cmin\setminus V_\fd$ satisfies properties \eqref{d:frozen.matched.i}-\eqref{d:frozen.matched.iii} in the absence of $V_\fd$. Therefore apply the leaf-removal procedure separately to each component $\mathfrak{c}$ of  $\cmin\setminus V_\fd$ or of $\cmin\setminus V_\dotf$: if $\rmin(\mathfrak{c})\ne\emptyset$ then it contains some $v\in\extY$ which must also appear in $\rmin=\rmin(\cmin)$. All the $\mathfrak{c}$ are disjoint, and $|\rmin\cap\extY|=1$ by assumption, so only a single component $\mathfrak{c}$ can be non-empty. This proves the claim, and we conclude $\cmin\in\trw_\of\cup\trw_\fo$ in case {\sc \ref{reduc.a}}, $\cmin\in\trw_\zf\cup\trw_\fz$ in case {\sc \ref{reduc.b}} (since we assumed $\cmin\setminus\extY\ne\emptyset$).

In case {\sc \ref{reduc.c}}, write $m_v\equiv2-\Ind{v\in\extY}$, and consider
\beq\label{e:d.minus.m}
\sum_{v\in\rmin}(d_v(\rmin)-m_v)=2\chi_{\cmin}+|\rmin\cap\extY|
\eeq
Any $v\in\rmin$ with $d_v(\rmin)-m_v<0$ must be in $V_\of\cup V_\fo$, and is compensated in the above sum by a matched partner $u\in\rmin$ with $d_u(\rmin)\ge3$. Therefore the sum restricted only to vertices belonging in $\fo$-$\fz$ or $\of$-$\zf$ matched pairs inside $\rmin$ is non-negative. Any other vertex gives a non-negative contribution to the sum. Thus $\chi_\cmin \ge-\tf12|\rmin\cap\extY|\ge -\tf12|\cmin\cap\extY|$, so in particular if $|\cmin\cap\extY|=1$ then (recalling we have assumed $\cmin\setminus\extY\ne\emptyset$) $\cmin$ cannot be a tree.

\medskip\noindent
\eqref{l:trees.w.b} Let $\mathbf{R}$ denote the union of the $\rmin(\cmin)$ over the components $\cmin$ of $\bdV_{\ne}\setminus\trw$ with $\cmin\setminus\extY\ne\emptyset$: it follows from the proof of \eqref{l:trees.w.a} that
$$0\le
\sum_{v\in\mathbf{R}}(d_v(\mathbf{R})-m_v)
=2\chi_{\mathbf{R}\setminus\extY}
	+|\mathbf{R}\cap\extY|
\le 
2A-|\intY_{\ne}\setminus\trw|
=2A''
\le 2A.$$
We furthermore decompose
\[2A''
\ge
\underbrace{\sum_{v\in V_\zo}(
	d_v(V_\oz)-m_v)}
	_{=
	\edges(V_\zo, V_\oz)-2\intN_\zo-\extN_\zo\ge0}
+\underbrace{\sum_{v\in V_{\oz}}
	(d_v(V_\zo)-m_v)}
	_{=
	\edges(V_\zo, V_\oz)
	-2\intN_\oz-\extN_\oz\ge0}
+\underbrace{\sum_{v\in V^\bullet_\zf}
	(d_v(\mathbf{R})-m_v)
	}_{\ge |V^\bullet_\zf|=n^\bullet_\zf }
+\underbrace{\sum_{v\in V^\bullet_\fz}(d_v(\mathbf{R})-m_v)
	}_{\ge |V^\bullet_\fz|=n^\bullet_\fz },\]
concluding the proof.
\epf
\elem

Suppose $(\uz^1,\uz^2)$ is a pair frozen model configuration on $\Gmin$ which is consistent with auxiliary configuration $\tauY\equiv(\usi^1_\YY,\usi^2_\YY)$ on the half-edges $\YY$ joining $\intY$ to the clauses $\bdF=\edges(\intY,\extY)$. Write $\tauY\equiv(\outm,\inm)\equiv(\outm^i,\inm^i)_{i=1}^2$ where $\outm^i\in\set{\zro,\one,\free}^\YY$ are the outgoing messages of $\usi^i_\YY$, and $\inm^i\in\set{\zro,\one,\free}^\YY$ the incoming. If $y\in\YY$ is at the root of a tree $T\in\trw$ then $\bm{o}_y\equiv\bm{o}^1_y\bm{o}^2_y$ is either $\oo$ ($T\in\trw_\of\cup\trw_\fo$) or $\ff$ ($T\in\trw_\zf\cup\trw_\fz$). Further, $v\in\intY$ is incident to $y\in\YY$ with outgoing messages $\bm{o}^1_y\ne\bm{o}^2_y$ only if either $v\in\intY_{\ne}\setminus\trw$, or $v\in\intY_\zz$ but is neighboring to $\intR_{\ne}$. Thus the Hamming distance $\hamming(\outm)$ between $\outm^1$ and $\outm^2$ satisfies
\beq\label{e:w.minus}
\hamming(\outm)
\le|\intY_{\ne}\setminus\trw|+k'
\quad\text{where }
k'\equiv|\set{v\in\intY_\zz:
	\edges(V_{\ne},v)\ne\emptyset}|.
\eeq

\newcommand{\JW}{\bm{J}}
\newcommand{\IU}{\bm{I}}

\blem\label{l:intensity}
Suppose $(\uz^1,\uz^2)$ is a pair frozen model configuration on $\Gmin$ contributing to $\idZZ^\pd_{n\al-(r^1,r^2)}[\tauY]$. For $i=1,2$ let $2j^i\equiv|\set{w\in\WW:\si^i_w=\oo}|$, $l^i\equiv|\set{u\in\UU:\si^i_u=\oo}|$, and $s^i\equiv r^i - j^i$. Then
$$|\JW-\IU|
	\le\hamming(\outm)+2A.
\quad\text{for }
	\JW\equiv j^1-j^2, \
	\IU\equiv s^1-l^1-(s^2-l^2).$$

\bpf
For $\om\in\set{\fo,\fz,\of,\zf,\ff}$ let $\bm{e}_\om$ count the number of vertices in $\intY$ with spin $\om$ and matched partner in $\extY$. For $i=1,2$ let $\bm{e}^i_{\ne}$ count the number of $\ffne$-vertices in $\intY$ whose matched partner in the $i$-th coordinate is in $\extY$. Note $\bm{e}_\fo$ also counts the set of $y\in\YY$ with $\bm{o}^1_y=\one=\bm{o}^2_y$ but $\bm{i}^1_y=\one\ne\bm{i}^2_y$, so
\beq\label{e:unmatched.hamming}
|(2j^1+l^1)-(2j^2+l^2)
	-(\bm{e}_\fo-\bm{e}_\of)|\le\hamming(\outm).
\eeq
Meanwhile $n\al-r^1 = n_\od +\tf12( n_\fd - \bm{e}_\fz-\bm{e}_\fo-\bm{e}_\ff-\bm{e}^1_{\ne} )$ and similarly for $n\al-r^2$. Recall $\numzf$ denotes the number of vertices in $V_\zf$ with matched partner in $V_\gf\cup\extY_\gf$, and let $\numzfint\le\numzf$ denote the number of vertices in $V_\zf$ with matched partner in $V_\gf$: then $n_\zf-n_\of=\bm{e}_\zf-\bm{e}_\of+\numzfint$, so taking the difference between the $n\al-r^i$ gives
$$r^1-r^2
=n_\zo-n_\oz
+\tf12( \bm{e}^1_{\ne}-\bm{e}^2_{\ne}
	+ \numzfint- \numfzint)
	+\bm{e}_\fo-\bm{e}_\of.
$$
Applying \eqref{e:Aprime.lbd} together with the trivial fact that $a\vee b\ge|a-b|$ for $a,b\ge0$ gives
\beq\label{e:diff.zo.oz}
|n_\oz-n_\zo|
\le|\intN_\oz+\tf12\extN_\oz-\intN_\zo-\tf12\extN_\zo|
	+\tf12|\extN_\oz-\extN_\zo|
\le A''+\tf12|\intY_{\ne}\setminus\trw|
= A,
\eeq
and similarly
$\tf12|\bm{e}^1_{\ne}-\bm{e}^2_{\ne}+ \numzfint- \numfzint|
\le \tf12[|\intY_{\ne}\setminus\trw| + \numzf\vee\numfz] \le A$,
so
\beq\label{e:diff.r}
|r^1-r^2-(\bm{e}_\fo-\bm{e}_\of)|\le 2A.
\eeq
Comparing \eqref{e:diff.r} with \eqref{e:unmatched.hamming} gives $|\bm{e}_\fo-\bm{e}_\of-[2s^1-l^1-(2s^2-l^2)]|\le\hamming(\outm)+4A$,
and combining with \eqref{e:unmatched.hamming} proves the lemma.
\epf
\elem

\newcommand{\projt}{\text{\textup{\textsf{proj}}}_\trw}

The following is our refinement of Propn.~\ref{p:rigidity} in the regime of small $\De$:

\bppn\label{p:improved.rigidity}
Suppose $\De\le n^{1/5}$, and
let $\tauY=(\outm,\inm)$ with $\JW$, $\IU$
as in Lem.~\ref{l:intensity}. Then
$$
\E[\idZZ^\pd_{n\al-(r^1,r^2)}[\tauY]]
\le
 n^{-\Xi}\,\E[\ZZ^\pd_{n\al-r^1}[\usi^1_\YY]]$$
for
$\Xi
\equiv \Xi(\outm,\JW,\IU) \equiv
	\tf{1}{20} \hamming(\outm)
+ 0\vee \tf{1}{20} (|\JW|-|\IU|)$.

\bpf
Recall \eqref{e:Z.boundary.decompose.A} that $\smash{\Ot^\pd_A[\pi|\tauY]}$ denotes the set of pair frozen model configurations $(\uz^1,\uz^2)$ on $\Gmin$ which are consistent with boundary conditions $\uta_\YY$, have empirical measure $\pi_\om=\intN_\om/\intN$ when restricted to $\intR$, and have $\intNE+A$ internal edges among the unequal spins. Analogously we now let $\smash{\Ot^\pd_{A,k'}[\pi|\tauY,\trw]}$ denote the set of pair frozen model configurations $(\uz^1,\uz^2)$ on $\Gmin$ which are consistent with boundary conditions $\uta_\YY$ \emph{and have the given $\trw$} (see Defn.~\ref{d:trw}), have empirical measure $\pi_\om=\intN_\om/\intN$ when restricted to $\intR$, have $\intNE+A$ internal edges among the unequal spins, and lastly have $k'$ as in \eqref{e:w.minus}. Let $\projt^1$ denote the projection mapping $(\uz^1,\uz^2)\mapsto(\uz^1,\uzprime^2)$ where $\uzprime^2$ equals $\uz^2$ on $\trw$ but equals $\uz^1$ on $V\setminus\trw$. We shall compare 
\[
\ZZt^\pd_{A,k'}[\pi|\tauY,\trw]
	\equiv
	\Big|\Ot^\pd_{A,k'}[\pi|\tauY,\trw]\Big|,\quad
\YYt^\pd[\pi|\tauY,\trw]
	\equiv
	\bigg|
	\bigcup_{A,k'}
	\projt^1\big(\Ot^\pd_{A,k'}[\pi|\tauY,\trw]\big)\bigg|.
\]
The ratio of their expected values can be upper bounded by repeating the derivation of \eqref{e:near.id.pi.A}, but now assuming $\De\le n^{1/5}$ and fixing the spins on $\trw$ as well as on $\extY\cup\intY$. Further accounting for the assignment of $A'\ge k'$ edges from $V_{\ne}$ to $\intY_\zz$ yields
$$
\RR^1_{A,k'}[\pi|\tauY,\trw]
\equiv
\f{\E[\ZZt^\pd_{A,k'}[\pi|\tauY,\trw]]}
	{\E[\YYt^\pd[\pi|\tauY,\trw]]}
\le
\f{1}{n^{k'/2}}
\f{(\log d)^{O(\De'')}d^{(\extNEucfd)''/4}}
	{d^{(\De_\zro/2+\intNE_\free/4)''}}
	\Big( \f{C (\De'')^2 d}{n \gm''} \Big)^A,
$$
where we use $''$ to denote counts off of $\trw$. Recall \eqref{e:diff.zo.oz} that $|n_\zo-n_\oz|\le A$, and note also $n^{}_\zo=n_\zo''$ and $n_\zf''\ge\intN_\of''$, so
\[\tf12\De_\zro''
	+\tf14\intNE_\free''
=\tf12(n^{}_\zo+n_\zf'')+\tf14\intNE_\free''
\ge \tf14
	[ n^{}_\zo + n^{}_\oz + n_\zf'' + \intN_\of''
	+\intNE^{}_\free
	]-\tf14 A \ge \tf14(\intNE''-A).\]
Meanwhile $(\extNEucfd)''\le\extNE''=|\intY_{\ne}\setminus\trw|\le2A$ by \eqref{e:Aprime}, so
$$
\RR^1_{A,k'}[\pi''|\tauY,\trw]
\le 
\f{(\log d)^{ O(\De'')}
	d^{ (\extNEucfd)''/4 }
	}{
		d^{\intNE''/4}
		n^{k'/2}} (d^{1/4}/n^{4/7})^A
\le \f{(d/n^{4/7})^A}{ d^{\intNE''/5} n^{k'/2}} 
\le \f{1}{ d^{\intNE''/5} n^{(A+k')/2}}.
$$
After applying $\projt^1$ on the spins of $V\setminus\trw$, the spins on $\trw$ are uniquely determined by the incoming messages $\inm$, therefore $\sum_\trw\YYt^\pd[\pi|\tauY,\trw]\le\ZZ^\pd[\pi^1|\usi^1_\YY]$. From \eqref{e:w.minus} and \eqref{e:Aprime}, $\hamming(\outm)\le2A+k'$. Combining Lem.~\ref{l:intensity} gives $|\JW|+\hamming(\outm)\le |\IU|+6(A+k')$, therefore
$$\tf12(A+k')\ge [\tf12 \hamming(\outm)]\vee [ \tf1{12} (\hamming(\outm)+ |\JW|-|\IU|) ]
\ge \tf{1}{12} \hamming(\outm)
+ \tf{1}{12} (|\JW|-|\IU|)\vee0.$$
The result follows by decomposing $\E[\idZZ^\pd_{n\al-(r^1,r^2)}[\tauY]]$ as a sum over $\E[\ZZt^\pd_{A,k'}[\pi|\tauY,\trw]]$.
\epf
\eppn

\subsection{Exponential cost of second matching}

\newcommand{\Sum}{\bm{S}}

For $i=1,2$ let $\Aw^i$ be a random matching on $\WW$, and write $\Aw^i\in\mathfrak{M}(\outm^i_\WW,j^i)$ if $\Aw^i$ forms exactly $j^i$ pairs $(w,w')$ with $\bm{o}^i_w=\bm{o}^i_{w'}=\one$. Note that $\outm_\WW$ and $\Aw$ together specify a unique
configuration $\inm(\outm_\WW,\Aw)$ of incoming messages. We therefore rewrite \eqref{e:vid} as
$$
\vid^\TT_{n\al,i}
=\sum_{\substack{r^1,j^1,\outm^1,\inm^1_\UU,\\
	\Aw^1 \in\mathfrak{M}(\outm^1_\WW,j^1)  }}
	\P(\Aw^1)\wt\ka_{r^1-j^1}(\usi^1_\UU)
\sum_{\substack{r^2,j^2,\outm^2,\inm^2_\UU,\\
	\Aw^2 \in\mathfrak{M}(\outm^2_\WW,j^2)  }}
	\P(\Aw^2)\wt\ka_{r^2-j^2}(\usi^2_\UU)
	\idZZ^\pd_{n\al-(r^1,r^2)}[ \tauY ].
$$
where $\tauY\equiv(\outm,\inm)$
with $\inm$ uniquely determined by $\inm_\UU$ and the matchings $\Aw^i$.

\newcommand{\xl}{x_\ell}
\newcommand{\xj}{x_j}

\bppn\label{p:second.matching}
It holds uniformly over all realizations $\TT$ of \eqref{e:TT} that
$$\ETT[\vid^\TT_{n\al,i}] \lesssim_{d,t} e^{-c_d\bw} \E\ZZ_{n\al}.$$

\bpf
By applying Propn.~\ref{p:improved.rigidity}
together with the fact that the total number of possibilities of
$\TT$, $r^i-j^i$, $\usi^1_\UU$, $\wt\ka$ is $\lesssim_{d,t}1$, we find
\beq
\label{e:vid.basic.ubd}
\ETT[\vid^\TT_{n\al,i}]
\lesssim_{d,t}
\sum_{\substack{
	j^1,\outm^1_\WW \\
	\Aw^1\in\mathfrak{M}(\outm^1_\WW,j^1)}}
	\P(\Aw^1)
	\E\ZZ^\WW_{n\al-j^1}[\usi^1_\WW]
\overbrace{
	\sum_{j^2,\outm^2_\WW}
	\f{
\ol\psi_{j^2}(\outm^{2}_\WW)
	}
	{
	n^{\hamming(\outm_\WW)/20
	+0\vee[|\JW|-c_{d,t}]/20}
	}
	}^{\Sum\equiv\Sum(j^1,\outm^1_\WW)}
\eeq
where $\ZZ^\WW$ refers to the partition function of the graph $G\setminus\Aw$ with $\bw\equiv|\WW|$ unmatched edges, $c_{d,t}$ denotes the maximum possible value of $|\IU|$, and $\ol\psi_{j^2}(\outm^{2}_\WW)\equiv\P(\Aw^2\in\mathfrak{M}(\outm^{2}_\WW,j^2))$. If $\outm^{2}_\WW$ equals $\one$ in $\ell$ coordinates, then
\[\ol\psi_j(\outm^{2}_\WW)
= \smb{\ell}{2j} \smf{(2j-1)!!
(\bw-\ell)_{\ell-2j}
(\bw-2\ell+2j-1)!!}{ (\bw-1)!! }
= 2^{\ell-2j} \smb{ \bw/2 }{ j,\ell-2j,\bw/2-\ell+j }
	\Big/\smb{\bw}{\ell}.\]
This is positive for $0\le\ell\le\bw$ and $0\vee(2\ell-\bw)\le 2j\le\ell$; and if $\smash{\psi^{\bw}_{\ell,j}>0}$ then crudely $\smash{\psi^{\bw}_{\ell',j'}/\psi^{\bw}_{\ell,j}\le \bw^{10[ |j'-j|+|\ell'-\ell| ] }}$ for any other $j',\ell'$. Since $\bw\lesssim_d\log n$, we see that the main effect in sum defining $\Sum\equiv\Sum(j^1,\outm^1_\WW)$ is the polynomial decay in $\hamming(\outm_\WW)$ and in $|\JW|\ge c_{d,t}$:
$$\Sum(j^1,\outm^1_\WW)
\lesssim\sum_{|j^2-j^1|\le c_{d,t}} \ol\psi_{j^2}(\outm^1_\WW).$$
Since $\Aw^1\in\mathfrak{M}(\outm^1_\WW,j^1)$,
$\usi^1_\WW$ must equal $\oo$ in exactly $2j^1$ coordinates,
and then \eqref{e:z.dangling.ones.comparison} gives
\beq
\label{e:z.w.comparison}
\E\ZZ^\WW_{n\al-j^1}[\usi^1_\WW]
\sim \wt{\bm{c}}^{-1}\, 
	\bq(\outm^1_\WW) \lm^{-j^1}
	\E\ZZ_{n\al}\eeq
for $\wt{\bm{c}}$ a proportionality constant depending on $n$, $|\WW|$, and $\lm$.  Therefore
$$
\f{\ETT[\vid^\TT_{n\al,i}]}{\E\ZZ_{n\al}}
\lesssim_{d,t}
\wt{\bm{c}}^{-1}
\sum_{\outm^1_\WW}\bq(\outm^1_\WW)
\sum_{j^1}
	\f{\ol\psi_{j^1}(\outm^1_\WW)}{ \lm^{j^1} }
\sum_{|j^1-j^2|\le c_{d,t}}
\ol\psi_{j^2}(\outm^1_\WW)
\le
\wt{\bm{c}}^{-1}
\sum_{\outm^1_\WW}\bq(\outm^1_\WW)
\sum_{j^1} \f{\ol\psi_{j^1}(\outm^1_\WW)}{ \lm^{j^1} }
\sim
1,
$$
where the last step is by a second application of \eqref{e:z.w.comparison}. Stirling's formula gives
$$
\psi^{\bw}_{\ell,j}
=\bw^{O(1)} \exp\{ \bw\, f( \ell/\bw,j/\bw ) \},\quad
1
\sim
\wt{\bm{c}}^{-1}\,\bw^{O(1)}
\sum_{( \ell/\bw,j/\bw )  \in\mathscr{C}}
\exp\{\bw\, g( \ell/\bw,j/\bw ) \}
$$
where $\mathscr{C}$ denotes the convex set $\set{0\le\xl\le1, 0\vee(2\xl-1)\le 2\xj\le \xl}$, and
\[\begin{array}{rl}
f(\xl,\xj)
	\hspace{-6pt}&\equiv
\xl\log 2-2\xj\log 2+\tf12H( 2\xj,2\xl-4\xj,1-2\xl+2\xj )-H(\xl),\\
g(\xl,\xj)
	\hspace{-6pt}&\equiv -H(\xl\,|\,\qo)
	- \xj\log\lm
	+ f(\xl,\xj).
\end{array}\]
Because of the entropy terms, the compact set $\argmax_x g(x)$ must lie some positive distance $\ep'_d$ away from the boundary of $\mathscr{C}$. By computing the Hessian, we find $g$ is strictly concave in the interior of $\mathscr{C}$, so it has a unique maximizer $(\xl^\lm,\xj^\lm)$.
Now compare
\[\xj^\lm(\xl) = \argmax_{\xj} [-\xj\log\lm+f(\xl,\xj)]
\quad\text{against}\quad
\xj^1(\xl)=\argmax_{\xj} f(\xl,\xj):\]
clearly they must be separated by some constant distance $5\ep_d$
uniformly over all $\xl$ near $\xl^\lm$.
Then, by concavity, there exists a constant $c_d>0$ such that $e^{\bw \,f(\xl,\xj)}\le e^{-2 c_d\bw}$ for all $\|(\xl,\xj)-(\xl^\lm,\xj^\lm)\|_1\le \ep_d$, while $\wt{\bm{c}}^{-1}\,e^{\bw\,g(\xl,\xj)} \le e^{-2 c_d\bw}$
for all $\|(\xl,\xj)-(\xl^\lm,\xj^\lm)\|_1\ge \ep_d$, so
$$\f{\ETT[\vid^\TT_{n\al,i}]}{\E\ZZ_{n\al}}
\lesssim_{d,t}
\bw^{O(1)}
\max_{(\xl,\xj)\in\mathscr{C}}
\wt{\bm{c}}^{-1}
\exp\{ \bw\,g(\xl,\xj)+\bw\,f(\xl,\xj) \}
\le e^{-c_d\bw}
$$
as claimed.
\epf
\eppn

\bpf[Proof of Propn.~\ref{p:var.log}]
Follows from Cor.~\ref{c:indep} and
Propn.~\ref{p:second.matching} with $m' = m - (2/c_d)\log n$.
\epf

\bibliographystyle{abbrv}
\bibliography{refs}
\end{document}